\documentclass[11pt,english]{smfart}

\usepackage[utf8]{inputenc}
\usepackage[T1]{fontenc}
\usepackage[english]{babel}

\usepackage{amsmath,amssymb,amscd,amsbsy,bm,mathrsfs,mathtools}
\usepackage{euscript}
\usepackage{faktor}
\usepackage{scalerel}
\usepackage{manfnt}
\usepackage{graphicx}
\usepackage{tikz-cd}
\usepackage{physics}
\usepackage{dsfont}
\usetikzlibrary{backgrounds,calc}

\usepackage[text={150mm,251mm},centering,marginparwidth=75pt]{geometry}
\usepackage[dvipsnames]{xcolor}
\usepackage{colortbl}
\usepackage{color}

\usepackage{enumitem}

\setlist[itemize]{
  wide=0pt,
  labelwidth=2em,
  labelsep*=0em,
  itemindent=0pt,
  leftmargin=\dimexpr\labelwidth + \labelsep\relax,
  noitemsep,
  topsep=1ex
}

\setlist[enumerate]{
  wide=0pt,
  labelwidth=2em,
  labelsep*=0em,
  itemindent=0pt,
  leftmargin=\dimexpr\labelwidth + \labelsep\relax,
  noitemsep,
  topsep=1ex,
  label=(\alph*)
}

\usepackage[
  colorlinks,
  linkcolor=BrickRed,
  citecolor=Green,
  urlcolor=Cerulean,
  hypertexnames=true
]{hyperref}

\usepackage{smfhyperref}
\usepackage[hyperpageref]{backref}

\usepackage{cleveref}
\usepackage{nameref}

\usepackage{datetime}
\usepackage{comment}
\usepackage{cite}
\usepackage[all]{xy}
\usepackage{marginnote}
\usepackage{interval}
\usepackage{letltxmacro}
\usepackage{xspace}

\allowdisplaybreaks

\makeatletter
\newcommand{\crefnames}[3]{%
  \@for\next:=#1\do{%
    \expandafter\crefname\expandafter{\next}{#2}{#3}%
  }%
}
\makeatother

\newtheorem{theorem}{Theorem}[section]
\newtheorem{lemma}[theorem]{Lemma}

\newtheorem{proposition}[theorem]{Proposition}
\newtheorem{corollary}[theorem]{Corollary}

\newtheorem{definition}[theorem]{Definition}
\newtheorem{remark}[theorem]{Remark}
\newtheorem{example}[theorem]{Example}

\newtheorem{mainthm}{Theorem}

\crefname{mainthm}{Theorem}{Theorems}
\Crefname{mainthm}{Theorem}{Theorems}

\newlist{thmlist}{enumerate}{1}

\setlist[thmlist]{
  font=\normalfont,
  label=(\roman*),
  ref=\thetheorem.(\roman{thmlisti})
}

\newlist{thmenum}{enumerate}{1}

\setlist[thmenum]{
  font=\normalfont,
  label=(\roman*),
  ref=\thetheorem.(\roman{thmenumi})
}

\newlist{corlist}{enumerate}{1}

\setlist[corlist]{
  font=\normalfont,
  label=(\roman*),
  ref=\thecorx.(\roman{corlisti})
}

\crefname{theorem}{Theorem}{Theorems}
\crefname{lemma}{Lemma}{Lemmas}
\crefname{claim}{Claim}{Claims}
\crefname{proposition}{Proposition}{Propositions}
\crefname{corollary}{Corollary}{Corollaries}
\crefname{corx}{Corollary}{Corollaries}
\crefname{conjecture}{Conjecture}{Conjectures}
\crefname{definition}{Definition}{Definitions}
\crefname{remark}{Remark}{Remarks}
\crefname{example}{Example}{Examples}
\crefname{assumption}{Assumption}{Assumptions}

\numberwithin{equation}{section}

\newcommand{\C}{\mathbb{C}}
\newcommand{\R}{\mathbb{R}}
\newcommand{\Q}{\mathbb{Q}}
\newcommand{\Z}{\mathbb{Z}}

\newcommand{\fkE}{\mathfrak{E}}
\newcommand{\fkS}{\mathfrak{S}}

\newcommand{\calq}{\mathcal{Q}}
\newcommand{\im}{\mathrm{im}}

\newcommand{\rms}{\mathrm{S}}

\newcommand{\ab}{\mathrm{ab}}
\newcommand{\id}{\mathrm{id}}

\newcommand{\bfv}{\mathbf{v}}
\newcommand{\calf}{\mathcal{F}}
\newcommand{\Trop}{\mathrm{Trop}}

\newcommand{\calo}{\mathcal{O}}

\newcommand{\bfm}{\mathbf{m}}

\newcommand{\orb}{\mathrm{orb}}
\newcommand{\N}{\mathbb{N}}

\newcommand{\Alb}{\mathrm{Alb}}
\newcommand{\interior}{\mathrm{int}}

\newcommand{\Pic}{\mathrm{Pic}}
\newcommand{\fkW}{\mathfrak{W}}
\newcommand{\fkT}{\mathfrak{T}}
\newcommand{\fkR}{\mathfrak{R}}

\newcommand{\Gr}{\mathrm{Gr}}
\newcommand{\bbS}{\mathbb{S}}
\newcommand{\red}{\mathrm{red}}

\newcommand{\BBP}{\mathbb{P}}

\begin{document}

\title[]{A Lefschetz Alternative for Real Characters on Smooth Quasi-Projective Varieties}

\author{Yuan Liu}

\email{lexliu@hku.hk}

\address{
Department of Mathematics,
The University of Hong Kong,
Pokfulam Road, Hong Kong.
}

\date{\today}

\begin{abstract}
We establish a Lefschetz alternative for real characters on smooth complex quasi-projective varieties, extending Simpson's classical theorem beyond the compact K\"ahler setting. More precisely, for a nonzero real character of the fundamental group of a smooth connected complex quasi-projective variety, either the character factors through a holomorphic fibration onto a Riemann orbifold with connected general fiber, or every level set of an associated primitive function on the corresponding covering is connected and its inclusion induces a surjection onto the fundamental group of the covering. Our approach combines logarithmic holomorphic $1$-forms with a geometric flow construction to overcome the difficulties arising from noncompactness. Meanwhile, we prove a quasi-projective analogue of the Napier--Ramachandran theorem, showing that an epimorphism onto $\Z $ whose kernel is not finitely generated must factor through a hyperbolic Riemann orbifold. As applications, we derive several new restrictions on the structure of quasi-projective groups.
\end{abstract}



\maketitle

\tableofcontents

\section{Introduction}\label{sec: introduction}

The present paper provides a generalization of Simpson's Lefschetz theorem from connected compact K\"ahler manifolds to connected smooth complex quasi-projective varieties, in the setting relevant to real characters of the fundamental group, and uses it to study quasi-projective groups, i.e., groups that can be realized as the fundamental group of a connected smooth complex quasi-projective variety. Simpson's original theorem applies to an arbitrary nonzero holomorphic $1$-form, as well as to its real part, which is a real harmonic $1$-form on a compact K\"ahler manifold. In contrast, here we start with a  nonzero real character
$$\chi:G=\pi_1(X)\to\R$$
and consider the real part of a suitably chosen logarithmic holomorphic $1$-form representing $\chi$. We say a nonzero $\chi$ is \textit{rational} if its image has rank $1$ over $\Z$; otherwise, we say it is \textit{irrational}.

We begin by recalling the following fundamental result of Simpson; see \cite[Theorem 1]{Simpson93Lef} or \cite[Theorem 2.1]{Delzant10MA}.
\begin{theorem}[Simpson]\label{thm: Simpson's Lefschetz}
Suppose $X$ is a connected compact K\"ahler manifold, and let $\omega$ be a nonzero (closed) holomorphic $1$-form. Let $p:\widehat{X}\to X$ be a regular covering such that $\Re p^*\omega$ is exact. Denote by
$$
g:\widehat{X}\to\R
$$
a primitive function of $\Re p^*\omega$. Then one of the following holds:
\begin{enumerate}[label=(\roman*)]
    \item There exists a Riemann surface $C$, a surjective holomorphic morphism
    $$
    f:X\to C
    $$
    with connected fibers, and a closed holomorphic $1$-form $\eta$ on $C$ such that
    $$
    \omega=f^*\eta.
    $$
    \item For every $t\in\R$, $g^{-1}(t)$ is connected and the inclusion induces a surjection
    $$
    \pi_1(g^{-1}(t))\twoheadrightarrow\pi_1(\widehat X).
    $$
\end{enumerate}
\end{theorem}
This theorem, originally stated for smooth projective varieties, along with other results, can be generalized to the compact K\"ahler case; see \cite[Remark~2.1.1]{Delzant10MA}. We refer to this result as \textit{Simpson's Lefschetz alternative}. It has played a central role in the study of K\"ahler groups, namely, groups that arise as fundamental groups of connected compact K\"ahler manifolds. For K\"ahler groups, the reader may refer to \cite{ArapuraBresslerRamachandran92,Arapura95survey,ABCKT,NapierRamachandran06,NapierRamachandran08,Delzant10MA,Burger11,Pybook} and so on. Remarkably, even its restriction to rational characters has significant geometric consequences. A closely related result was established independently by Napier and Ramachandran through a quite different approach, as part of their work on proper holomorphic fibrations from complete noncompact K\"ahler manifolds to Riemann surfaces; see \cite{NapierRamachandran95,Ramachandran96,NapierRamachandran97,NapierRamachandran01,NapierRamachandran08}.

\begin{theorem}\textup{\textbf{(Napier--Ramachandran, \cite[Theorem 4.3]{NapierRamachandran01})}}\label{thm: NR on Z covering}
Suppose $X$ is a connected compact K\"ahler manifold. Assume that there exists a surjective homomorphism
$$
\chi:\pi_1(X)\twoheadrightarrow\Z
$$
such that its kernel is not finitely generated. Then $X$ admits a surjective holomorphic mapping onto a Riemann surface $C$, and the homomorphism factors through $\pi_1(C)$.
\end{theorem}
This result follows from Simpson's Lefschetz alternative. Indeed, for a rational character $\chi$, choose a holomorphic $1$-form $\omega$ whose real part $\alpha=\Re\omega$ represents
$\chi$. Then the primitive function of $\alpha$ on the corresponding $\Z$-covering is proper, so every regular level set is compact and smooth. If the second alternative of \cref{thm: Simpson's Lefschetz} held, the fundamental group of such a level set would be finitely generated and would surject onto $\ker\chi$, forcing $\ker\chi$ to be finitely generated. Thus the first alternative must occur. By contrast, the essential geometric argument underlying \cref{thm: NR on Z covering} accounts for much of the technical difficulty in proving the alternative \cref{thm: Simpson's Lefschetz} itself, as is evident from \cref{sec: Lefschetz alternative rational case,sec: Lefschetz alternative general case}.

This naturally raises the following\\
\textit{Question:} To what extent does Simpson's Lefschetz theorem \cite{Simpson93Lef} and the thoerem of Napier--Ramachandran extend to smooth complex quasi-projective varieties?

The first question has been studied from different perspectives. Zhou \cite{Zhou} proved a Lefschetz-type theorem for weakly $1$-complete complex manifolds under an additional properness assumption, while Corlette and Simpson \cite[Theorem~5.11]{CorletteSimpson08} treated the quasi-projective case when the relevant holomorphic $1$-form extends holomorphically across a smooth projective compactification.

Here we start with a nonzero real character $\chi:G\to\R$ and consider a specific logarithmic holomorphic $1$-form whose real part represents $\chi$. In this setting, we show that Simpson's Lefschetz alternative continues to hold for smooth complex quasi-projective varieties. In full generality, however, one should perhaps not expect the same statement to hold for an arbitrary holomorphic $1$-form or real harmonic $1$-form. For example, the holomorphic $1$-form $\sin(1/z)\,dz$ on $\mathbb C^*$ has zeros accumulating at the boundary; see \cref{rmk: accumulating zero sets}.

\begin{figure}[htbp]
\centering
\small
$$
\begin{tikzcd}[
  row sep={2cm,between origins}, % Vertical spacing
  column sep=1.5em,             % Default horizontal spacing
  cells={nodes={font=\small}},
  execute at end picture={
    \coordinate (upper) at ($(caseI)!0.5!(both)$);
    \coordinate (lower) at ($(both)!0.5!(caseII)$);
    % Ellipse radii: horizontal 1.65cm, vertical 1.7cm.
    % Change all occurrences together.
    \begin{scope}[on background layer]
      \fill[gray!4] (upper) ellipse (1.65cm and 1.7cm);
      \fill[gray!4] (lower) ellipse (1.65cm and 1.7cm);
      \begin{scope}
        \clip (upper) ellipse (1.65cm and 1.7cm);
        \fill[gray!20] (lower) ellipse (1.65cm and 1.7cm);
      \end{scope}
      \draw[gray!65] (upper) ellipse (1.65cm and 1.7cm);
      \draw[gray!65] (lower) ellipse (1.65cm and 1.7cm);
    \end{scope}
  }
]
\dim\psi(X)=1
  &[6mm] % Increase this value to lengthen BOTH LEFT arrows.
    |[alias=caseI]| \textup{(i)}
    \arrow[l,Leftrightarrow]
  &[6mm] % Increase this value to lengthen the OVERLAP arrow.
    \textup{(i)}\setminus\textup{(ii)}
    \arrow[l,Rightarrow]
    \arrow[r,Leftrightarrow]
  & \textup{(i)}_{\mathrm{hyp}}
    \arrow[r,Leftrightarrow]
  & \ker\chi\text{ not f.g.}
\\
  & |[alias=both]| \boxed{\textup{(i)}\cap\textup{(ii)}}
    \arrow[u,Rightarrow]
    \arrow[d,Rightarrow]
    \arrow[r,Leftrightarrow]
  & \textup{(i)}_{\mathrm{non-hyp}}
    \arrow[rr,Leftrightarrow]
  &
  % Middle-right statement: more negative xshift moves it left.
  & |[overlay,xshift=-8mm]| \begin{gathered}
      \dim\psi(X)=1,\\[-2pt]
      \ker\chi\text{ f.g.}
    \end{gathered}
\\
\begin{gathered}
  \dim\psi(X)\ge2\\[-2pt]
  \Updownarrow\\[-2pt]
  \textup{(ii)}\setminus\textup{(i)}
\end{gathered}
    \arrow[r,Rightarrow]
  & |[alias=caseII]| \textup{(ii)}
    \arrow[rrr,Leftrightarrow]
  &
  &
  % Lower-right statement: more negative xshift moves it left.
  & |[overlay,xshift=-25mm]| \ker\chi\text{ f.g.}
\end{tikzcd}
$$
\caption{logical relations for rational Lefschetz alternative: here (ii) refers to the specific $\alpha$ constructed in \cref{subsec: Choice of logarithmic holomorphic one-form}.}
\label{fig: logical relations for Lefschetz alternative}
\end{figure}
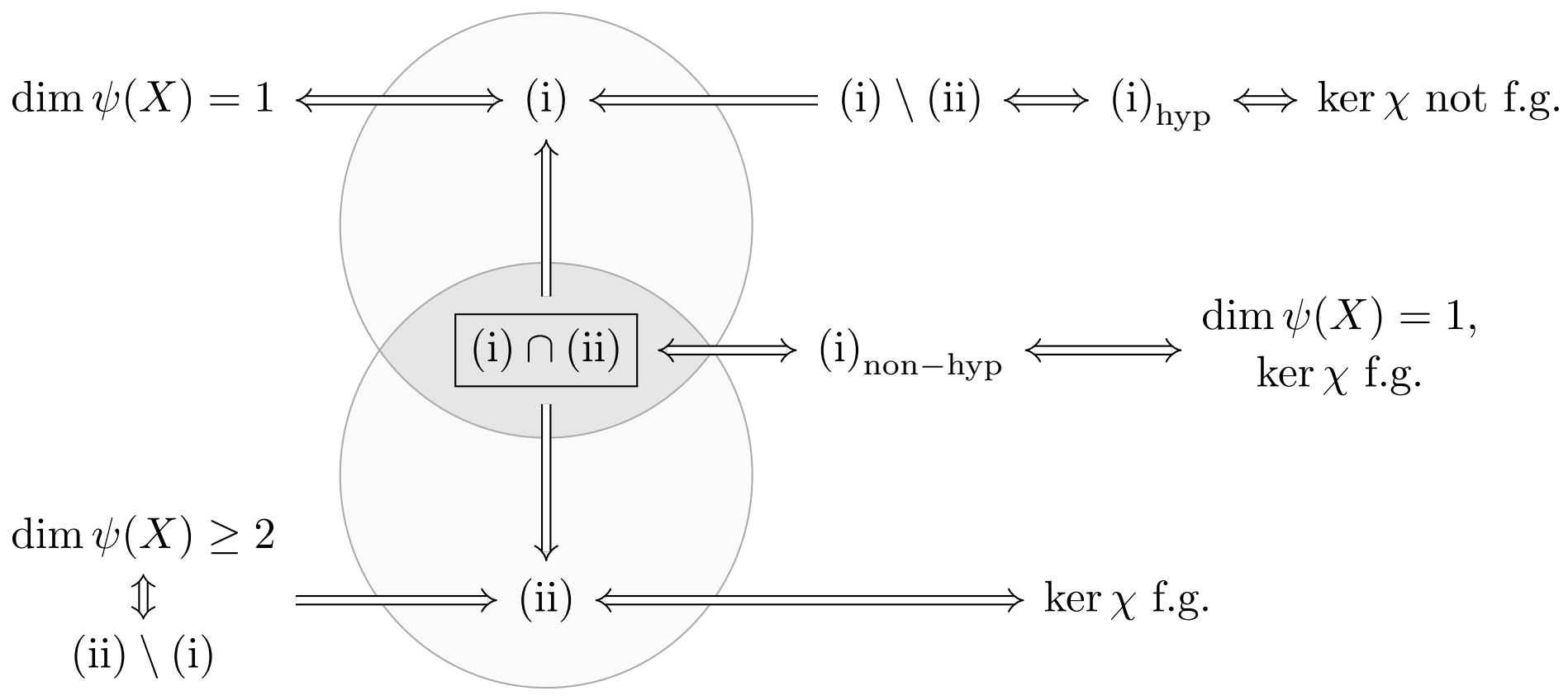

\begin{mainthm}[Lefschetz alternative, see \cref{thm: Lefschetz alternative general case} for this general version and \cref{thm: Lefschetz alternative rational case} for the rational version]\label{mainthm: lefschetz alternative}
Let $X$ be a smooth connected complex quasi-projective variety and let $\chi:\pi_1(X)\to\R$
be a nonzero group homomorphism. Let $p:\widehat X\to X$ be the corresponding $\chi$-covering, with $\pi_1(\widehat X)\cong\ker\chi$. Then one of the following holds:
\begin{enumerate}[label=(\roman*)]
\item There exists a surjective algebraic morphism with connected general fiber $f:X\to C$ onto a smooth quasi-projective curve such that, after equipping $C$ with the orbifold structure determined by the gcd multiplicities $\bfm$ of the fibers of $f$, there exists a homomorphism
$$\chi_C^\orb:\pi_1^\orb(C,\bfm)\to\R$$
satisfying
$$\chi=\chi_C^\orb\circ f_*^\orb.$$
\item There exists a smooth closed real $1$-form $\alpha$ on $X$ representing $\chi$ such that, for any primitive function
$$u:\widehat X\to\R,\quad du=p^*\alpha,$$
every level set $L_t=u^{-1}(t)$ is connected and the inclusion induces a surjection
$$\pi_1(L_t)\twoheadrightarrow\pi_1(\widehat X).$$
Equivalently, the pair $(\widehat X,L_t)$ is $1$-connected for every $t\in\R$.
\end{enumerate}
\end{mainthm}
\noindent We call it a \textit{Lefschetz alternative} since it is not a dichotomy; see \cref{rmk: overlap for Lefschetz alternative} and \cref{fig: logical relations for Lefschetz alternative} for the rational case. Here in (ii), the smooth closed 1-form $\alpha$ is the real part of a logarithmic holomorphic 1-form on $X$ with nonzero purely imaginary residues at infinity, constructed in \cref{subsec: Choice of logarithmic holomorphic one-form,subsec: Modifications before lifting} for rational and general case respectively; more specifically, see \cref{eq: log omega,eq: b pull-back of omega,eq: alpha real case}.

Correspondingly, we have the following analogue of \cref{thm: NR on Z covering} in the smooth quasi-projective setting:

\begin{mainthm}[\cref{thm: quasi-proj ker not fg implies fact}]\label{mainthm: quasi-proj ker not fg implies fact}
Let $X$ be a smooth connected complex quasi-projective variety and let $\chi:\pi_1(X)\twoheadrightarrow\Z$
be an epimorphism. If $\ker\chi$ is not finitely generated, then there exists a surjective algebraic morphism with connected general fiber $f:X\to C$
onto a smooth quasi-projective curve such that, after equipping $C$ with the orbifold structure determined by the gcd multiplicities $\mathbf m$ of the fibers of $f$, $(C,\mathbf m)$ is hyperbolic and there exists an epimorphism
$$\chi_C^\orb:\pi_1^\orb(C,\mathbf m)\twoheadrightarrow\Z$$
such that
$$\chi=\chi_C^\orb\circ f_*^\orb.$$
\end{mainthm}

As in the compact case, the geometric argument underlying this factorization accounts for much of the work required for the quasi-projective extension of \cref{thm: Simpson's Lefschetz}, say \cref{mainthm: lefschetz alternative}. The proof of rational case combines ideas of Simpson \cite{Simpson93Lef}, especially the local picture developed there, with ideas of Zhou \cite{Zhou}, in particular the construction of a flow used to push a path in a prescribed direction. The main point is to ensure that the deformed path remains in the same relative homotopy class as the original one, while being pushed in the direction determined by the primitive function $u$ on the regular $\Z$-covering associated with $\chi$. We carry out this argument in \cref{sec: Lefschetz alternative rational case} and then adapt it to the general $\chi$-cover in \cref{sec: Lefschetz alternative general case}. An outline of the proof strategy appears at the beginning of \cref{sec: Lefschetz alternative rational case}.

A cautious reader may notice that \cref{mainthm: quasi-proj ker not fg implies fact} involves a \textit{hyperbolic} holomorphic fibration, whereas \cref{mainthm: lefschetz alternative} refers only to a holomorphic fibration. This distinction is one of the main reason why the two alternatives in \cref{mainthm: lefschetz alternative} need not be mutually exclusive, see \cref{fig: logical relations for Lefschetz alternative}. We explain this phenomenon in terms of the BNS invariant in \cref{subsec: Lefschetz alternative revisited via the BNS invariant}.

Given the central role of Simpson's Lefschetz theorem in the study of K\"ahler groups, one may expect its extension to the smooth quasi-projective setting to have a number of consequences for quasi-projective groups. Here we mention only those consequences that are most closely related to the previous work of the author and his collaborator \cite{LiuLiu25,Liu26CR}. First, we describe the complement of the BNS invariant for quasi-projective groups, in analogy with Delzant's result for K\"ahler groups \cite{Delzant10MA}.

\begin{mainthm}[\cref{thm: BNS quasiprojective}]
\label{thm: BNS quasiprojective introduction}
Let $X$ be a connected smooth complex quasi-projective variety with $G=\pi_1(X)$. Then
$$\Sigma^1(G)^\complement=
\rms\left(\bigcup_f\im\left(f^*:H^1(C,\mathbf m;\R)\to H^1(X;\R)\right)\right),$$
where the union runs over all hyperbolic orbifold fibrations $f:X\to(C,\mathbf m)$, which is a finite union. Here $\rms$ means the projection onto the unit sphere in $H^1(X;\R)$.
\end{mainthm}

For restrictions on quasi-projective groups, we have
\begin{proposition}
\begin{enumerate}
    \item (\cref{thm: HNN}) Any group admitting a properly ascending (or descending) HNN decomposition with finitely generated base group is not a quasi-projective group.
    \item (\cref{prop: Thompson not quasi-proj}) Any Thompson--Stein group $F\left(\ell,\Z\left[\frac{1}{n_1\cdots n_k}\right],\langle n_1,\ldots,n_k\rangle\right)$ is not quasi-projective.
\end{enumerate}
\end{proposition}
We also remark that \textup{(a)} and \textup{(b)} with $k=1$ require only the rational Lefschetz alternative.

\begin{remark}
A Lefschetz-type theorem for smooth quasi-projective varieties was also claimed by Rogov in \cite[Appendix A]{Rogov}. The main issue we see in the proof occurs in the local trivialization of \cite[Lemma~A.13]{Rogov}. There, the strict monotonicity of $\Re g$ along a complete vector field is used to conclude that every trajectory meets the inverse image of a prescribed real line in $\C$ exactly once. In the noncompact setting, however, completeness and strict monotonicity alone do not imply such a global crossing property: a complete trajectory may escape to infinity while $\Re g$ remains bounded. Thus an additional uniform estimate forcing the relevant trajectories to cross the prescribed hypersurface would be needed, but such an estimate is not established. Consequently, the local trivialization in \cite[Lemma~A.13]{Rogov} is not justified by the stated argument.

%This issue is closely related to the difficulty of extending Simpson's global picture and patching procedure to the quasi-projective setting; see \cref{rmk: on patching difficulty}. Our proof avoids this global difficulty by using the flow only on compact flow neighborhoods of individual paths and treating the critical and boundary regions locally.

%We also note a separate issue in the moving-surface argument of \cite[Lemma~A.5]{Rogov}, where the required equidimensional surface is assumed to be choosable through arbitrary prescribed points; this may fail when a prescribed point lies on a positive-dimensional fiber while $\omega$ is nonzero there.

Our proof differs substantially from his, both conceptually and technically. In particular, we replace the patching argument and global construction in \cite{Simpson93Lef} with a modification of the geometric flow constructed in \cite{Zhou}.
%Our proof is substantially different from his, both conceptually and technically. For example, we do not use the pathching process and global picture in \cite{Simpson93Lef}; instead, we modify the geometric flow constructed in \cite{Zhou} instead.
\end{remark}

We also emphasize in this paper the relation between the Lefschetz-type theorem and the theorem of Napier--Ramachandran, which is one of the reasons for treating the rational case first. As the proof shows, although the lifting picture in the rational case is easier to describe, it already contains most of the essential difficulties of the argument.

There is also a second motivation coming from cohomology jump loci. As shown by Delzant in \cite{Delzant08GAFA}, the rational part of the complement of the BNS invariant alone is sufficient to obtain a geometric description of the rank-one cohomology jump loci over any algebraically closed field of K\"ahler groups, while this rational part is closely related to the rational case of Simpson's Lefschetz theorem. These considerations motivate our emphasis on the rational case.

In an upcoming paper with Yongqiang Liu \cite{LiuLiuTBD}, we use the rational Lefschetz alternative, \cref{thm: Lefschetz alternative rational case}, to study twisted rank-one cohomology jump loci (see \cite{LiuSuciu26}) over any algebraically closed field in the smooth quasi-projective setting. Our approach there follows a path different from that of \cite{Delzant10MA}: rather than using the BNS invariant, we exploit the tropical geometry of cohomology jump loci, continuing the approach developed in \cite{LiuLiu25}.

The paper is organized as follows. In \cref{sec: Riemann orbifolds and fibration}, we introduce the basic notation concerning hyperbolic Riemann orbifolds and holomorphic fibrations. In \cref{sec: Lefschetz alternative rational case}, we prove the rational version of the Lefschetz alternative in the smooth quasi-projective setting; in \cref{sec: Lefschetz alternative general case}, we prove \cref{mainthm: lefschetz alternative} for arbitrary nonzero real characters and explain the essential modifications needed to pass from the rational to the general case. Finally, in \cref{sec: applications}, we discuss applications of \cref{mainthm: lefschetz alternative} to quasi-projective groups and, in \cref{subsec: Lefschetz alternative revisited via the BNS invariant}, compare \cref{mainthm: lefschetz alternative} with \cref{mainthm: quasi-proj ker not fg implies fact} from the viewpoint of the BNS invariant.
\section{Preliminaries on Riemann orbifolds and fibration}\label{sec: Riemann orbifolds and fibration}
In this section, we recall the notion of a Riemann
orbifold, namely a Riemann surface equipped with
an orbifold structure allowing finitely many points
with nontrivial finite isotropy groups. For proofs
of the results presented below, we refer the reader
to \cite[\S 4]{Delzant08GAFA}, \cite[\S 5]{NapierRamachandran08} and the references
therein.

Let $X$ be a connected complex manifold, and let $C$ be a Riemann surface. Let $h:X\to C$ be a surjective holomorphic map with connected general fibers, which we call a \textit{holomorphic fibration}. For a subset $A\subset C$, denote $X_A\coloneqq h^{-1}(A)$; in particular, for a point $z\in C$, denote $X_z=h^{-1}(z)$. Assume that the set of critical values of $h$ is locally finite. For each critical value $p$, let $m_p$ be the greatest common divisor of the multiplicities of the irreducible components of the divisor $h^{-1}(p)$. Let
$$P\coloneqq\{p\in C:m_p\ge2\}.$$
For each $p\in P$, let $\gamma_p$ be a positively oriented loop tracing the boundary circle of a coordinate disk $D_p$ in $C$ centered at $p$, chosen so that $\overline{D_p}\cap P=\{p\}$. Then the orbifold fundamental group is given by
$$
\pi_1^\orb(C,\bfm)\coloneqq \pi_1(C\setminus P)/N
$$
where $N$ is the normal subgroup of $\pi_1(C\setminus P)$ generated by the loops $\{\gamma_p^{m_p}\}$ with $m_p\ge 2$, and $\bfm\subset {P}^{\N}$ is the data of critical points with multiplicities. For a regular value $z\in C\setminus P$ of $h$, we have
\begin{equation}\label{eq: ses of rimann fibration}
\pi_1(X_z)\to \pi_1(X)\xrightarrow{h_*\coloneqq h_*^\orb} \pi_1^\orb(C)\to 1
\end{equation}
When $h$ is proper or algebraic, we have the kernel of $\pi$ is finitely generated (see, for example, \cite{Simpson89}).

For a compact orbifold $(C,\mathbf m)$, defines it orbifold Euler characteristic as:
$$\chi^{\mathrm{orb}}(C,\bfm)
=2-2g(C)-\sum_{p\in P}
\left(1-\frac1{m_p}\right)$$
where $g(C)$ is the genus of its smooth compactification. If $\chi^\orb(C,\bfm)<0$, we say it is \textit{hyperbolic}.

When $X$ is compact, so is $C$, then $P$ is a finite set of points. The following finiteness of holomorphic fibrations for $X$ fixed is well known, for example, see \cite[Theorem 2]{Delzant08GAFA}. Note that a holomorphic map $h: X\to (C,\bfm)$ means a holomrphic map such that, for each orbifold point $p_i$, the multiplicity of every irreducible component of $h^{-1}(p_i)$ is divisible by $m_i$. Two holomorphic fibrations $h_i:X\to\Sigma_i, i=1,2$ are considered equivalent if there exists an orbifold isomorphism $\psi:\Sigma_1\to\Sigma_2$ such that $h_2=\psi\circ h_1$.
\begin{theorem}\label{thm: finiteness of hyperbolic fibration when compact}
Let $X$ be a compact complex manifold, then there exists, up to equivalence, a finite set of pairs $(C_i,h_i)$, where $C_i=(C_i,\bfm_i)$ is a hyperbolic Riemann orbifold, and $h_i: X\to (C_i,\bfm)$ is a holomorphic map with connected fibers.
\end{theorem}
In the compact case, by applying the Stein factorization theorem and replacing the target by the corresponding intermediate Riemann surface if necessary, we may assume that every fiber of a holomorphic fibration is connected. Similarly, in the noncompact case, if $h:X\to C$ is proper, the same argument shows that, after replacing $C$ by the corresponding intermediate Riemann surface, we may assume that every fiber is connected. A related factorization result for smooth fiber bundle without the properness assumption, under the additional assumptions that $h$ is a holomorphic submersion and that every fiber has finitely many connected components, is given in \cite[Theorem~2]{Zhou}; Zhou refers to this as a noncompact Stein factorization.

However, a holomorphic surjection from a noncompact quasi-projective manifold to an open Riemann surface need not be proper, and the hypotheses of \cite[Theorem~2]{Zhou} need not be satisfied. Thus, in the general quasi-projective setting, one can in general require only that the general fiber be connected. From now on, by a holomorphic fibration from a compact complex manifold we mean a surjective holomorphic map with connected fibers, whereas in the noncompact case we require only that the general fiber be connected.

For the case when $X$ and $C$ are noncompact, and the Riemann surface $C$ is a compact Riemann surface with finitely many points removed, say
$$
(C\setminus\{p_1,\cdots,p_r\},\bfm)
$$
with $\bfm$ the data of orbifold points in $C\setminus \{p_1,\cdots,p_r\}$; we denote this data as a triple $(C,r,\bfm)$. The Euler characteristic of such an open Riemann surface is
\begin{equation}\label{eq: char of orbifold}
\chi^\orb(C,r,\bfm)=2-2g(C)-r-\sum_{p\in P}
\left(1-\frac1{m_p}\right).
\end{equation}
If $\chi^\orb(C,r,\bfm)<0$, we say it is \textit{hyperbolic}.

For our applications, the real part of a holomorphic $1$-form on $C$ must represent a nonzero class in $H^1(C;\R)$. We therefore assume that its dimensions $b_1(C)\ge 1$. If $C$ is compact, this is equivalent to $g\coloneqq g(C)\ge 1$; if $C$ is noncompact, it is equivalent to $2g+r-1\ge 1$. We also need to pay attention to the case when $\chi^\orb=0$ when $2g+r\ge 2$, namely
\begin{equation}\label{eq: euler char is zero}
2g+r=2,\quad \text{ and we have no multiply fibers.}
\end{equation}
We have the following two cases:
\begin{itemize}
    \item $g=1,r=0$: then $C$ is an elliptic curve, and $\pi_1^\orb(C)=\pi_1(C)\cong \Z^2$;
    \item $g=0,r=2$: then $C\cong \C^*$, and $\pi_1^\orb(C)=\pi_1(C)\cong \Z$.
\end{itemize}
In either case, the (orbifold) fundamental group is free abelian.

Under the condition $2g+r\ge 2$, the finiteness of hyperbolic fibrations in the noncompact case follows from \cite[Theorem~1.6]{Arapura97}; see also \cite[\S 4.1]{DuttaHaoLiu24}. We state the following results which is needed in \cref{thm: BNS quasiprojective}.
\begin{theorem}\label{thm: finiteness of hyperbolic fibration when noncompact}
Let $X$ be a smooth connected complex quasi-projective variety. Then, up to equivalence, there exist only finitely many pairs $(C_i,h_i)$, where $C_i^\orb\coloneqq (C_i,r_i,\mathbf m_i)$ is a hyperbolic Riemann orbifold with $2g(C_i)+r_i\ge 2$, and $h_i:X\to C_i$ is a surjective algebraic morphism with general fiber connected .
\end{theorem}

We finish this section by a result on the BNS-invariant of hyperbolic Riemann orbifolds (See \cref{subsec: Applications on BNS invariants} for essential definitions).
\begin{proposition}\label{prop: bns empty for hyperbolic RS}
Let $(C,\bfm)$ be a hyperbolic Riemann orbifold, open or closed, then we have $\Sigma^1(C,\bfm)=\varnothing$.
\end{proposition}
\begin{proof}
This result is well known. For the compact case, for example, we can use \cite[Proposition~5.5]{LiuLiu25}; for the noncompact case, a similar calculation on rank one jump loci is carried out in \cite{LiLiu21}, and the same tropical upper bound of $\Sigma^1(G)$ (see \cite[Theorem~1.1]{LiuLiu25} and \cref{thm: BNS quasiprojective}) gives the same result.

The emptiness of $\Sigma^1(G)$ for $(C,\bfm)$ noncompact can also be seen directly as follows. For $(C,r,\bfm)$ a noncompact Riemann orbifold, it is well-known that
$$
\pi_1^{\orb}(C,r,\mathbf m) \cong F_{2g+r-1}* \Z/m_1\Z\cdots *\Z/m_s\Z,
$$
We only need to consider the case in which the character sphere is nonempty, namely $2g+r-1\ge1$. If $2g+r\ge 3$, then the nonabelian free group $F_{2g+r-1}$ has empty BNS invariant, and so does $\pi_1^\orb(C,r,\bfm)$ by \cref{prop: of bns}(a); if $2g+r=2$, then there exists multiply fibers, and $\pi_1^\orb(C,r,\bfm)\cong \Z*\Z/m_1\Z\cdots *\Z/m_s\Z$, whose BNS-invariant is empty by direct computation.
\end{proof}

\section{Lefschetz alternative of smooth quasi-projective manifolds in the rational case}\label{sec: Lefschetz alternative rational case}

In this section, we prove a version of the Lefschetz alternative for connected smooth complex quasi-projective varieties under the assumption that the real character is \textit{rational}. We begin with this special case because of its relative clarity, and because its proof illustrates several of the main ideas that will also enter the proof of the general case.

%Without the rationality assumption, this result was proved by Simpson in \cite{Simpson93Lef} for smooth projective manifolds, as well as in a quasi-projective setting when the boundary divisor is contained in the zero locus of the one-form. There is also a parallel statement formulated directly for the holomorphic one-form itself, rather than for its real part. Simpson's argument extends with only minor modifications to compact K\"ahler manifolds; see \cite[Remark 2.1.1]{Delzant10MA}. For noncompact weakly $1$-complete complex manifolds, a corresponding result was proved by Zhou in \cite[Theorem 1]{Zhou} under an additional assumption of properness of the primitive function.

Let $X$ be a smooth quasi-projective manifold, and let $G=\pi_1(X)$. Up to rescaling, let $\chi:G\twoheadrightarrow \Z $ be a nonzero rational character and set $N=\ker\chi$. Let
$p:\widehat X\to X$ be the corresponding regular $\Z $-cover, so that $\pi_1(\widehat X)=N$ and the deck transformation group is naturally identified with $\Z $. Denote by $T$ the deck transformation corresponding to $1\in\Z $.

Choose a smooth closed $1$-form $\theta$ on $X$ representing $\chi$, in the sense that
$$\int_\gamma\theta=\chi(\gamma)$$
for every $\gamma\in\pi_1(X)$. Since $p^*\theta$ has trivial periods on $\widehat X$, it is exact. Hence there exists a smooth function $u:\widehat X\to\mathbb R$ such that $du=p^*\theta$. For a fixed choice of $\theta$, the function $u$ is unique up to an additive constant. Moreover,
$$u(Tx)=u(x)+1$$
for every $x\in\widehat X$. We call $u$ a \textit{primitive function} associated with $\chi$ on the $\Z $-cover.

\begin{theorem}\label{thm: Lefschetz alternative rational case}
Let $X$ be a smooth connected complex quasi-projective variety and let $\chi:\pi_1(X)\twoheadrightarrow \Z $ be an epimorphism. Let $p:\widehat X\to X$ be the corresponding $\Z$-covering. Then one of the following holds:
\begin{enumerate}[label=(\roman*)]
    \item There exists a surjective algebraic morphism with connected general fiber $f:X\to  C$ onto a smooth quasi-projective curve such that, after equipping $C$ with the orbifold structure determined by the gcd multiplicities of the fibers of $f$, there exists a homomorphism
    $$\chi_C^\orb:\pi_1^\orb(C,\bfm)\to \Z $$
    satisfying
    $$\chi=\chi_C^\orb\circ f_*^\orb.$$
    \item There exists a smooth closed $1$-form $\alpha$ on $X$ representing $\chi$ such that, for a primitive function
    $$u:\widehat X\to\R,\quad du=p^*\alpha,$$
    every level set $L_t=u^{-1}(t)$ is connected and the inclusion induces a surjection
    $$\pi_1(L_t)\twoheadrightarrow\pi_1(\widehat X).$$
    Equivalently, the pair $(\widehat X,L_t)$ is $1$-connected for every $t\in\R$.
\end{enumerate}
\end{theorem}

\begin{remark}\label{rmk: overlap for Lefschetz alternative}
The two alternatives are not mutually exclusive. For example, take $X=\C^*$ and let $\chi:\pi_1(X)\cong\Z\to\Z$ be the identity. Then both \textup{(i)} and \textup{(ii)} hold. For an example of compact case, take $X$ to be an elliptic curve, with $\chi:\pi_1(X)\to \Z$ any nonzero character.
\end{remark}

%If we further assume that $\ker(\chi)$ is not finitely generated, then a similar argument as in the introduction can be used to show that \textup{(i)} holds. 

%Indeed, if \textup{(i)} failed, our construction on a Lefschetz surface $Y\subset X$, with $\pi_1(Y)\cong\pi_1(X)$, would yield a logarithmic representative whose primitive extends to a proper function $u:\widehat M\to\R$, where $M$ is the compact model constructed in \cref{subsec: Reduction to surfaces}. For a regular value $t$, the compactified level set $\overline L_t$ is compact. Since $L_t=\overline L_t\cap\widehat Y$ is obtained by removing finitely many real-analytic subsets, it has the homotopy type of a finite CW complex. In particular, $\pi_1(L_t)$ is finitely generated. The surjection supplied by our construction, $$\pi_1(L_t)\twoheadrightarrow\pi_1(\widehat Y)\cong\ker(\chi),$$ would then imply that $\ker(\chi)$ is finitely generated, a contradiction. Hence \textup{(i)} holds. Moreover, the orbicurve $C$ is hyperbolic, as shown in the proof of \cref{thm: quasi-proj ker not fg implies fact} at the end of \cref{subsec: Construction of a flow on the infinite cyclic covering}. This provides a generalization of \cref{thm: NR on Z covering} to the smooth quasi-projective setting.

\begin{theorem}\label{thm: quasi-proj ker not fg implies fact}
Let $X$ be a smooth connected complex quasi-projective variety and $\chi:\pi_1(X)\twoheadrightarrow\Z$ an epimorphism. If $\ker(\chi)$ is not finitely generated, there is a surjective algebraic morphism with connected general fiber $f:X\to  C$ onto a smooth quasi-projective curve such that, after equipping $C$ with the gcd multiplicities $\mathbf m$ of the fibers of $f$, $(C,\mathbf m)$ is hyperbolic and
$$\chi=\chi_C\circ f_*^\orb,\quad\chi_C:\pi_1^\orb(C,\mathbf m)\twoheadrightarrow\Z.$$
\end{theorem}

\begin{remark}
If in \cref{thm: Lefschetz alternative rational case}(i), we additionally require the orbifold curve $(C,\mathbf m)$ to be hyperbolic, then we have
$$(i) \Longleftrightarrow \ker(\chi) \text{ is not finitely generated.}$$
This will be explained using the BNS invariant in \cref{subsec: Lefschetz alternative revisited via the BNS invariant}.
\end{remark}

We further transfer this non-finitely generated property into a topological condition.

\begin{proposition}\label{prop: non fg implies level set disconnect}
Let $M$ be a compact connected smooth manifold, and let
$$\chi: G\coloneqq \pi_1(M)\twoheadrightarrow \Z$$
be a group epimorphism. Denote the coverings 
$$\widetilde{M}\xrightarrow{q} \widehat{M}\xrightarrow{p} M$$
where $\widehat{M}$ is the regular $\Z$-covering and $\widetilde{M}$ is the universal covering. Take $N=\ker(\chi)$, and $u$ be the primitive function of $\chi$ on $\widehat{M}$. If $N$ is infinitely generated, then, for every regular value $t$ of $u$, the level set of $q^*u$
$$\widetilde{M}_t\coloneqq \{x\in \widetilde{M}: u\circ q(x)=t\}$$
has infinitely many connected components.
\end{proposition}
\begin{proof}
By Sard theorem, we can find a regular value $t$ of $u$, such that the conclusion of this lemma is not vacuous. Due to $M$ is compact and $u(T.x)=u(x)+1$, the function $u:\widehat{M}\to \R$ is proper. Thus the level set $\widehat{M}_t\coloneqq\{x\in \widehat{M}:u(x)=t\}$ is compact, and we decompose it into the finite union of connected compact components:
$$\widehat{M}_t=M_1\sqcup M_2\sqcup \cdots \sqcup M_k.$$
Fix a component $M_i$, and write $G_i=\im[\pi_1(M_i)\to \pi_1(M)]$. Due to $M_i$ is compact and smooth, $G_i$ is finitely generated. If $[N: G_i]<\infty$, then $N$ would be finitely generated, contradicting the assumption. Hence we have $[N:G_i]=\infty$, which equals the number of connected components of $q^{-1}(M_i)\subset \widetilde{M}_t$. This shows that $\widetilde{M}_t$ has infinitely many connected components.
\end{proof}

The proof of \cref{thm: Lefschetz alternative rational case} for \textit{regular} level sets, together with the proof of \cref{thm: quasi-proj ker not fg implies fact}, will be given in \cref{subsec: On semi-Albanese map determined by a character,subsec: Choice of logarithmic holomorphic one-form,subsec: Reduction to surfaces,subsec: Construction of two types of critical neighborhoods,subsec: Construction of a flow on the infinite cyclic covering}. The extension of \cref{thm: Lefschetz alternative rational case} to arbitrary level sets will be given in \cref{subsec: From regular level set to any level set}.
The main ideas of the proof are as follows:
\begin{enumerate}[label=(\arabic*)]
\item First, the character $\chi$ factors through the quasi-Albanese morphism. We construct an associated morphism
$$\psi:X\to A$$
to a semi-abelian variety $A$ such that the induced character does not vanish on any positive-dimensional connected algebraic subgroup of $A$. In the rational case, the torus part of $A$ has dimension at most one. We then take a natural projective compactification of $A$ whose boundary consists of two disjoint divisors, and use Deligne's logarithmic de Rham theory to choose a logarithmic $1$-form $\omega$ whose real part representing $\chi$, with nonzero purely imaginary residues along the boundary. This is carried out in \cref{subsec: On semi-Albanese map determined by a character,subsec: Choice of logarithmic holomorphic one-form}.
\item When $\dim\psi(X)\ge2$, the Lefschetz hyperplane theorem allows us to reduce to a smooth quasi-projective surface $Y\subset X$ such that $\pi_1(Y)\cong\pi_1(X)$ and $\psi|_Y$ is generically finite onto its image. After resolving indeterminacies and the boundary, we choose smooth projective compactifications $Y\subset\bar Y$, $A\subset\bar A$ such that $\psi|_Y$ extends to a morphism $\bar\psi:\bar Y\to\bar A$, and the relevant boundary divisors have simple normal crossings. We define the polar boundary by $P=\bigl(\bar\psi^{-1}(D_A)\bigr)_{\red}$, and let $Q$ be the union of the remaining irreducible components of $\bar Y\setminus Y$. Thus $\bar Y\setminus Y=P\cup Q$, and denote $Y^\circ=\bar Y\setminus P$. The pullback $\xi=\bar\psi^*\omega$ is holomorphic on $Y^\circ$, and the nonzero residues of $\omega$ give uniform control near $P$. We implement this control by passing to the real-oriented blowup $M=\operatorname{Bl}^{\mathrm{ro}}_P(\bar Y)$,
on which the angular variables remain well defined along the boundary; another way to see this is to use a Poincar\'e-type metric. This is carried out in \cref{subsec: Reduction to surfaces}.
\item The nonpolar boundary $Q$ gives rise to two types of critical neighborhoods. Write $Q=Q'\cup Q''$, where $\xi|_{Q_j\setminus P}\equiv0$ for components $Q_j\subset Q'$, while $\xi|_{Q_j\setminus P}\not\equiv0$ for components $Q_j\subset Q''$. We put $Q'$ together with the zero set $Z(\xi)$ and construct Type I critical neighborhoods, following the local construction of Simpson \cite{Simpson93Lef}. For $Q''$, we introduce Type II critical neighborhoods around the points where the tangential restriction of $\xi$ vanishes, as well as around the transverse intersections of components of $Q''$ which cannot be handled by a globally tangent vector field. When components of $Q''$ meet a Type I neighborhood, we incorporate them into the corresponding mixed Type I situation. We then establish the required local $1$-connectedness statements for all these neighborhoods, always after deleting the relevant part of $Q$, denoted as $W^*_i,\fkW_l^*,W_i^\dagger$ for the type I, type II and mixed type respectively. This is carried out in \cref{subsec: Construction of two types of critical neighborhoods}.
\item We then pass to the regular $\Z$-covering $p:\widehat M\to M$ corresponding to $\chi$, and write $\widehat Y=p^{-1}(Y)$. The primitive function
$u:\widehat M\to\R,du=p^*\Re\xi$, is \textit{proper} in the rational case. We construct a uniformly controlled flow on $\widehat M$ which increases $u$ outside the lifted critical neighborhoods, is tangent to the polar boundary, and preserves the complement of the lifted nonpolar boundary. Given a path in $\widehat Y$ whose endpoints lie in a level set $L_t\subset \widehat Y$, the flow moves the portions outside the critical neighborhoods toward $L_t$. Inside the critical neighborhoods, the local $1$-connectedness in \cref{prop: local 1-connectedness for type I,prop: local 1-connectedness for type II,prop: local 1-connectedness for mixed,prop: local sup level 1 connected} allow us to replace the corresponding portions by paths in the same relative homotopy class while preserving the required control of $u$. The construction is arranged throughout so that the deformed path never crosses either $P$ or $Q$, and hence remains in $\widehat Y$. This is carried out in \cref{subsec: Construction of a flow on the infinite cyclic covering}.
\end{enumerate}

This proves that, when $\dim\psi(X)\ge2$, every regular level set of $u$ in $\widehat Y$ is connected and its inclusion induces a surjection onto $\pi_1(\widehat Y)$. Using $\pi_1(Y)\cong\pi_1(X)$ and the corresponding induced coverings, we transfer this conclusion from $\widehat Y$ to $\widehat X$, obtaining condition \textup{(ii)} of \cref{thm: Lefschetz alternative rational case} for regular levels. Additional local $1$-connectedness results (\cref{prop: local 1-connectedness central real level I,prop: local 1-connectedness central real level II,prop: local 1-connectedness central real level mixed}), together with the same flow argument, extend the conclusion to arbitrary level sets in \cref{subsec: From regular level set to any level set}. The resulting statement is then transferred from $\widehat Y$ to $\widehat X$ in the same way. On the other hand, when $\dim\psi(X)=1$, the morphism $\psi$ yields the factorization in condition \textup{(i)}.

\subsection{On semi-Albanese map determined by a real character}\label{subsec: On semi-Albanese map determined by a character}

This section is the quasi-projective analogue of Simpson's construction of the Albanese map determined by a holomorphic one-form in \cite{Simpson93Lef}. We work instead with the quasi-Albanese map and semi-abelian varieties. We refer the reader to \cite{Fujino25,Fujino25Correction} for the relevant properties used below.

Let $X$ be a smooth connected quasi-projective variety over $\C$, and let $a:X\to A_0$ be its quasi-Albanese map. The induced map
$a_*:H_1(X;\Z)/\operatorname{Tor}\to H_1(A_0;\Z)$
is an isomorphism; see, for example, \cite[Lemma~3.12]{Fujino25}. Thus there exits a unique integral class $\lambda_0\in H^1(A_0;\Z)$ such that $a^*\lambda_0=\chi$.

Let $B\subset A_0$ be the sum of all connected algebraic subgroups of $A_0$ where the restriction of $\lambda_0$ vanishes, which is a connected algebraic subgroup of $A_0$ and $\lambda_0|_B=0$. Set $A=A_0/B$, and $q: A_0\to A$ is the natural quotient map. Since $B$ is connected, $q_*: \pi_1(A_0)\to \pi_1(A)$ is surjective, and we have a unique (surjective) $\lambda\in H^1(A;\Z)$ such that $q^*\lambda=\lambda_0$. Set $\psi=q\circ a: X\to A$ be the \textit{quasi-Albanese map associated to $\chi$}. Then $\chi=\psi^*\lambda$.
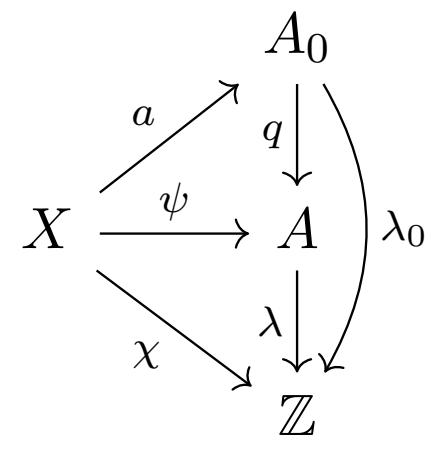
\[
\begin{tikzcd}
&A_0\ar[d,"q"']\ar[dd,bend left=30, "\lambda_0"]\\
X\ar[ur,"a"]\ar[dr,"\chi"']\ar[r,"\psi"] &A\ar[d,"\lambda"'] \\
&\Z
\end{tikzcd}
\]
We claim that $\lambda$ is nonzero on every positive-dimensional connected algebraic subgroup of $A$. Suppose that $H\subset A$ connected of positive-dimension, and $\lambda|_H\equiv 0$, then its inverse image $\widetilde{H}=q^{-1}(H)$ is a connected algebraic subgroup of $A_0$ containing $B$, and $\lambda_0|_{\tilde{H}} \equiv 0$; we have $B\subsetneq \widetilde{H}$ due to $\dim H>0$, contradicting the definition of $B$. Thus we have shown

\begin{lemma}\label{lem: nonvanishing on postive subgp after quotient}
For every positive dimensional connected algebraic subgroup $H\subset A$, the restriction of $\lambda$ to $H$ is not identically zero.
\end{lemma}

Moreover, for the semi-abelian variety $A$, we have the following decomposition:
$$1\to T\to A\to A_{\ab}\to 1$$
where $T\cong (\C^*)^\ell$ is an algebraic torus and $A_{\ab}$ is an abelian variety. For $\lambda|_T\in H^1(T;\Z)$, if it is identically zero, then by \cref{lem: nonvanishing on postive subgp after quotient}, we get $\ell=0$; otherwise, under the natural isomorphism $H^1(T;\Z)\cong \mathrm{Hom}(\pi_1(T);\Z)$, assume that $\lambda|_T$ is
$$\lambda|_T: (a_1,\cdots,a_\ell)\in \Z^\ell\mapsto a_1n_1+\cdots +a_\ell n_\ell, \quad  (n_1,\cdots, n_\ell)\in \Z^\ell\text{ not all zeros.}$$
Consider the associated algebraic character
$\tau: T\to\C^*$ given by $\tau(z_1,\ldots,z_\ell)=z_1^{n_1}\cdots z_\ell^{n_\ell}$, then $\tau_*:\pi_1(T)\to \Z$ equals $\lambda|_T$. If $\ell>1$, the identity component $H=(\ker\tau)^0$ is a positive-dimensional connected algebraic subgroup of dimension $\ell-1$, and $\lambda|_H\equiv 0$, contradicting \cref{lem: nonvanishing on postive subgp after quotient}. We therefore obtain:
\begin{lemma}\label{lem: rank of torus le 1}
The torus part of $A$ has dimension $\ell\le 1$.
\end{lemma}

\subsection{Choice of logarithmic holomorphic one-form}\label{subsec: Choice of logarithmic holomorphic one-form}
In case $\ell=1$, we can write $A=L^\times$ with $L$ a topological trivial line bundle (due to its translation-invariance) with fiber $\C^*$ on $A_\ab$. Set
\begin{equation}\label{eq: compact of A rational case}
    \bar{A}=\mathbb{P}_{A_\ab}(\calo_{A_\ab}\oplus L), \quad D_A\coloneqq \bar{A}\setminus A=D_0\sqcup D_\infty
\end{equation}
with $D_0$ and $D_\infty$ are zero and infinity sections respectively. For $\ell=0$, we take $\bar{A}=A$ and $D_A=\varnothing$. 

Recall that the mixed Hodge structure \cite{DeligneHodgeII} on $H^1(A;\C)$ has only Hodge types
$(1,0)$, $(0,1)$, and $(1,1)$. The weight-one part is pulled back from the abelian quotient $A_\ab$, while the weight-two part is contributed by the torus $T$. Consequently the real-part map
$$F^1H^1(A;\C)\to  H^1(A;\R),\quad\alpha\longmapsto\Re \alpha,$$
is surjective. We may now apply Deligne's logarithmic de Rham theory \cite{DeligneHodgeII} to choose a logarithmic holomorphic one-form whose real part represents $\lambda$ on $A$ and has controlled behavior along the boundary divisor $D_A$.

\begin{lemma}\label{lem: good choice of holomorphic one-form}
\begin{enumerate}
\item There is a closed algebraic one-form $\omega$ on $A$ such that $[\Re\omega]=\lambda$, and $\omega$ extends logarithmically to the standard projective compactification of $A$ with purely imaginary residues. 
\item For every nonconstant algebraic map $u:C\to A$ from a smooth algebraic curve, $u^*\omega\not\equiv0$.
\end{enumerate}
\end{lemma}
\begin{proof}
(a) In case $\ell=0$, we take $\omega\in H^0(A;\Omega_A^1)$ a holomorphic 1-form such that $[\Re \omega]=\lambda$. In case $\ell=1$, consider 
$$1\to \C^*\to A\xrightarrow{\pi}A_\ab\to1,$$
we have $W_1 H^1(A;\C)=p^* H^1(B;\C)$ and $\mathrm{Gr}^W_2H^1(A;\C)\cong \C[-1]$.
Consider the residue sequence 
$$0\to \Omega^1_{\bar A}
\to \Omega^1_{\bar A}(\log D_A )
\xrightarrow{\operatorname{Res}}
\mathcal O_{D_0}\oplus\mathcal O_{D_\infty}
\to 0,$$
then the toric part $\C^*$ is detected by the residue map
$$\Res:
H^0\!\left(\bar A,\Omega_{\bar A}^1(\log D_A )\right)
\to 
\{(c,-c)\in\C^2:c\in\C\}.$$
%Indeed, since $D_0$ is connected and compact, $\Res_{D_0}\omega$ is constant, say $c_0$, on $D_0$; similarly, write $\Res_{D_\infty}\omega=c_\infty$. Restricting $\omega$ to a fiber $\mathbb{P}^1$ of $\bar A\to A_{\ab}$, the residue theorem gives that $c_0+c_\infty=0$. 
It is standard that there exists a global logarithmic $1$-form $\eta$ on $A$ extending logarithmically to $\bar A$ and satisfying
$$\Res_{D_0}\eta=1,\quad \Res_{D_\infty}\eta=-1.$$
%For example, this follows from the existence of a holomorphic flat connection on the line bundle $L\in\Pic^0(A_{\ab})$.
Locally near a point of $D_0$, if $z$ is local defining function of $D_0$, we can write $\eta=\frac{dz}{z}+\beta$, with $\beta$ holomorphic.

Let $\gamma_0$ be a small positively oriented loop around $D_0$, and set $\lambda(\gamma_0)=n\in \Z$, we have
$$\int_\mu \frac{n}{2\pi \sqrt{-1}}\eta=n=\lambda(\mu)$$
and $\lambda-\left[\Re \frac{n}{2\pi \sqrt{-1}}\eta\right]$ has trivial image in $\Gr^W_2(H^1(A;\C))$, thus there exists $\lambda_\ab\in H^1(A_\ab;\R)$ such that 
$$\lambda-\left[\Re \frac{n}{2\pi\sqrt{-1}}\eta\right]=\pi^*\lambda_\ab$$
Choose a holomorphic one form $\omega_\ab$ such that $[\Re \omega_\ab]=\pi^*\lambda_\ab$, then
\begin{equation}\label{eq: log omega}
\omega=\frac{n}{2\pi\sqrt{-1}}\eta+\omega_\ab
\end{equation}
is a logarithmic one form on $A$ such that $[\Re \omega]=\lambda$, and $\omega$ has purely imaginary residues at divisors $D_0\sqcup D_\infty$:
$$\Res_{D_0}\omega=
\frac{n}{2\pi\sqrt{-1}},
\quad\Res _{D_\infty}\omega
=-\frac{n}{2\pi\sqrt{-1}}.$$

(b) Assume that $u^*\omega=0$. Up to translation, the universal property of the quasi-Albanese morphism \cite[Theorem~1.1]{Fujino25} gives a factorization $u=\phi\circ a_C$, where $\phi:\Alb(C)\to A$ is a homomorphism and $a_C$ is the quasi-Albanese morphism of $C$. Let
$H=\operatorname{im}(\phi)$. Since $u$ is nonconstant, $H$ is a positive-dimensional connected semi-abelian subvariety of $A$. The map $\Alb(C)\to H$ induces a surjection on $H_1(-;\R)$, and $a_C$ also induces a surjection on $H_1$. Since
$u^*\lambda=[\Re(u^*\omega)]=0$, it follows that $\lambda|_H=0$, contradicting \cref{lem: nonvanishing on postive subgp after quotient}.
\end{proof}

\begin{remark}\label{rmk: residue not zero rational case}
Whenever we use the expression \cref{eq: log omega}, we have $n\ne0$. Indeed, in the case $\ell=1$, the torus $T\simeq\C^*$ is a positive-dimensional connected algebraic subgroup of $A$, and hence $\lambda|_T\ne0$ by \cref{lem: nonvanishing on postive subgp after quotient}. Thus the corresponding residues of $\omega$ along $D_0$ is nonzero.
\end{remark}

\begin{remark}\label{rmk: accumulating zero sets}
A good choice such as a logarithmic holomorphic one-form is essential for the argument below, since an arbitrary holomorphic $1$-form on a smooth quasi-projective variety may have zeros accumulating toward infinity. For example, the holomorphic one-form $\sin(1/z)\,dz$ on $\C^*$ has zeros at $z=1/(k\pi)$, $k\in\Z\setminus\{0\}$, which accumulate at $0$.

Such accumulation prevents us, in general, from choosing critical neighborhoods with a uniform positive width in \cref{subsec: Construction of two types of critical neighborhoods}, and therefore creates a serious difficulty for the flow argument used in \cref{subsec: Construction of a flow on the infinite cyclic covering}. The same difficulty also obstructs a direct extension of Simpson's global picture and patching argument to the noncompact setting, since that construction likewise requires uniform control of the critical neighborhoods.
\end{remark}

\subsection{Reduction to surfaces and real-oriented blowup}\label{subsec: Reduction to surfaces}
As in \cite{Simpson93Lef}, we consider the following two cases: $\dim \psi(X)=1$ and $\dim \psi(X)\ge 2$. They will corresponding to (i) and (ii) of \cref{thm: Lefschetz alternative rational case} correspondingly. For now, we consider first the case $\dim \psi(X)\ge 2$.

We may first reduce to the case where $X$ is a smooth quasi-projective surface as follows. By taking a general complete intersection of hyperplanes, due to the Lefschetz hyperplane theorem (for example, see \cite[Theorem 1.1.3]{HammLe85}) we obtain a smooth quasi-projective surface $Y\subset X$ such that the inclusion induces an isomorphism: 
\begin{equation}\label{eq: iso on pi1}
\pi_1(Y)\xrightarrow{\cong}\pi_1(X).
\end{equation}
Moreover, since $\dim\psi(X)\ge2$, the hyperplanes may be chosen sufficiently generally so that, at a general point $x\in Y$ at which $d\psi_x$ has maximal rank,
$T_xY\cap\ker(d\psi_x)=0$. Thus $\dim\psi(Y)=2$ and $\psi|_Y$ is generically finite onto its image. More precisely, we have
\begin{lemma}\label{lem: surface through prescribed points}
Let $X$ be a smooth connected complex quasi-projective variety, and let $\psi:X\to A$ be an algebraic morphism with $\dim\psi(X)\ge2$. For every finite subset $S\subset X$, there exists a smooth connected quasi-projective surface $Y\subset X$ containing $S$ such that
$$\pi_1(Y)\xrightarrow{\cong}\pi_1(X)
\quad\text{and}\quad
\dim\psi(Y)=2.$$
\end{lemma}
We omit the standard proof, which uses a suitable very ample line bundle together with Bertini's theorem and the quasi-projective Lefschetz hyperplane theorem.

%\begin{proof}
%If $\dim X=2$, take $Y=X$. Suppose that $n=\dim X\ge3$. Enlarge $S$, if necessary, by adding a point $q$ at which $d\psi_q$ has rank at least two.

%Choose a smooth projective compactification $\bar X$ of $X$, and let $b:X'=\operatorname{Bl}_S(\bar X)\to\bar X$ be the blowup at $S$, with exceptional divisor $E$. Put $\widetilde X=b^{-1}(X)$. If $H$ is a very ample line bundle on $\bar X$, then, for $m$ sufficiently large, \[\mathcal L_m=b^*H^{\otimes m}\otimes\mathcal O_{X'}(-E)\] is very ample. We may also choose $m$ sufficiently large that sections of $H^{\otimes m}$ vanishing on $S$ have arbitrary first-order terms at the points of $S$.

%Take a sufficiently general complete intersection $Y'$ of $n-2$ members of $|\mathcal L_m|$, and set $\widetilde Y=Y'\cap\widetilde X$. By successive applications of the quasi-projective Lefschetz hyperplane theorem \cite[Theorem~1.1.3(ii)]{HammLe85}, the inclusion induces $\pi_1(\widetilde Y)\xrightarrow{\cong}\pi_1(\widetilde X)$.

%The corresponding divisors on $\bar X$ have independent linear terms at every point of $S$. Consequently, $Y=b(\widetilde Y)$ is smooth and contains $S$, and $\widetilde Y\to Y$ is the blowup of $Y$ at $S$. Since blowing up finitely many smooth points does not change the fundamental group, the preceding isomorphism gives $\pi_1(Y)\cong\pi_1(X)$.

%Finally, the first-order terms at $q$ may be chosen so that $d(\psi|_Y)_q$ has rank two. This is an open condition compatible with the preceding generality conditions. Hence $\dim\psi(Y)=2$.
%\end{proof}

We may choose smooth projective compactifications $Y\subset\bar Y$ and $A\subset\bar A$ such that $\psi:Y\to A$ extends to a morphism $\bar\psi:\bar Y\to\bar A$, and such that
$$D_A=\bar A\setminus A,\quad
D_Y=\bar Y\setminus Y,\quad
P=\bigl(\bar\psi^{-1}(D_A)\bigr)_{\red}$$
are simple normal crossings divisors, with $P$ a union of irreducible components of $D_Y$. Let $Q$ be the reduced union of the remaining irreducible components of $D_Y$, so that $D_Y=P\cup Q$.
%Equivalently, an irreducible component $D\subset D_Y$ belongs to $P$ if and only if $$\bar\psi(D)\subset D_A,$$ whereas for every irreducible component $Q_j\subset Q$, the generic point of $Q_j$ is mapped into $A$. In particular, $$Q_j\setminus P\subset\bar\psi^{-1}(A),$$ while $Q_j\cap P$ is a finite set.

Since $\psi:Y\to A$ is generically finite onto its two-dimensional image, the extension $\bar\psi:\bar Y\to\bar A$ is also generically finite onto its image. Let $Y^\circ =\bar Y\setminus P$ and $\psi^\circ=\bar \psi|_{Y^\circ}$. Let $E$ denote the union of the positive-dimensional fibers of $\psi^\circ:Y^\circ\to A$. Then $E$ is a finite union of irreducible curves in $Y^\circ$, and it contains both the closures in $Y^\circ$ of the exceptional curves of $\psi$ in $Y$ and any irreducible component of $Q$ contracted by $\bar\psi$. We call any irreducible curve $C\subset E$ an \textit{exceptional curve}, and call $E$ the \textit{exceptional set}.
%Each connected component of $E$ is mapped by $\psi^\circ$ to a single point $a\in A$.

Due to \cref{lem: good choice of holomorphic one-form}(b), we have
\begin{lemma}\label{lem: out of excep set omega does not vanish}
For any irreducible curve $C\subset Y^\circ$ such that $C\not\subset E $,
$$\bar\psi^*\omega|_C\not\equiv0.$$
\end{lemma}

Recall that for the logarithmic holomorphic $1$-form $\omega$ on $A$ defined in \cref{subsec: Choice of logarithmic holomorphic one-form}, locally near a smooth point of a divisor $D_0\subset D_A$ written as $\{z=0\}$, it has the form
$$\omega=\frac{n}{2\pi\sqrt{-1}}\frac{dz}{z}+dh$$
where $h$ is non-vanishing holomorphic and $n\ne0$ (see \cref{rmk: residue not zero rational case}).

For $x_0\in A$ fixed, define a function on $A$:
$$\Phi_A(x)=\exp\left(2\pi\sqrt{-1}\int_{x_0}^x \Re \omega\right),\quad \forall x\in A,$$
which is well defined because $\lambda=[\Re\omega]\in H^1(A;\Z)$ is integral. Usually, this function cannot be extended to $\bar A$, for near $D_0$ we have
$$\Re \omega=\Re\left(\frac{n}{2\pi\sqrt{-1}}\left(\frac{dr}{r}+\sqrt{-1}d\theta\right)+dh\right)=\frac{n}{2\pi}d\theta+\Re(dh),$$
where $z=re^{\sqrt{-1}\theta}$ and $x_0$ is fixed (here we take $1$ for example), and
$$\Phi_A(x)=\exp(\sqrt{-1}n\theta)\cdot \exp\left(2\pi\sqrt{-1}\int_{x_0}^x\Re(dh)\right).$$
When $r$ tends to $0$, the second factor tends to a limit while the first factor is not well defined on $\bar A$. To remedy this, we use the real-oriented blowup; then the function $\Phi_A$ extends to the boundary. In short, the real-oriented blowup allows the angular coordinate $\theta$ to remain well defined even when $r=0$. Of course, the resulting manifold is no longer a complex manifold in general. This construction is also the Kato--Nakayama space \cite[(1.2.3)]{Kato-Nakayama} associated with the log structure defined by the SNC divisor $P$. For simplicity, in our special case we describe the local construction and gluing process explicitly below. If the torus part is trivial in $A$, i.e., $\ell=0$, then $\bar A=A$, $D_A=P=\varnothing$, and $M=\bar Y$. For the following, we assume $\ell=1$.

Since $P\subset\bar Y$ is a simple normal crossing divisor, we define its \textit{real-oriented blowup}
$$M=\operatorname{Bl}^{\mathrm{ro}}_P(\bar Y),\quad \beta:M\to\bar Y$$
locally as follows. If $P=\{z_1\cdots z_k=0\}$, $k=1,2$, write $z_j=r_je^{\sqrt{-1}\theta_j}$, with $r_j\ge0$, and replace each normal coordinate $z_j$ by the pair $(r_j,\theta_j)$. Thus locally
$$M\cong(\mathbb R_{\ge0}\times \bbS^1)^k\times\mathbb C^{2-k},$$
where $\bbS^1$ is the unit circle in $\C$. These local models glue canonically. Indeed, on the overlap of two SNC coordinate charts, after possibly reordering, their defining coordinates satisfy $z_j'=u_jz_j$ for nowhere-vanishing holomorphic functions $u_j$. Hence, writing $z_j=r_je^{\sqrt{-1}\theta_j}$,
$$r_j'=|u_j|r_j,\quad e^{\sqrt{-1}\theta_j'}=\frac{u_j}{|u_j|}e^{\sqrt{-1}\theta_j}\in\bbS^1.$$
These transition maps extend smoothly to $r_j=0$, and therefore the local oriented blowups glue to a manifold with boundary. The boundary with corners can be described as follows. In a coordinate neighborhood near $\{z_1z_2=0\}$, the boundary is of the form $\{r_1=0\}\cup\{r_2=0\}$, and the corner $C$ is of the form $\{r_1=r_2=0\}$ and $C\cong(\bbS^1)^2$.

Composing with $\bar\psi:\bar Y\to\bar A$, we get $\Phi:Y^\circ \to\bbS^1$ defined by
$$\Phi(x)=\exp\left(2\pi\sqrt{-1}\int_{x_0}^x\Re\bar \psi^* \omega \right).$$
On $Y$, this is simply
$$\Phi(x)=\exp\left(2\pi\sqrt{-1}\int_{x_0}^x\Re\psi^*\omega\right).$$

\begin{proposition}\label{prop: no critical at boundary}
The $\bbS^1$-valued map $\Phi:Y^\circ \to\bbS^1$ extends smoothly to the boundary of $M$. Moreover, its restriction to each boundary component $\beta^{-1}(P)$ and each boundary corner of $M$ is a submersion.
\end{proposition}

\begin{proof}
For the case $\ell=0$ and $M=\bar Y$, there is nothing to prove. For the case $\ell=1$, by \cref{eq: compact of A rational case}, take $x\in P$ with $\bar\psi(x)\in D_A$; without loss of generality, assume $\bar\psi(x)\in D_0$. Let $t$ be a local coordinate around $\bar\psi(x)$ such that $D_0=\{t=0\}$ and $\{z_1,z_2\}$ be local coordinates around $x$, with $P=\{z_1\cdots z_k=0\}$, $k=1,2$. Then locally we have
$$\bar\psi^*t=vz_1^{n_1}\cdots z_k^{n_k}$$
with $v$ a non-vanishing holomorphic function, and $n_j$, $1\le j\le k$, positive integers. Recall that from \cref{eq: log omega}, locally we write
$$\omega=\pm\frac{n}{2\pi\sqrt{-1}}\frac{dt}{t}+\omega_\ab,\quad n\ne0,$$
with $\pm$ for $D_0$ and $D_\infty$, respectively. Then
\begin{equation}\label{eq: b pull-back of omega}
\bar\psi^*\omega=\pm \frac{n}{2\pi\sqrt{-1}}\left(\frac{dv}{v}+\sum_{j=1}^k n_j\frac{dz_j}{z_j}\right)+\bar\psi^*\omega_\ab.
\end{equation}
The part $\pm \frac{n}{2\pi\sqrt{-1}}\frac{dv}{v}+\bar\psi^*\omega_\ab$ is a closed holomorphic $1$-form on $\bar Y$ near $x$, which can be written as $dg$ locally with $g$ a holomorphic function. Then
\begin{equation}\label{eq: Phi downstairs}
\Phi=c\exp\left(\sqrt{-1}\sum_{j=1}^k m_j\theta_j\right)\exp(2\pi\sqrt{-1}\Re g),
\end{equation}
with $m_j=\pm nn_j\ne0$ and $c\in\bbS^1$ a constant depending on the base point $x_0$. This is a smooth function in the polar coordinates $(r_1,\theta_1;\cdots;r_k,\theta_k)$, $k=1,2$, with $r_j\ge0$, and therefore extends smoothly to the boundary of $M$.

It remains to show that the extended map does not have any critical point on the boundary $\beta^{-1}(P)\subset M$. Away from its corner, we can write its local equation as $r_j=0$, then we have
$$\frac{\partial}{\partial\theta_j}\Re(g)=\frac12\left(\frac{\partial g}{\partial z_j}\sqrt{-1}z_j+\frac{\partial\bar g}{\partial\bar z_j}(-\sqrt{-1}\bar z_j)\right)=\Re\left(\frac{\partial g}{\partial z_j}\sqrt{-1}z_j\right)=0,\quad\text{when }r_j=0,$$
thus
$$\frac{\partial}{\partial\theta_j}(\arg\Phi)=m_j\ne0.$$
Thus $\Phi$ is a submersion on each boundary away from the corner. Near the corner locally defined by $r_1=r_2=0$, we have
$$\frac{\partial}{\partial\theta_1}(\arg\Phi)=m_1\ne0,\quad \frac{\partial}{\partial\theta_2}(\arg\Phi)=m_2\ne0.$$
Any one of these two equations implies that $\Phi$ restricted to the corner is a submersion.
\end{proof}

%\begin{remark}
%Recall that $$D_Y=P\cup Q.$$ Since $Q$ is nonpolar, we have $Q\setminus P\subset\bar Y\setminus P=\bar\psi^{-1}(A)$, where $\Phi$ is already defined. Moreover, the exceptional curves $E$ lie in $Y\subset\bar Y\setminus P$. Thus, in extending $\Phi$ to $M$, it suffices to consider the polar divisor $P$.
%\end{remark}

\begin{remark}\label{rmk: poincare metric replacement}
The use of the real-oriented blowup along the polar divisor $P$ is convenient but not essential for the argument. One may instead work directly on $Y^\circ$ and equip it with a complete metric of Poincar\'e type near $P$. In local coordinates $(z_1,z_2)$, if $P=\{z_1=0\}$, a Poincar\'e-type metric can be written as
$$g_P\sim \frac{|dz_1|^2}{|z_1|^2(\log|z_1|)^2}+|dz_2|^2.$$
Hence one may choose the vector field used in \cref{subsec: Construction of a flow on the infinite cyclic covering} to have uniformly bounded norm and to satisfy $d\arg\Phi(\bfv)=2\pi$ away from the critical neighborhoods defined in \cref{subsec: Construction of two types of critical neighborhoods}. Since a uniformly bounded smooth vector field on a complete Riemannian manifold has complete flow (\cite[Proposition 1.4.4(v)]{CTCWall}), no integral curve can reach the divisor $P$ in finite time. Thus the role played above by the submersion of $\Phi$ on the boundary of $\beta^{-1}(P)$ can be replaced by completeness with respect to a Poincar\'e-type metric.
\end{remark}

\subsection{Construction of two types of critical neighborhoods and local 1-connectedness}\label{subsec: Construction of two types of critical neighborhoods}
Set
$$\xi\coloneqq\bar\psi^*\omega|_{Y^\circ}\mathrel{}=:(\psi^\circ)^*\omega,$$
which is a holomorphic $1$-form on $Y^\circ$. We now consider its zero set $Z(\xi)\subset Y^\circ$. By \cref{lem: out of excep set omega does not vanish}, every positive-dimensional irreducible component of $Z(\xi)$ must be contained in $E$. Hence $Z(\xi)\setminus E$ is zero-dimensional, which is a finite set of points. Write
$$Z(\xi)=A_1\sqcup\cdots\sqcup A_n,$$
where each $A_i$ is either a point or a connected union of irreducible curves contained in $E$. 

Since we are now dealing with the non-compact case, we need to be extremely careful about the boundary. Recall from \cref{subsec: Reduction to surfaces} that
$$\bar Y\setminus Y=P\cup Q\quad\text{and}\quad Y^\circ=\bar Y\setminus P=Y\cup(Q\setminus P).$$
So far in the present section, we have been working on $Y^\circ$, while our ultimate goal is to obtain the corresponding conclusion on $Y$. Thus, in the flow construction later in \cref{subsec: Construction of a flow on the infinite cyclic covering}, we need to ensure that the flow preserves the complement of $Q$, so that a trajectory starting in $Y$ remains in $Y$ (or for the suitable $\Z$-coverings). This \textit{cannot} always be achieved simply by requiring the flow to be tangent to $Q$. Write $Q=Q_1\cup\cdots\cup Q_l$ as the union of its irreducible components. For each $Q_j$, there are two possibilities:
\begin{itemize}
    \item $\xi|_{Q_j\setminus P}\equiv0$;
    \item $\xi|_{Q_j\setminus P}\not\equiv0$.
\end{itemize}
We need to deal with these two cases separately. For the first case, this was considered by Simpson in \cite[Theorem 1]{Simpson93Lef} (defined as the deleted subset $D$ there). For this part, we shall put these components together with $Z(\xi)$ and construct the critical neighborhoods as done by Simpson; we will call these \textit{critical neighborhoods of type I}. For the second case, the zero set of $\xi|_{Q_j}$ is discrete. Around each of these finitely many points, we choose a suitable neighborhood such that outside these neighborhoods we can choose the flow tangent to $Q$. We call these \textit{critical neighborhoods of type II}.

Denote the union of the $Q_j$ in the first case by $Q'$, and the union of those in the second case by $Q''$. Set
$$\Sigma=Z(\xi)\cup Q'$$
and decompose it into connected components as
$\Sigma=\bigsqcup_{i=1}^N\Sigma_i$. Note that for every $Q_j\subset Q'$, the restriction
$\bar\psi|_{Q_j\cap Y^\circ}$ is constant by
\cref{lem: out of excep set omega does not vanish}. Hence
$Q_j\subset E$ and $\bar\psi(Q_j)=\{a_j\}$ for some $a_j\in A$.
In particular, $Q'\cap P=\varnothing$, and the first condition is
reduced to $\xi|_{Q_j}\equiv 0$.

Since \cref{prop: no critical at boundary} implies that
$\bar\psi^*\omega$ has no zero in a neighborhood of the polar divisor
$P$, the zero locus of $\xi=\bar\psi^*\omega|_{Y^\circ}$ is contained
in a compact subset of $Y^\circ$. Thus the problem of constructing
critical neighborhoods of type I around $\Sigma$ reduces to the
compact situation considered in \cite{Simpson93Lef}.

For each $\Sigma_i$, let $a_i=\psi^\circ(\Sigma_i)\in A$.
Choose a sufficiently small simply connected holomorphic coordinate
neighborhood $B_i$ of $a_i$ in $A$. Since $\omega$ is closed, there
exists a holomorphic function $h_i$ on $B_i$ such that
$dh_i=\omega$ and $h_i(a_i)=0$. Choose pairwise disjoint relatively
compact neighborhoods $V_i'$ of $\Sigma_i$ such that
$$
\overline{V_i'}\subset(\psi^\circ)^{-1}(B_i)
\quad\text{and}\quad
\overline{V_i'}\cap\Sigma=\Sigma_i.
$$
Define $f_i=h_i\circ\psi^\circ:V_i'\to\C$. Then
$df_i=\xi|_{V_i'}$ and $f_i|_{\Sigma_i}=0$.

Choose a smooth nonnegative function $\rho_i$ on $V_i'$ such that $\rho_i^{-1}(0)=\Sigma_i$ and its sufficiently small sublevel sets are compactly contained in $V_i'$ \footnote{For example, assume $\Sigma_i$ is defined locally as $f_i=0,1\le i\le n$; take $\sum_{i=1}^n|f_i|^2$; using a smooth partition of unity to patch these function to get $\rho_i$, with $\rho_i=1$ outside a relatively compact neighborhood of $\Sigma_i$ in $V_i'$.}. Since all the critical points of $f_i$ in $V_i'$ belong to $\Sigma_i$, the set $f_i^{-1}(0)\setminus\Sigma_i$ is smooth. By Sard's theorem, we may choose $\varepsilon>0$ sufficiently small to be a regular value both of $\rho_i$ and of its restriction to $f_i^{-1}(0)\setminus\Sigma_i$, simultaneously for all $i$, $1\le i \le N$.

Let $U_i$ be the connected component of
$\rho_i^{-1}([0,\varepsilon])$ containing $\Sigma_i$. Then $U_i$ is
a compact smooth manifold with smooth boundary, and $f_i^{-1}(0)$
intersects $\partial U_i$ transversely whenever their intersection
is nonempty. Moreover, $U_i\cap U_j=\varnothing$ when $i\ne j$,
and $\Sigma_i$ is the unique connected component of $\Sigma$
contained in $U_i$, for each $1\le i\le N$.

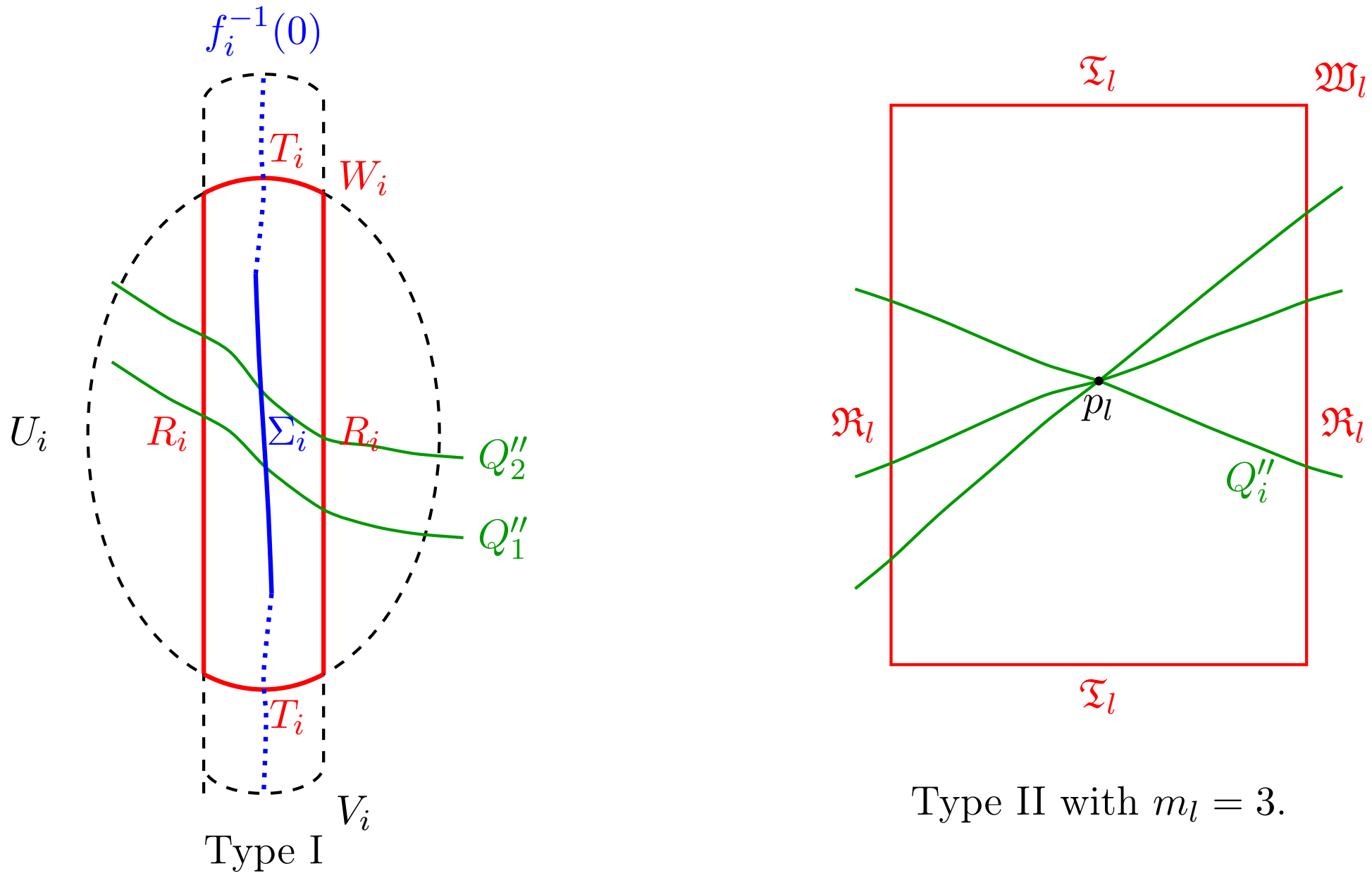
\begin{figure}[htp]
    \centering

    %-------------------- Type I --------------------
    \begin{minipage}[c]{0.48\textwidth}
    \centering
    \begin{tikzpicture}[scale=.75]

    % U_i: ball-shaped neighborhood
    \draw[dashed,thick]
        (0,0) ellipse [x radius=2.2, y radius=3.2];
    \node[left] at (-2.5,0) {$U_i$};

    % V_i: long cylindrical neighborhood
    \draw[dashed,thick]
        (-0.75,-4.5) -- (-0.75,4.15)
        arc[start angle=180,end angle=0,
            x radius=.75,y radius=.35]
        -- (0.75,-4.15)
        arc[start angle=0,end angle=-180,
            x radius=.75,y radius=.35];

    \node[below right] at (0.72,-4.35) {$V_i$};

    % Boundary of W_i = U_i \cap V_i
    % Vertical pieces inherited from V_i
    \draw[red,very thick]
        (-0.75,-3.0) -- (-0.75,3.0);
    \draw[red,very thick]
        (0.75,-3.0) -- (0.75,3.0);

    % Upper and lower pieces inherited from U_i
    \draw[red,very thick]
        (0.752,3.007)
        arc[start angle=70,end angle=110,
            x radius=2.2,y radius=3.2];

    \draw[red,very thick]
        (-0.752,-3.007)
        arc[start angle=250,end angle=290,
            x radius=2.2,y radius=3.2];

    % Three nearby fibers / leaves

    \draw[green!60!black,thick]
        plot[smooth] coordinates {
            (-1.9,0.9) (-1.2,0.45) (-0.45,0.05) (0.06,-0.45) (0.75,-0.95) (1.35,-1.15) (1.9,-1.25) (2.5, -1.3)
        };
    \node[green!60!black, right] at (2.5,-1.3) {$Q_{1}''$};
    \draw[green!60!black,thick]
        plot[smooth] coordinates {
            (-1.9,1.9) (-1.2,1.45) (-0.45,1.05) (0.06,.45) (0.75,-.05) (1.35,-0.15) (1.9,-.25) (2.5, -.3)
        };
    \node[green!60!black, right] at (2.5,-.3) {$Q_{2}''$};

    % Central fiber f_i^{-1}(0)
    \draw[blue,very thick]
        plot[smooth] coordinates {
            (-0.10,2.0)
            (-0.06,1.0)
            (0,0)
            (0.06,-1.0)
            (0.10,-2.0)
        };

    % Dotted continuation of the central fiber
    \draw[blue,very thick,dotted]
        plot[smooth] coordinates {
            (-0.10,2.0)
            (0,3.0)
            (-0.03,3.5)
            (-0.02,4.0)
            (0,4.5)
        };

    \draw[blue,very thick,dotted]
        plot[smooth] coordinates {
            (0.10,-2.0)
            (0,-3.0)
            (0.03,-3.5)
            (0.02,-4.0)
            (0,-4.5)
        };

    % Labels
    \node[blue,above] at (0,4.55) {$f_i^{-1}(0)$};

    \node[right,red] at (0.75,3.20) {$W_i$};

    % T_i: boundary pieces inherited from U_i
    \node[red,above] at (0.3,3.15) {$T_i$};
    \node[red,below] at (0.3,-3.15) {$T_i$};
    \node[blue] at (0.3,0) {$\Sigma_i$};

    % R_i: boundary pieces inherited from V_i
    \node[red,left] at (-0.75,0) {$R_i$};
    \node[red,right] at (0.75,0) {$R_i$};

    \node at (0,-5.25) {Type I};

    \end{tikzpicture}
    \end{minipage}
    \hfill
    %-------------------- Type II --------------------
    \begin{minipage}[c]{0.48\textwidth}
    \centering
    \begin{tikzpicture}[scale=.75]

    % Type II neighborhood
    \draw[red,thick]
        (-2.6,-3.5) rectangle (2.6,3.5);

    \node[above right, red] at (2.52,3.42)
        {${\fkW}_l$};

    % Boundary labels
    \node[red,above] at (0,3.5) {$\fkT_l$};
    \node[red,below] at (0,-3.5) {$\fkT_l$};

    \node[red,left] at (-2.6,-0.5) {$\fkR_l$};
    \node[red,right] at (2.6,-0.5) {$\fkR_l$};

        % Three smooth branches of Q through p_l.
    % Each branch crosses the R_l boundaries and avoids the T_l boundaries.
    % First branch
    \draw[green!60!black, thick]
        plot[smooth] coordinates {
            (-3.05,1.20)
            (-2.6,1.05)
            (-2.0,0.82)
            (-1.3,0.52)
            (-0.65,0.25)
            (0,0.05)
            (0.65,-0.22)
            (1.35,-0.52)
            (2.0,-0.78)
            (2.6,-1.02)
            (3.05,-1.15)
        };
    \node[green!60!black,below] at (1.9,-0.75) {$Q_i''$};

    % Second branch
    \draw[green!60!black, thick]
        plot[smooth] coordinates {
            (-3.05,-1.15)
            (-2.6,-0.98)
            (-2.0,-0.72)
            (-1.3,-0.40)
            (-0.65,-0.12)
            (0,0.05)
            (0.65,0.28)
            (1.35,0.58)
            (2.0,0.82)
            (2.6,1.05)
            (3.05,1.18)
        };

    % Third branch
    \draw[green!60!black, thick]
        plot[smooth] coordinates {
            (-3.05,-2.55)
            (-2.6,-2.18)
            (-2.0,-1.62)
            (-1.3,-1.02)
            (-0.65,-0.45)
            (0,0.05)
            (0.65,0.58)
            (1.3,1.12)
            (2.0,1.68)
            (2.6,2.15)
            (3.05,2.48)
        };
    % Type II point p_l
    \node[below] at (0,0.05) {$p_l$};
    \fill (0,0.05) circle (1.6pt);

    \node at (0,-5.25) {Type II with $m_l=3$.};

    \end{tikzpicture}
    \end{minipage}

    \caption{Critical neighborhoods of Types I and II (with possible mixed cases).}
    \label{fig: typeI-typeII-critical-neighborhoods}
\end{figure}

We take a small number $\delta>0$ and define $V_i=f_i^{-1}(\overline{\Delta(0;\delta)})$, where $\Delta(0;\delta)$ is a disk in $\C$ of radius $\delta$. By choosing $\delta$ sufficiently small, we may assume that $f_i^{-1}(v)$ intersects $\partial U_i$ transversely for every $v\in\overline{\Delta(0;\delta)}$, whenever $\partial U_i\cap f_i^{-1}(0)\ne \varnothing$. If $\partial U_i\cap f_i^{-1}(0)=\varnothing$, we may further arrange that the connected component of $V_i$ containing $\Sigma_i$ is contained in $\interior(U_i)$. Then $U_i\cap V_i$ is a smooth compact manifold with corners. We denote by $W_i$ the unique connected component of $U_i\cap V_i$ containing $\Sigma_i$. Write
$$T_i=\partial U_i\cap V_i\cap W_i=\partial U_i\cap f_i^{-1}(\overline{\Delta(0;\delta)})\cap W_i$$
and
$$R_i=\partial V_i\cap U_i\cap W_i=f_i^{-1}(\{z:|z|=\delta\})\cap U_i\cap W_i.$$
Their intersection $R_i\cap T_i$ is a smooth corner. We call these
$W_i$ \textit{critical neighborhoods of type I}; see
\cref{fig: typeI-typeII-critical-neighborhoods}.

There are now two mutually exclusive possibilities:\footnote{If $\Sigma=Z(\xi)\cup Q'=\varnothing$, there are no critical neighborhoods of type I. Condition \ref{it: condition A} is then false, while
condition \ref{it: condition B} holds vacuously. In this case, the steps involving type I neighborhoods are simply omitted.}
\begin{enumerate}[label=(\Alph*)]
    \item \label{it: condition A} For some $i$, $f_i^{-1}(0)\cap \partial U_i\cap W_i=\varnothing$;
    \item \label{it: condition B} For each $i$, $f_i^{-1}(0)\cap \partial U_i\cap W_i\ne  \varnothing$;
\end{enumerate}

In \cref{subsec: Reduction to surfaces}, under the assumption $\dim\psi(X)\ge2$, we reduced to the case
$\dim \bar\psi(Y)=2$. We claim that Case \ref{it: condition A} cannot occur. Assume not, up to shrinking $\delta>0$ if necessary, we have $f_i^{-1}(\overline{\Delta(0;\delta)})\cap W_i$ is compactly contained in $\interior(U_i)$. Hence, for every sufficiently small nonzero regular value $v$ of $f_i$, each connected component $C\subset f_i^{-1}(v)\cap W_i$ is a compact smooth complex curve. Since the exceptional set $E$ has only finitely many irreducible components and $f_i$ is constant on each irreducible component of $E\cap W_i$, the set $f_i(E\cap W_i)$ is finite. Choose a nonempty open subset $V\subset\Delta^*(0;\delta)\setminus f_i(E\cap W_i)$ consisting of regular values of $f_i$. Then $f_i^{-1}(V)\cap W_i\cap E=\varnothing$. For any $v\in V$ and any connected component $C\subset f_i^{-1}(v)\cap W_i$, we therefore have $C\not\subset E$. By \cref{lem: out of excep set omega does not vanish}, the restriction $(\bar\psi|_C)^*\omega\not\equiv0$. On the other hand, since $df_i=\xi=\bar\psi^*\omega$, we have $\xi|_C=0$, a contradiction. Therefore Case \ref{it: condition A} cannot occur.

From now on, we are under case \ref{it: condition B}, and we have the following local 1-connectedness property which is \cite[Lemma 8]{Simpson93Lef}. We include its proof here for the convenience of readers.
\begin{proposition}[Simpson]\label{prop: local 1-connectedness for type I}
Assume that Case \ref{it: condition B} holds. Up to shrinking $\delta$, assume that $f_i^{-1}(z)\cap W_i$ meets $\partial U_i$ transversely for each $z\in\overline{\Delta(0;\delta)}$. Let 
$$W_i^*\coloneqq W_i\setminus Q',$$
and let $f_i^{-1}(v)$ be a nearby fiber with $v\ne0$. Then pair
$$(W_i^*,f_i^{-1}(v)\cap W_i^*)=(W_i^*,f_i^{-1}(v)\cap W_i)$$
is $1$-connected.
\end{proposition}
\begin{proof}
Set
$$W_i^{**}=W_i\setminus f_i^{-1}(0)=W_i^*\setminus f_i^{-1}(0).$$
Since $f_i^{-1}(0)$ is a complex analytic hypersurface in the connected manifold $W_i$, its complement $W_i^{**}$ is connected. Moreover, the inclusions
$W_i^{**}\hookrightarrow W_i$ and $W_i^{**}\hookrightarrow W_i^*$ both induce surjections on fundamental groups. By compactness of $W_i$ and transversality along its boundary, the restriction
$$f_i:W_i^{**}\to \Delta^*$$
is a proper smooth submersion. By Erhesmann's fibration theorem
\begin{equation}\label{eq: fibration with center removed}
f_i|_{W_i^{**}}:W_i^{**}\to \Delta^*,
\end{equation}
is a locally trivial fibration. 

Consider the following long exact sequence on homotopy groups of the fibration (\ref{eq: fibration with center removed}):
\begin{equation}\label{eq: ses of homotopy}
\begin{aligned}
\pi_1(f_i^{-1}(v)\cap W_i)
\to \pi_1(W_i^{**})
\xrightarrow{(f_i)_*}\pi_1(\Delta^*)\to \pi_0(f_i^{-1}(v)\cap W_i)\to \pi_0(W_i^{**})=1
\end{aligned}
\end{equation}
with $0<|v|<\delta$.

By Case \ref{it: condition B}, choose $y\in f_i^{-1}(0)\cap T_i$, we have $y\notin\Sigma_i$ and $f_i^{-1}(0)$ meets $\partial U_i$ transversely at $y$. Hence we can choose a regular point $x\in f_i^{-1}(0)\cap\interior(W_i)\setminus\Sigma$ sufficiently close to $y$; in particular, $x\in W_i^*$. By the implicit function theorem, there exists a small holomorphic disc $\sigma: {\Delta(0,\delta_1)}\to W_i^*$ with $0<\delta_1< \delta$ such that $f_i\circ\sigma=\id$, thus $\im(\sigma)\cap f_i^{-1}(0)=\{x\}$. Let $\alpha$ be the image under $\sigma_*$ of a sufficiently small circle $\{|v|=\delta_2\}$ with $\delta_2<\delta_1$ in this disc, then $\alpha$ is a loop in $W_i^{**}$ such that $(f_i)_*[\alpha]$ generates $\pi_1(\Delta^*)\cong\Z$, and $\alpha$ is null-homotopic in $W_i^*$. In particular, $(f_i)_*:\pi_1(W_i^{**})\to\pi_1(\Delta^*)$ is surjective, then \cref{eq: ses of homotopy} implies that 
\begin{equation}\label{eq: nearby fiber connectedness}
\text{every fiber $f_i^{-1}(v)\cap W_i$, $0<|v|<\delta$,
is connected.}
\end{equation}
Besides, for any $[\beta]\in\pi_1(W_i^{**})$, let 
$n=(f_i)_*[\beta]\in \pi_1(\Delta^*)$, then $[\beta][\alpha]^{-n}
\in\ker(f_i)_*$, and hence belongs to the image of
$\pi_1(f_i^{-1}(v)\cap W_i)$. Since $\alpha$ is null-homotopic in $W_i^*$ and $\pi_1(W_i^{**})\to\pi_1(W_i^*)$ is surjective, the inclusion $f_i^{-1}(v)\cap W_i\hookrightarrow W_i^*$ induces a surjection on fundamental groups:
\begin{equation}\label{eq: surjection on pi 1 of nearby fiber}
\pi_1(f_i^{-1}(v)\cap W_i)\twoheadrightarrow \pi_1(W_i^*).
\end{equation}
\end{proof}

\begin{remark}
Similarly,
%due to the same argument using Erhesmann's fibration theorem, 
we have a locally-trivial fibration
\begin{equation}\label{eq: fibration on the boundary T}
f_i|_{T_i}: T_i\to \Delta(0;\delta).
\end{equation} 
\end{remark}

Next, consider an irreducible component $Q_i\subset Q''$. Set
$$Q_i^\circ\coloneqq Q_i\setminus P.$$
Since $\xi|_{Q_i^\circ}\not\equiv0$ and $\xi|_{Q_i^\circ}$ extends meromorphically across the finite set $Q_i\cap P$, its zero set on $Q_i^\circ$ is finite. Let $D_i$ denote the finite set of points $p\in Q_i^\circ$ such that
$$\xi(p)\ne0,\quad \xi|_{T_pQ_i}=0.$$
Away from $D_i$ and the Type I critical neighborhoods, we may choose the vector field contructed in \cref{subsec: Construction of a flow on the infinite cyclic covering} to be tangent to $Q_i$ (or its lift) and to satisfy $\Re\xi(v)=1$, so that its flow preserves $Q_i$ and hence also its complement. 

We therefore need additional local neighborhoods around the points of $D_i$, as well as around the finitely many transverse intersection points of distinct components of $Q''$. 
\begin{enumerate}
    \item[$(\alpha)$] If $p\in Q_i\cap Q_j$ and $\xi$ vanishes on both $T_pQ_i$ and $T_pQ_j$, then
    %, since the two tangent directions span $T_pY^\circ$, we have $\xi(p)=0$. Hence 
    $p\in Z(\xi)$, and this point is already contained in a type I critical neighborhood. 
    \item[($\beta$)] If $\xi$ vanishes on only one of the two tangent directions, then we can take a type II critical neighborhood.
    \item[($\gamma$)]\label{it: intersection of Q both nonvanishing} If neither tangential restriction vanishes, we still include $p$ among the type II critical neighborhoods.
\end{enumerate}
The argument for case $(\gamma)$ is given in
\cref{rmk: type II}; in this case $f_l|_{Q_i}$ and $f_l|_{Q_j}$ are both submersions near $p$, so over the punctured disc the intersections with $Q_i\cup Q_j$ vary as two distinct unramified branches. The same local $1$-connectedness argument as for type II neighborhoods applies. The case $(\beta)$ is analogous, with $\fkW_l^*=\fkW_l\setminus(Q_i\cup Q_j)$.

At each point $p_l\in D_i$, choose a sufficiently small simply connected neighborhood $V_l'\subset Y^\circ$. Since $\xi$ is closed, the holomorphic Poincar\'e lemma gives a holomorphic function $f_l$ on $V_l'$ such that
$$df_l=\xi,\quad f_l(p_l)=0.$$
Since $\xi(p_l)\neq0$, after shrinking $V_l'$ we may assume that $\xi$ is nowhere vanishing on $V_l'$. Thus $f_l$ is a holomorphic submersion. Choose local holomorphic coordinates $(z_1,z_2)$ centered at $p_l$ such that $f_l=z_1$. Since $Q_i$ is smooth, the curve $Q_i$ can be written locally as 
$$Q_i=\{z_1=\phi(z_2)\},$$
where $\phi(0)=\phi'(0)=0$ and $\phi$ is not identically zero. After a further holomorphic change of the coordinate $z_2$, leaving $z_1$ unchanged, we may therefore write $\phi(z_2)=z_2^{m_l}$ for some integer $m_l\ge2$. Thus, locally,
$$Q_i=\{z_1=z_2^{m_l}\}.$$

Choose $\varepsilon_l>0$ sufficiently small and then choose $0<\delta_l<\varepsilon_l^{m_l}$. Define
$$\fkW_l=\left\{(z_1,z_2):|z_1|\le\delta_l,\ |z_2|\le\varepsilon_l\right\}.$$
such that $\fkW_l$ is a smooth compact manifold with smooth corners. Its boundary consists of
$$\fkT_l=\left\{|z_1|\le\delta_l,\ |z_2|=\varepsilon_l\right\},\quad \fkR_l=\left\{|z_1|=\delta_l,\ |z_2|\le\varepsilon_l\right\}.$$
By the choice $\delta_l<\varepsilon_l^{m_l}$, we have
$$Q_i\cap\fkT_l=\varnothing.$$
Thus $Q_i$ exits $\fkW_l$ only through $\fkR_l$. After shrinking these neighborhoods if necessary, we may assume that the $\fkW_l$ are pairwise disjoint. We call the resulting neighborhoods $\fkW_l$ the \textit{critical neighborhoods of Type II}, see \cref{fig: typeI-typeII-critical-neighborhoods}.

For a critical neighborhood of Type II, we have a local $1$-connectedness property which is considerably simpler than \cref{prop: local 1-connectedness for type I}. The essential point is that $\xi$ itself is nowhere vanishing on $\fkW_l$. Set
$$\fkW_l^{*}=\fkW_l\setminus Q_i\quad \text{and}\quad \fkW_l^{**}=f_l^{-1}(\Delta^*(0;\delta_l))\cap(\fkW_l\setminus Q_i).$$
\begin{proposition}\label{prop: local 1-connectedness for type II}
With the above notations, for every sufficiently small nonzero $v\in\Delta(0;\delta_l)$, the pair
$$\left(\fkW_l^*,\,f_l^{-1}(v)\cap \fkW_l^*\right)$$
is $1$-connected.
\end{proposition}

\begin{proof}
By the local description above, for every $0<|v|<\delta_l$, the fiber $f_l^{-1}(v)\cap \fkW_l^*$ is a disc with the $m_l$ points satisfying $z_2^{m_l}=v$ removed. In particular, it is connected.

Since the map $z_2\mapsto z_2^{m_l}$ is unramified over $\Delta^*(0;\delta_l)$, the $m_l$ deleted points vary locally as $m_l$ disjoint holomorphic branches. Hence
$$f_l:\fkW_l^{**}\to\Delta^*(0;\delta_l)$$
is a locally trivial fibration. Its homotopy exact sequence gives
$$\pi_1(f_l^{-1}(v)\cap \fkW_l^*)\to\pi_1(\fkW_l^{**})\to\pi_1(\Delta^*(0;\delta_l))\to1.$$
Moreover, the inclusion $\fkW_l^{**}\hookrightarrow\fkW_l^*$ induces a surjection on fundamental groups, since $\fkW_l^{**}$ is obtained by removing the central complex hypersurface from $\fkW_l^*$.

Choose $a\in\C$ such that $\delta_l^{1/m_l}<|a|<\varepsilon_l$. Then $s(v)=(v,a)$ defines a section of $f_l$ over the disc $\Delta(0;\delta_l)$ whose image is disjoint from $Q_i$. For a loop $\gamma$ representing a generator of $\pi_1(\Delta^*(0;\delta_l))$, the loop $s(\gamma)$ lies in $\fkW_l^{**}$ and is null-homotopic in $\fkW_l^*$.
%, since it bounds the disc $s(\Delta(0;\delta_l))$.

Therefore, for any class in $\pi_1(\fkW^*)$, we may deform it to a class in $\pi_1(\fkW_l^{**})$, and then multiply by a suitable power of $[s(\gamma)]$ so that the resulting class lies in the kernel of
$$\pi_1(\fkW_l^{**})\to\pi_1(\Delta^*(0;\delta_l)).$$
thus lies in the image of $\pi_1(f_l^{-1}(v)\cap(\fkW^*))$. Since $[s(\gamma)]$ is trivial in $\pi_1(\fkW^*)$, it follows that
$$\pi_1(f_l^{-1}(v)\cap(\fkW^*))\twoheadrightarrow\pi_1(\fkW^*).$$
Together with the connectedness of the nearby fiber, this proves that the pair is $1$-connected.
\end{proof}
\begin{remark}\label{rmk: type II}
For a type II critical neighborhood, Ehresmann's fibration theorem is not directly applicable, since after deleting $Q_i$ the resulting map need not be proper. Instead, the explicit local description above, in which the $m_l$ branches are unramified away from $p_l$, gives the required locally trivial fibration.

In the case $(\gamma)$, let $p\in Q_i\cap Q_j\subset Q''$. Since neither $\xi|_{T_pQ_i}$ nor $\xi|_{T_pQ_j}$ vanishes, we may choose local coordinates $(z_1,z_2)$ centered at $p$ such that $f=z_1$ and
$$Q_i=\{z_2=\phi_i(z_1)\},\quad Q_j=\{z_2=\phi_j(z_1)\},$$
where $\phi_i(v)\ne\phi_j(v)$ for $v\ne0$ sufficiently small. Taking
$$\fkW^{*}=\fkW\setminus(Q_i\cup Q_j),$$
the fibers of $f:\fkW^*\to\Delta$ away from $0$ are discs with two points removed. The same section and fibration argument as above therefore gives the required local $1$-connectedness, so this case may also be treated as a Type II critical neighborhood. $(\beta)$ is proved similarly, thus we shall not repeat.
\end{remark}

The two types of critical neighborhoods may intersect, as illustrated on the left in \cref{fig: typeI-typeII-critical-neighborhoods}. Recall that the $Q'$-part is already contained in $\Sigma$. Thus, any non-polar component of $Q$ shown intersecting $\Sigma$ in the picture on the left belongs to $Q''$. We use the notation $Q_i''$ to emphasize this in the picture. The same convention applies to the picture on the right.

This kind of mixed situation can be remedied as follows. Fix $\Sigma_i$. There are only finitely many irreducible components of $Q''$ intersecting $\Sigma_i$; relabel them as $Q_1^i,\ldots,Q_{J_i}^i$. Since $Q_j^i\subset Q''$, we have $\xi|_{Q_j^i}\not\equiv0$. Hence $Q_j^i$ is not contained in $\Sigma_i$, and therefore $Q_j^i\cap\Sigma_i$ is a finite set. Define the union of these $Q^i_j$ as $\mathcal{Q}_i''$.

We can moreover arrange the parameter $\delta$ in the construction of the type I neighborhood so that each $Q_j^i$ exits through $R_i$, and not through $T_i$. Indeed, since $df_i=\xi$ and $\xi|_{Q_j^i}\not\equiv0$, the restriction $f_i|_{Q_j^i}$ is nonconstant. Hence its zero set is discrete. By shrinking the neighborhood of $\Sigma_i$ and choosing $\varepsilon$ sufficiently small, we may assume that
$$Q_j^i\cap\partial U_i\cap f_i^{-1}(0)=\varnothing$$
for every $1\le j\le J_i$. For each $j$ such that $Q_j^i\cap\partial U_i\ne\varnothing$, set
$$c_j=\min_{z\in Q_j^i\cap\partial U_i}|f_i(z)|>0.$$
We then further require, in the construction of the type I neighborhood, that $\delta<\min_j c_j$, where the minimum is taken over those $j$ for which $Q_j^i\cap\partial U_i\ne\varnothing$. It follows that $Q_j^i\cap T_i=\varnothing$ for every $j$. Thus the components $Q_j^i$ can meet the boundary of $W_i$ only through $R_i$, and not through $T_i$. In the mixed case, we always choose the type I neighborhoods with this additional property.

After decreasing $\delta$ further if necessary, we may also assume that $f_i|_{Q_j^i}$ has no critical point over $\Delta^*(0;\delta)$ for every $j$. Thus, over $v$ in the punctured disc, the finite set $f_i^{-1}(v)\cap Q_i$ varies locally as a finite collection of unramified branches, as in the case of type II critical neighborhoods.

Set
$$W_i^\dagger=W_i\setminus(Q'\cup \calq_i'') \quad \text{and}\quad W_i^{\dagger\dagger}=W_i^\dagger\cap f_i^{-1}(\Delta^*(0;\delta)).$$
Since $Q'\cap W_i\subset f_i^{-1}(0)$, we have
$$W_i^{\dagger\dagger}=W_i^{**}\setminus \calq_i''.$$
The map
$$f_i:W_i^{**}\to\Delta^*(0;\delta)$$
is a locally trivial fibration by \cref{prop: local 1-connectedness for type I}. Since the intersections with the components of $\calq_i''$ vary locally as unramified branches, the restriction
$$f_i:W_i^{\dagger\dagger}\to\Delta^*(0;\delta)$$
is again a locally trivial fibration. Its fiber over $v\ne0$ is $F_v^\dagger\coloneqq f_i^{-1}(v)\cap W_i^\dagger$. It is obtained from the connected smooth complex curve $f_i^{-1}(v)\cap W_i$ by removing finitely many points, namely $f_i^{-1}(v)\cap\calq_i''$, and is therefore connected.

Consider the associated homotopy exact sequence:
$$\pi_1(F_v^\dagger)\to\pi_1(W_i^{\dagger\dagger})\to\pi_1(\Delta^*(0;\delta))\to\pi_0(F_v^\dagger)=1.$$
By the construction above, both $Q'$ and the components $Q_j^i\subset \calq_i''$ are disjoint from $T_i$. Hence, as in the proof of \cref{prop: local 1-connectedness for type I}, we may choose a regular point
$$x\in f_i^{-1}(0)\cap\left(\interior(W_i)\setminus(Q'\cup \calq_i'')\right)$$
sufficiently close to $T_i$. Up to shrinking $\delta$, there is a holomorphic section $\sigma:\Delta(0;\delta)\to W_i^\dagger$ with $f_i\circ\sigma=\id$. Hence the loop obtained by applying $\sigma$ to a generator of $\pi_1(\Delta^*(0;\delta))$ is null-homotopic in $W_i^\dagger$. Finally, the inclusion $W_i^{\dagger\dagger}\hookrightarrow W_i^\dagger$ induces a surjection on fundamental groups.
%, since its complement is the remaining part of the central complex hypersurface $f_i^{-1}(0)$, which has real codimension two. 
The same argument as in the proof of \cref{prop: local 1-connectedness for type I} therefore gives
$$\pi_1(F_v^\dagger)\twoheadrightarrow\pi_1(W_i^\dagger).$$
Thus the same local $1$-connectedness conclusion holds in the mixed case. In conclusion, we have proved
\begin{proposition}\label{prop: local 1-connectedness for mixed}
In the mixed case, we still have the pair
$$(W_i^{\dagger},F_v^\dagger)$$
is 1-connected, for any $v\ne 0$ in $\Delta(0;\delta)$.
\end{proposition}

\begin{remark}
There is a useful analogy with the Lefschetz-type argument of Eyssidieux in \cite[\S 5.1]{Eyssidieux04}. There, after decomposing the degeneracy locus of a system of holomorphic $1$-forms into irreducible components, one restricts the system to each such component and imposes a nondegeneracy condition on the restriction. This is quite similar to the mixed case considered here.
\end{remark}

We deduce one more fact from the above $1$-connectedness results that will be used later. The critical neighborhood can be of type I, type II, or mixed type. After decreasing $\delta>0$ if necessary, we assume that the same $\delta$ works for all type I, type II, and mixed-type critical neighborhoods.

\begin{proposition}\label{prop: local sup level 1 connected}
For any $-\delta<\eta<\delta$, 
\begin{enumerate}
\item For a type I critical neighborhood $W_i$, denote
$$S^*_\eta=W_i^*\cap \{\eta\le \Re f_i\},$$
then $S_\eta^*$ is connected and $\pi_1(S_\eta^*)\to \pi_1(W_i^*)$ is surjective.
\item For a type II critical neighborhood $\fkW_l$, denote
$$\fkS^*_\eta=\fkW_l^*\cap \{\eta\le \Re f_l\},$$
then $\fkS_\eta^*$ is connected and $\pi_1(\fkS_\eta^*)\to \pi_1(\fkW_l^*)$ is surjective.
\item For a mixed-type critical neighborhood $W_i$, denote
$$S^\dagger_\eta=W_i^\dagger\cap \{\eta\le \Re f_i\},$$
then $S_\eta^\dagger$ is connected and $\pi_1(S_\eta^\dagger)\to \pi_1(W_i^\dagger)$ is surjective.
\end{enumerate}
A similar conclusion holds if we consider intersection with $\{\Re f_i\le \eta\}$ instead. 

More generally, let $-\delta<\eta_1 \le  \eta_2<\delta$, with either $\eta_1<\eta_2$ or $\eta_1=\eta_2\ne0$. Then the same conclusions hold for the corresponding bands $W_i^*\cap\{\eta_1\le\Re f_i\le\eta_2\}$, $\fkW_l^* \cap\{\eta_1\le\Re f_l\le\eta_2\}$, and  $W_i^\dagger\cap\{\eta_1\le\Re f_i\le\eta_2\}$.
\end{proposition}
\begin{proof}
For (a)-(c), we abuse $W_i$ for critical neighborhood of type I or II. Choose $v\in\Delta(0;\delta)$ such that $v\ne0$ and $\Re v\ge\eta$. Then
$$f_i^{-1}(v)\cap W_i^*\subset S_\eta^*\subset W_i^*.$$
By \cref{prop: local 1-connectedness for type I,prop: local 1-connectedness for type II},
$\pi_1(f_i^{-1}(v)\cap W_i^*)\to\pi_1(W_i^*)$
is surjective. Hence the factorization
$$\pi_1(f_i^{-1}(v)\cap W_i^*)\to\pi_1(S_\eta^*)\to\pi_1(W_i^*)$$
implies that $\pi_1(S_\eta^*)\to\pi_1(W_i^*)$ is surjective.
%In the mixed case (c), choose again $v\ne0$ with $\Re v\ge\eta$. Then $F_v^\dagger\subset S_\eta^\dagger\subset W_i^\dagger$. By \cref{prop: local 1-connectedness for mixed}, $\pi_1(F_v^\dagger)\to\pi_1(W_i^\dagger)$ is surjective. Therefore the factorization $$\pi_1(F_v^\dagger)\to\pi_1(S_\eta^\dagger)\to\pi_1(W_i^\dagger)$$ shows that $\pi_1(S_\eta^\dagger)\to\pi_1(W_i^\dagger)$ is surjective.

For the connectedness statements, in each of the three cases, the corresponding punctured superlevel set is a locally trivial fibration with connected fibers over
$\{z\in\Delta^*(0;\delta):\Re z\ge\eta\}$, which is connected. Hence the punctured superlevel set is connected. If $\eta>0$, this punctured superlevel set is already the desired superlevel set. If $\eta\le0$, the desired superlevel set is its closure in $W_i^*$, $\fkW_l^*$, or $W_i^\dagger$, respectively, obtained by adjoining the remaining part of the central fiber. Hence $S_\eta^*$, $\fkS_\eta^*$, and $S_\eta^\dagger$ are connected.

The sublevel and band assertions follow by almost the same argument; which we shall not repeat.

%For the band case, set $$B_{\eta_1,\eta_2}=\{z\in\Delta^*(0;\delta):\eta_1\le\Re z\le\eta_2\}.$$  The assumptions on $\eta_1,\eta_2$ ensure that $B$ is nonempty and connected. The same fibration argument therefore shows that the corresponding punctured band is connected. By holomorphic openness and the boundary transversality established above, its closure in the corresponding deleted neighborhood is the full band, which is consequently connected.

%Choose $v\in B_{\eta_1,\eta_2}$. The nearby fiber over $v$ is contained in the band, and its inclusion into the corresponding deleted neighborhood induces a surjection on fundamental groups by the preceding local $1$-connectedness results. Factoring this inclusion through the band proves the required surjectivity on $\pi_1$.
\end{proof}
%Consequently, a path in $W_i^*$, $\fkW_l^*$, or $W_i^\dagger$, respectively, joining two points $x,y$ in $S_\eta^*$, $\fkS_\eta^*$, or $S_\eta^\dagger$, respectively, can be deformed, with $x$ and $y$ fixed, into a path contained in the corresponding superlevel set. 
From now on, we include the mixed-type critical neighborhoods among the type I critical neighborhoods. Thus, unless otherwise specified, a type I critical neighborhood may refer either to an ordinary type I neighborhood $W_i^*$ or to a mixed-type neighborhood $W_i^\dagger$.

\subsection{Construction of a flow on the infinite cyclic covering}\label{subsec: Construction of a flow on the infinite cyclic covering}

In this section, we introduce and modify the flow used by Zhou in \cite{Zhou} to deform paths whose endpoints lie in a level set into that level set. Under the assumption $\dim\psi(X)\ge2$, we have shown that case \ref{it: condition A} cannot occur, so every critical neighborhood of type I (if any) satisfies case \ref{it: condition B}. We will show that this leads to \textup{(ii)} of \cref{thm: Lefschetz alternative rational case}.

The main point is to control the flow away from the critical neighborhoods so that it remains in $Y$, equivalently, so that it does not cross the boundary $P\cup Q$. Near the polar divisor $P$, this will follow from the boundary behavior established in \cref{subsec: Choice of logarithmic holomorphic one-form}. Along $Q''$, away from the critical neighborhoods, we choose the flow tangent to $Q''$, so that it preserves its complement. Inside critical neighborhoods of type I, type II, or mixed type, such a flow can still be constructed, but it may fail to have the desired increasing property for the primitive function, then we replace the relevant portions of the path using the local $1$-connectedness results established in \cref{prop: local 1-connectedness for type I,prop: local 1-connectedness for type II,prop: local 1-connectedness for mixed,prop: local sup level 1 connected}.

We now lift all the objects constructed in \cref{subsec: Construction of two types of critical neighborhoods} to the regular $\Z $-cover
$$
p:\widehat M\to M,
$$
corresponding to $\chi$, which can be defined explicitly as
\begin{equation}\label{eq: def of covering}
\widehat M\coloneqq M\times_{\bbS^1}\R
=\{(x,t)\in M\times\R:\Phi(x)=e^{2\pi\sqrt{-1}t}\},
\end{equation}
where $M$ is the real-oriented blowup of $\bar Y$ along the polar divisor $P$. To simplify notation, we also use $P$ to denote the polar boundary $\beta^{-1}(P)\subset M$. We write $T$ for a fixed generator of the deck transformation group of $\widehat M$. For any connected subset $B\subset M$, we use $\widehat B$ to denote one chosen connected component of $p^{-1}(B)$. Its other connected lifts, when distinct, are obtained from $\widehat B$ by applying suitable powers of $T$. Thus, throughout the following discussion, a hat indicates a chosen connected component of the inverse image under $p$.

We first consider the critical neighborhoods of the two types. On every critical neighborhood under consideration, the restriction of $\xi$ admits a single-valued holomorphic primitive: namely, $\xi=df_i$ on a type I neighborhood $W_i$, and $\xi=df_l$ on a type II neighborhood $\fkW_l$. Hence $\Re\xi$ is exact on each such neighborhood, so $\int_\gamma\Re\xi=0$ for every loop $\gamma$ contained in it. Therefore the character $\chi$ vanishes on the fundamental group of every critical neighborhood. It follows that the inverse image of each critical neighborhood in the $\Z $-cover is a disjoint union of copies of that neighborhood. More precisely,
$$p^{-1}(W_i)=\bigsqcup_{n\in\Z }\widehat W_{i,n} =\bigsqcup_{n\in\Z }T^n\widehat W_{i,0},$$
where each $\widehat W_{i,n}$ is mapped homeomorphically onto $W_i$; similarly, for type II,
$$p^{-1}(\fkW_l)=\bigsqcup_{n\in\Z}\widehat{\fkW}_{l,n}=\bigsqcup_{n\in\Z }T^n\widehat{\fkW}_{l,0},$$
where each $\widehat{\fkW}_{l,n}$ is mapped homeomorphically onto $\fkW_l$.

Finally, we consider the lifting of the polar boundary $P$. Recall that, near a component $P_i\subset P$, the logarithmic form $\bar\psi^*\omega$ has nonzero residue along $P_i$. More precisely, on a punctured tubular neighborhood of $P_i$, we have $\xi=\bar\psi^*\omega$,
and for a small positively oriented loop $\gamma$ around $P_i$,
$$\int_\gamma\Re\xi=m_i\ne0.$$

The lift of a simply connected local component of $P_i$ (along its tubular neighborhood $N_i$) contains $m_i$-many components by standard theory of covering space; indeed, we have $\int \Re\xi:\pi_1(N_i)\to \Z$ has image $m_i\Z$, and its lift to $\widehat{M}$ has $|\Z/m_i\Z|=m_i$ many connected components, labeled as $\widehat P_{i}^{j}, 1\le j\le m_i$. We emphasize that this is a local picture. For a tubular neighborhood of the entire component $P_i$, loops in $P_i$ may contribute additional $\xi$-periods. Thus the image of $\int\Re\xi:\pi_1(N_i)\to\Z $
may be $d_i\Z $ for some $d_i\mid m_i$. Correspondingly, the inverse image of the global tubular neighborhood, as well as $P_i$, has $d_i$-many connected components.

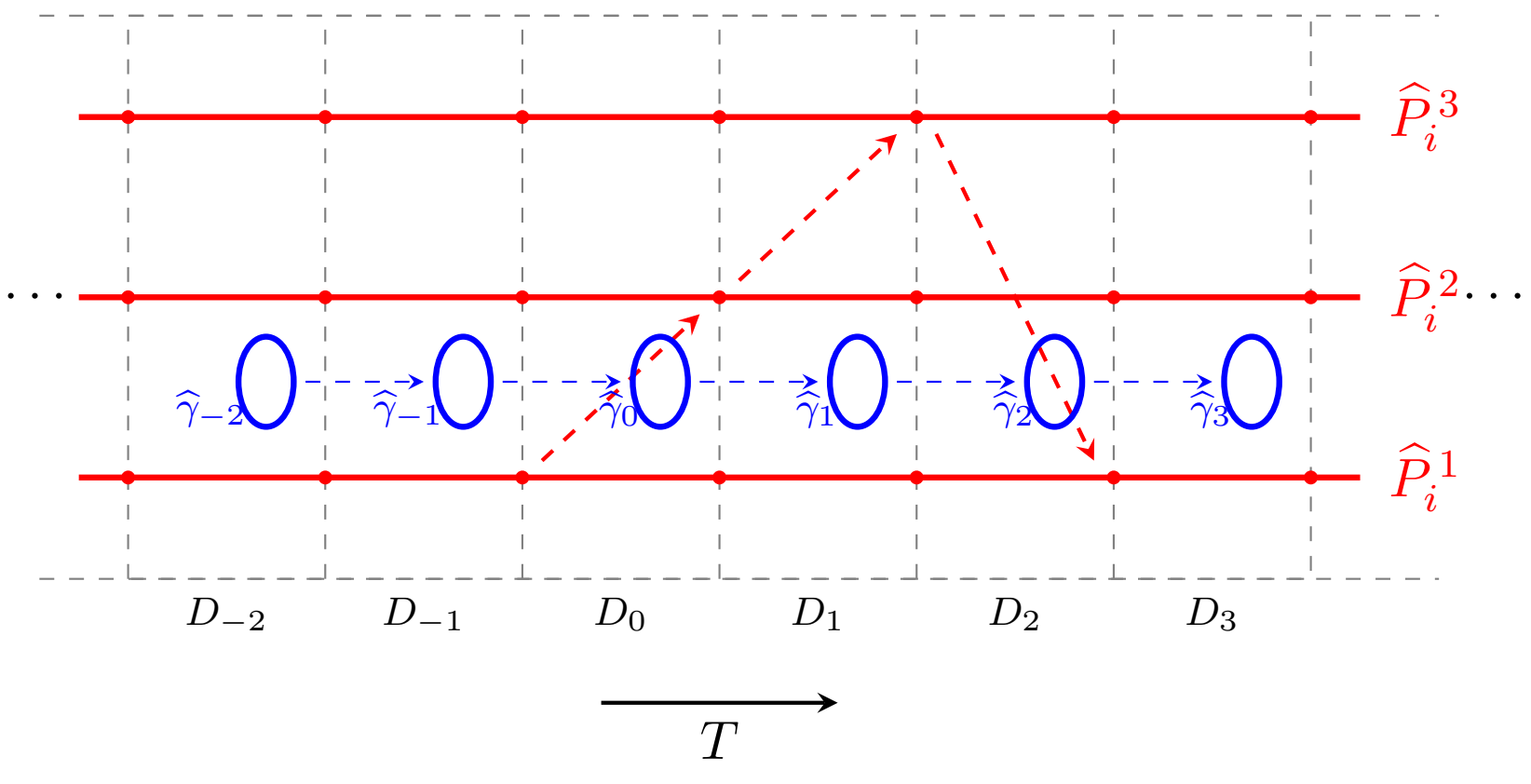
\begin{figure}[ht!]
    \centering
\begin{tikzpicture}[x=1.4cm,y=4cm,>=stealth]

\def\W{6}
\def\H{1}

%------------------------------------------------
% Outer rectangle: flattened cylinder
%------------------------------------------------
\draw[dashed,gray] (0,0) rectangle (\W,\H);

% Top/bottom identification
%\node at (\W/2,\H+0.16)   {\small Top and bottom edges are identified};

%------------------------------------------------
% Fundamental domains
%------------------------------------------------
\foreach \x in {1,...,5}{
    \draw[dashed,gray] (\x,0)--(\x,\H);
}

% Labels for fundamental domains
\foreach \x/\k in {
    0.5/{-2},
    1.5/{-1},
    2.5/{0},
    3.5/{1},
    4.5/{2},
    5.5/{3}
}{
    \node[below] at (\x,0)
        {\scriptsize $D_{\k}$};
}

% Indication of infinitely many fundamental domains
\node at (-0.45,0.5) {$\cdots$};
\node at (\W+0.95,0.5) {$\cdots$};

% Extend top and bottom edges slightly
\draw[dashed,gray] (-0.45,1)--(\W+0.65,1);
\draw[dashed,gray] (-0.45,0)--(\W+0.65,0);

%------------------------------------------------
% Three local lifts of the polar component
% Example m_i = 3
%------------------------------------------------
\foreach \y/\lab in {
    0.18/{\widehat P_i^{\,1}},
    0.50/{\widehat P_i^{\,2}},
    0.82/{\widehat P_i^{\,3}}
}{
    \draw[very thick,red]
        (-0.25,\y)--(\W+0.25,\y);

    \node[right,red]
        at (\W+0.30,\y)
        {$\lab$};
}

%------------------------------------------------
% Red points showing the action of T on the
% lifted polar components
%------------------------------------------------
\foreach \x in {0,...,6}{
    \foreach \y in {0.18,0.50,0.82}{
        \fill[red] (\x,\y) circle (1.4pt);
    }
}

% Deck transformation permutes the three
% local lifts of the polar component
\draw[->,red,dashed,thick]
    (2.10,0.21) -- (2.90,0.47);

\draw[->,red,dashed,thick]
    (3.10,0.53) -- (3.90,0.79);

\draw[->,red,dashed,thick]
    (4.10,0.79) -- (4.90,0.21);

%------------------------------------------------
% Blue elliptic objects:
% lifts of zero-period loops / models for
% Type I and Type II critical neighborhoods.
%------------------------------------------------
\foreach \x/\k in {
    0.70/{-2},
    1.70/{-1},
    2.70/{0},
    3.70/{1},
    4.70/{2},
    5.70/{3}
}{
    % ellipse
    \draw[very thick,blue] (\x,0.35) ellipse [x radius=0.14, y radius=0.08];

    % label on the left of each ellipse
    \node[blue,left]
        at (\x,0.3)
        {\scriptsize $\widehat{\gamma}_{\k}$};
}

%------------------------------------------------
% Blue dashed arrows indicating the T-action
% on the zero-period lifts
%------------------------------------------------
\foreach \x in {0.90,1.90,2.90,3.90,4.90}{
    \draw[->,blue,dashed]
        (\x,0.35)--(\x+0.60,0.35);
}

%------------------------------------------------
% Main deck transformation arrow
%------------------------------------------------
\draw[->,thick]
    (2.40,-0.22)--(3.60,-0.22);

\node[below]
    at (3.00,-0.22)
    {$T$};
\end{tikzpicture}
\caption{Lifting of polar boundary}
\label{fig: zcover-cylinder}
\end{figure}
\noindent The \cref{fig: zcover-cylinder} shows the case when $d_i=m_i=3$, with red dashed arrows meaning the orbit of a point on $\widehat{P}_i$ under action of $T$; in particular, $T$ does not need to preserve each connected component. Those $\widehat{\gamma}_i$ corresponds to loops with $\xi$-period 0, and is a model of the lift of two types critical neighborhoods; each $D_i$ means the fundamental domain.

\begin{remark}\label{rmk: on rationality and why Z cover}
This explains why we work on the $\Z$-covering rather than on the universal covering. Our argument only needs to control the real part of a holomorphic $1$-form $\xi$ corresponding to $\chi$. Passing to the $\ker\chi$-cover kills exactly the periods of $\Re\xi$ detected by $\chi$, so that
$$du=p^*\Re\xi$$
for a globally defined primitive function $u$, while the deck group records precisely the translations detected by $\chi$. Passing instead to the universal covering would introduce the full fundamental group as deck group and substantially more information than is needed for the one-dimensional dynamic behavior determined by $\xi$.
\end{remark}
Now we take care of the $1$-connectedness of the level sets of the primitive function
$$u:\widehat M\to\R,$$
defined by
$$du=p^*\Re\xi,\quad u(Tx)=u(x)+1,\quad \forall x\in\widehat M.$$
With $\Phi$ defined downstairs in \cref{eq: Phi downstairs} and $\widehat M$ defined in \cref{eq: def of covering}, we have
$$e^{2\pi\sqrt{-1}u}=\Phi\circ p.$$

Recall that $Y$ is a dense open subset of $M$. By abuse of notation, we continue to denote the oriented boundary $\beta^{-1}(P)$ by $P$, and also denote the proper transform of $Q$ in $M$ by $Q$. Then
$$M\setminus Y=P\cup Q=\partial M\cup Q.$$
Write $\widehat Y=p^{-1}(Y)$ and, for $t\in\R$ regular value of $u$, put
$$L_t\coloneqq u^{-1}(t)\cap\widehat Y.$$
We will show that any smooth path $\gamma:[0,1]\to\widehat Y$ whose endpoints $p_0=\gamma(0)$, $p_1=\gamma(1)$ lie in $L_t$ can be deformed relative to its endpoints into a path entirely contained in $L_t$. This shows that $L_t$ is connected and that $\pi_1(L_t)\twoheadrightarrow\pi_1(\widehat Y)$
is surjective.

First, after perturbing $\gamma$ slightly relative to its endpoints, we may assume that it meets $L_t$ transversely, including at its
endpoint germs. Then $\gamma$ intersects $L_t$ at finitely many points, and we can divide $\gamma$ into finitely many subpaths accordingly, whose endpoints lie in $L_t$ and whose interiors are disjoint from $L_t$; in other words, each subpath lies entirely either below or above $L_t$, except for its two endpoints.

We may therefore treat these finitely many subpaths separately. It suffices to consider the case in which
$\gamma(0),\gamma(1)\in L_t$, and $u(\gamma(a))<t$
for all $a\in(0,1)$; the case $u(\gamma(a))>t$ is treated in the same way using the negative flow. For $\gamma$ below the regular level set, let
\begin{equation}\label{eq: def of r0}
r_0=\min_{a\in[0,1]}u(\gamma(a))<t.
\end{equation}

Next, we construct a vector field $\widehat{\bfv}$ which is the lift of a suitable smooth vector field $\bfv$ on $M$. The flow defined by $\widehat{\bfv}$ will push $\gamma$ upwards until it lies entirely inside the regular level set $L_t$.

For the critical neighborhoods of both types, let $\delta>0$ be the uniform constant chosen above. Choose $c>0$ sufficiently small such that $5c<\delta$.
We moreover choose $c$ so that, for every type I critical neighborhood, all values $z$ with $0<|z|\le5c$ are regular values of $f_i$, and, in the mixed case, the restrictions of $f_i$ to the components of $\calq_i''$ have no critical points over this punctured disc. For a type II critical neighborhood, $f_l$ is already a submersion, so the corresponding regularity is automatic.

For each type I critical neighborhood $W_i$, choose relatively compact neighborhoods of $\Sigma_i$:
$$W_i'''\Subset W_i''\Subset W_i\cap f_i^{-1}(\overline{\Delta(0;c)})\Subset W_i'\Subset W_i,$$
where $W_i'=W_i\cap f_i^{-1}(\overline{\Delta(0;3c)})$. Define
$$T_i'=T_i\cap W_i'\quad \text{
and}\quad R_i'=W_i\cap\{|f_i|=3c\}.$$

Similarly, for each type II critical neighborhood $\fkW_l$, choose relatively compact neighborhoods of $p_l$:
$$\fkW_l'''\Subset\fkW_l''\Subset\fkW_l\cap f_l^{-1}(\overline{\Delta(0;c)})\Subset\fkW_l'\Subset\fkW_l,$$
where $\fkW_l'=\fkW_l\cap f_l^{-1}(\overline{\Delta(0;3c)})$. Define
$$\fkT_l'=\fkT_l\cap\fkW_l'\quad\text{and}\quad\fkR_l'=\fkW_l\cap\{|f_l|=3c\}.$$

In the mixed case, which from now on is included among the type I critical neighborhoods, the same construction is used. By our previous choice of the type I neighborhoods, every component of $\calq_i''$ is disjoint from $T_i$ and exits $W_i$ only through the $R_i$-side. The same property therefore holds for the smaller neighborhood $W_i'$ with $R_i'$ in place of $R_i$.

Set $\alpha=\Re\xi$. We choose a smooth vector field $\bfv$ on $M$ such that:
\begin{enumerate}
\item \label{it: v vanishes neaby singular set} $\bfv$ vanishes in the neighborhoods $W_i'''$ and $\fkW_l'''$ for each $i$ and $l$.

\item \label{it: out of critical increase} Outside
$\bigsqcup_{i=1}^m W_i''\cup\bigsqcup_{l=1}^n\fkW_l''$,
we have
$$d \arg(\Phi)(\bfv)=2\pi\Longleftrightarrow \alpha(\bfv)=1.$$
\item \label{it: parallel to P} $\bfv$ is tangent to every boundary component $P_i$ of the real-oriented blowup $M$. For example, in the local coordinates used in \cref{subsec: Reduction to surfaces}, it may be chosen to be $\frac{2\pi}{m_i}\frac{\partial}{\partial\theta_i}$ due to \cref{prop: no critical at boundary}; around the corner, we can mix them to take $\bfv$ of the form $\frac{2\pi}{m_1^2+m_2^2}
\left(m_1\frac{\partial}{\partial\theta_1}+m_2\frac{\partial}{\partial\theta_2}
\right)$. Consequently, its flow preserves the boundary and $M\setminus\partial M$.
\item \label{it: tanget to T for tyep I} On each boundary $T_i$ of a type I critical neighborhood $W_i$, $\bfv$ is tangent to $T_i$ and satisfies $\alpha(\bfv)=1$.
\item \label{it: tangent to T for type II} On each boundary $\fkT_l$ of a type II critical neighborhood $\fkW_l$, $\bfv$ is tangent to $\fkT_l$ and satisfies $\alpha(\bfv)=1$.
\item \label{it: tangent to Q} On $Q''$ away from the type I and type II critical neighborhoods, $\bfv$ is tangent to $Q''$ and satisfies $\alpha(\bfv)=1$.
\item \label{it: horizontal near critical boundary} On the annular regions $W_i\cap f_i^{-1}(\overline{\Delta(0;5c)})\setminus W_i''$ and $\fkW_i\cap f_i^{-1}(\overline{\Delta(0;5c)})\setminus\fkW_i''$, the vector field may be chosen so that
$$df_i(\bfv)=1\Longleftrightarrow d\Re f_i(\bfv)=1,d\Im f_i(\bfv)=0.$$
\end{enumerate}
\begin{proof}
Conditions (a)--(g) can be satisfied simultaneously as follows. For each critical neighborhood $W$ (denoting either $W_i$ or $\fkW_l$), the bad loci are contained in $W'''$. Away from these loci, choose local fields tangent to $Q$ and the prescribed boundary pieces $T$, with $\alpha(v)=1$. On $W'\setminus W''$, the local primitive and its restrictions to the $T$ (or $Q''$ outside of critical neighborhoods) along which tangency is required are submersions. Hence these fields
can also be chosen to satisfy $df_i(v)=1$ there, consistently with $\alpha(v)=1$ since $d\Re f_i=\alpha$.

Multiply these fields by smooth cutoff functions that vanish on $W'''$ and equal $1$ on $W\setminus W''$. The fields therefore remain unchanged on $W'\setminus W''$, preserving (g). Patching these fields with the normalized fields outside the
critical neighborhoods by a partition of unity preserves the prescribed tangency and normalization conditions. The resulting field is smooth on $M$, tangent to every component of $Q''$ and $T$, vanishes on the sets $W'''$, and satisfies $\alpha(v)=1$ outside the union of the sets $W''$.
%We also control $d\Im f_i$ in (g) so that $|f_i|$ cannot decrease too rapidly along the flow, ensuring that the flow does not enter $W_i''$, and hence $W_i'''$, within a short time.
\end{proof}

Let $\widehat{\bfv}$ denote the lift of $\bfv$ to $\widehat M$. Equip $M$ with a Riemannian metric and $\widehat M$ with its lifted metric. Since $M$ is compact, $\bfv$ has bounded length, and hence so does $\widehat{\bfv}$. Moreover, the lifted metric on $\widehat M$ is complete, and $\widehat{\bfv}$ is tangent to the lifted boundary
$$\partial\widehat M=\bigcup \widehat P_i^j.$$
Therefore $\widehat{\bfv}$ is complete, and its flow, denoted by $\Theta_s$, is defined for all $s\in\R$.

Denote $\gamma_s(a)\coloneqq \Theta_s(\gamma(a))$. Since
\begin{equation}\label{eq: compact subset asso to gamma and c}
K_\gamma\coloneqq \{\gamma_s(a): 0\le s\le c, 0\le a\le 1\}
\end{equation}
is compact, it intersects only finitely many lifts of the type I and type II critical neighborhoods. We relabel those which intersect $K_\gamma$ as
$$\widehat{W}_1,\cdots,\widehat{W}_{n_{\mathrm{I}}};\quad
\widehat{\fkW}_1,\cdots,\widehat{\fkW}_{n_{\mathrm{II}}}.$$
Denote their unions by $\widehat{W}$ and $\widehat{\fkW}$, respectively. We use the same convention for $\widehat{W}_l'$, $\widehat{W}_l''$, $\widehat{W}_l''',\widehat R_i,\widehat T_i$, et cetera.

Put $\widehat Y_{\le t}\coloneqq \widehat Y\cap\{u\le t\}$. For the flow $\Theta_s$, we have the following.
\begin{proposition}\label{prop: of flow}
\begin{enumerate}[label=(\roman*)]
\item A trajectory cannot enter a critical neighborhood of type I, say $\widehat{W}_i'$, through $\widehat{T}_i'$.

\item A trajectory cannot enter a critical neighborhood of type II, say $\widehat{\fkW}_l'$, through $\widehat{\fkT}_l'$.

\item For every $s\ge0$ and $x\in\widehat Y_{\le t}\setminus\bigl(\widehat W\cup\widehat{\fkW}\bigr)$ such that $\Theta_{[0,s]}(x)\cap\bigl(\widehat W\cup\widehat{\fkW}\bigr)=\varnothing$, we have
$$\Theta_s(x)\in\widehat Y$$
and
\begin{equation}\label{eq: flow increase}
u(\Theta_s(x))=u(x)+s.
\end{equation}
In particular, the same \cref{eq: flow increase} holds for $x$ belonging $\widehat{T}_i$ and $\widehat{\fkT}_l$.
\end{enumerate}
\end{proposition}

\begin{proof}
Statements \textup{(i)} and \textup{(ii)} follow from \ref{it: tanget to T for tyep I} and \ref{it: tangent to T for type II}, respectively, since the vector field is tangent to $\widehat T_i$ and $\widehat{\fkT}_l$. For \textup{(iii)}, by \ref{it: parallel to P}, the flow preserves the lift of the polar boundary $p^{-1}(P)$ and hence also its complement. Thus a trajectory starting in $\widehat Y$ cannot cross the polar boundary. Moreover, $Q'$ is contained in the type I critical neighborhoods, while all transverse intersection points of distinct components of $Q''$ which require special treatment are contained in mixed type I or type II critical neighborhoods. Away from these critical neighborhoods, $Q''$ is smooth and, by \ref{it: tangent to Q}, the vector field is tangent to $Q''$. Hence a trajectory satisfying the hypothesis in \textup{(iii)} cannot cross the lift of $Q$ either. Therefore, $\Theta_s(x)\in\widehat Y$.

For \cref{eq: flow increase}, recall that on the cyclic cover $2\pi\,du=p^*d\arg\Phi$. Outside $\widehat W\cup\widehat{\fkW}$, we have $d\arg\Phi(\bfv)=2\pi$, and therefore
$\frac{d}{ds}u(\Theta_s(x))=du(\widehat{\bfv})=1$. Integrating from $0$ to $s$ gives $u(\Theta_s(x))=u(x)+s$. The same computation applies along $\widehat T_i$ and $\widehat{\fkT}_l$ by \ref{it: tanget to T for tyep I} and \ref{it: tangent to T for type II}.
\end{proof}

The same short-time control applies to critical neighborhoods of both types. Suppose first that a trajectory enters a type I critical neighborhood through $\widehat R_i'$. Let $s_0$ be the first time at which it meets $\widehat R_i'$, and put $y=\Theta_{s_0}(x)$. Thus $|f_i\circ p(y)|=3c$. For $0\le s\le c$, it follows that
$$c<2c\le |f_i\circ p(\Theta_s(y))|\le4c<5c.$$
Thus the trajectory remains in the region $c<|f_i|<5c$ for $0\le s\le c$, and in particular cannot enter $\widehat W_i''\subset\{|f_i|\le c\}$ during this time; thus we have $f_i\circ p(\Theta_s(y))=f_i\circ p(y)+s$. Exactly the same argument applies if a trajectory enters a type II critical neighborhood through $\widehat{\fkR}_l'$.

Moreover, on each connected lift $\widehat W_i$ or $\widehat{\fkW}_i$, there is a real constant such that $\Re(f_i\circ p)$ or $\Re (f_l\circ p)$ differs from $u$ by that constant. Hence in either case, $u(\Theta_s(y))=u(y)+s$ for $0\le s\le c$. In particular, for $y\in\widehat R_i'$ or $y\in\widehat{\fkR}_l'$,
\begin{equation}\label{eq: u increase for point on Ri}
u(\Theta_s(y))=u(y)+s,\quad 0\le s\le c.
\end{equation}
%The choice $5c<\delta$ ensures that these short trajectory segments remain in the regions where the corresponding local primitives and vector-field constructions are valid. By choosing the critical neighborhoods and these collars sufficiently small, they are also disjoint from the other critical neighborhoods, so no interaction with another critical neighborhood occurs during this short-time flow.

We now consider the finitely many lifted critical neighborhoods intersecting $K_\gamma$. Since the same argument applies to both types, we write $\widehat W_i$ for either type, with $f_i$, $\widehat R_i'$, and $\widehat T_i'$ denoting the corresponding objects. After a small perturbation, assume that $\gamma$ is transverse to their boundaries.

Let $\gamma|_{[a,b]}$ be a component lying in $\widehat W_i'$, with $x=\gamma(a)$ and $y=\gamma(b)$. Then
$x,y\in\widehat R_i'\cup\widehat T_i'$. By \cref{eq: flow increase,eq: u increase for point on Ri}, both points may be flowed for time $c$ with desired $u$-increase:
$$u(\Theta_c(x))=u(x)+c,\quad u(\Theta_c(y))=u(y)+c.$$
The flowed subpath $\Theta_c(\gamma|_{[a,b]})$ is homotopic relative to its endpoints to
$$\sigma_x^{-1}*\gamma|_{[a,b]}*\sigma_y,$$
where $\sigma_x(s)=\Theta_s(x)$ and $\sigma_y(s)=\Theta_s(y)$. This path lies in the corresponding critical neighborhood. Since on $\widehat W_i$ we have $\Re(f_i\circ p)=u+c_i$ for some $c_i\in\R$, \cref{prop: local sup level 1 connected} allows us to replace it, relative to its endpoints, by a path on which
$$r_0+c\le u\le t.$$

We first apply $\Theta_c$ to the whole path $\gamma$, and then perform this replacement on every portion passing through a critical neighborhood. Outside these neighborhoods, \cref{prop: of flow} gives an increase of $u$ by $c$. Hence the resulting path is homotopic to $\gamma$ and has minimum $u$-value at least $r_0+c$. Meanwhile, outside the critical neighborhoods, the choice of vector field prevents the flow from crossing $\widehat P$ or $\widehat Q$, while inside the critical neighborhoods the replacements are carried out entirely in the corresponding neighborhoods with $Q$ deleted, namely $\widehat W_i^*$, $\widehat{\fkW}_l^*$, or $\widehat W_i^\dagger$. Hence the resulting path remains in $\widehat Y$.

To deform $\gamma$ into $L_t$ without ever pushing any point above the level $t$, we allow the flow time to vary with the initial point. For $x\in\widehat Y_{\le t}$, define
$\tau_c(x)=\min\{c,t-u(x)\}$. Thus $0\le\tau_c(x)\le c$, and $\tau_c(x)=0$ for $x\in L_t$. On the portions of the path where the flow is controlled, we deform by
$$H(s,a)=\Theta_{s\tau_c(\gamma(a))}(\gamma(a)),
\quad 0\le s,a\le1,$$
so that $u(H(1,a))=u(\gamma(a))+\tau_c(\gamma(a))= \min\{u(\gamma(a))+c,t\}$. After the local replacements inside the critical neighborhoods, denote the resulting path by $\gamma_c$, we obtain
$$\min_{a\in[0,1]}u(\gamma_c(a)) \ge \min\{r_0+c,t\}, \qquad u(\gamma_c(a))\le t
\quad\text{for all }a\in[0,1].$$

Consequently, for a fixed constant $c>0$, after finitely many iterations the whole path lies in $L_t$, and the endpoints are fixed throughout. During this process, at successive iterations we consider the compact sets $K_{\gamma},K_{\gamma_c},K_{\gamma_{2c}},\ldots$; although the lifted critical neighborhoods intersecting these compact sets may vary from one iteration to the next, only finitely many occur at each step, and we can apply the replacement in \cref{prop: local sup level 1 connected} for finitely many times.

Now recall that $u:\widehat M\to\R$ is proper. For a regular value $t$, the full level $u^{-1}(t)\subset \widehat M$ is therefore compact. Moreover, by the local description of the lift of $P\cup Q$, the level $L_t=u^{-1}(t)\cap\widehat Y$ has finite topological type: it is obtained from the compact level $u^{-1}(t)$ by deleting the lifted boundary $\widehat P$ and $\widehat Q$. Although the lifted boundary may have infinitely many components, these form a locally finite family, so only finitely many meet $L_t$. The compact pair $\left(\overline L_t,\,\overline L_t\cap(\widehat P\cup\widehat Q)\right)$ admits a finite compatible triangulation. Consequently, $L_t=\overline L_t\setminus(\widehat P\cup\widehat Q)$ has the homotopy type of a finite CW complex. In particular, $\pi_1(L_t)$ is finitely generated.
We can now prove \cref{thm: quasi-proj ker not fg implies fact}.

\begin{proof}[Proof of \cref{thm: quasi-proj ker not fg implies fact}]
Suppose first that $\dim\psi(X)\ge2$. By the preceding argument, for every regular value $t$ the inclusion $L_t\hookrightarrow\widehat Y$
induces a surjection on fundamental groups. Since $\pi_1(L_t)$ is finitely generated, so is $\pi_1(\widehat Y)\cong \ker \chi$, contradicting the assumption. Hence $\dim\psi(X)\le1$. Since $\chi\ne0$, the image of $\psi$ cannot be a point, and therefore $\dim\psi(X)=1$.

Let $Z=\overline{\psi(X)}^{\mathrm{Zar}}\subset A$ and let $C_0\to Z$ be the normalization. Then $\psi$ factors through $C_0$. Taking the Stein factorization after passing to proper compactifications and then restricting back to $X$, we obtain a surjective algebraic morphism $f:X\to C$ onto a smooth quasi-projective curve $C$ with connected general fiber, together with a homomorphism
$\chi_C:\pi_1^{\orb}(C)\to\Z$
such that $\chi=\chi^\orb_C\circ f_*^{\orb}$.

We claim that the orbifold curve associated with $f$ is hyperbolic. Let $g$ be the genus of a smooth projective completion of $C$, let $r$ be the number of points removed, and let $m_i$ be the multiplicities of the multiple fibers (see \cref{sec: Riemann orbifolds and fibration}). Since $\chi\ne0$ factors through $f_*^{\orb}$, the orbifold fundamental group admits a nontrivial homomorphism to $\Z$, and in particular the underlying curve has nontrivial first homology. Thus $2g+r-1\ge1$. Hence $\chi^{\orb}(C)\le0$.

Suppose that $\chi^{\orb}(C)=0$, then
\cref{eq: euler char is zero} implies that
$\pi_1^\orb(C)$ is free abelian, and hence $\ker\chi_C$ is finitely generated. Moreover, the orbifold homotopy exact sequence \cref{eq: ses of rimann fibration} shows that $\ker(f_*^\orb)$ is a quotient of the finitely generated fundamental group of a general fiber. Thus $\ker\chi$, being an extension of $\ker\chi_C$ by $\ker(f_*^\orb)$, is finitely generated, a contradiction. Therefore
$\chi^{\orb}(C)<0$, and the orbifold curve associated with $f$ is hyperbolic.
\end{proof}

We now prove \textup{(ii)} of \cref{thm: Lefschetz alternative rational case} for each regular level set under the assumption that \textup{(i)} does not hold. If $\dim\psi(X)=1$, the preceding argument gives precisely the factorization in \textup{(i)}. Hence we must have $\dim\psi(X)\ge2$.

We now pass from the surface case to $X$. Throughout this argument, we keep the form $\alpha=\Re(\psi^*\omega)$ on $X$ fixed. Let $p:\widehat X\to X$ be the covering associated with $\ker\chi$, and choose a primitive $u:\widehat X\to\R$ satisfying $du=p^*\alpha$. For $t\in\R$, put $L_t=u^{-1}(t)\subset \widehat X$.

Let $x_1,x_2\in L_t$. By \cref{lem: surface through prescribed points}, we may choose a smooth connected surface $Y\subset X$ containing $p(x_1)$ and $p(x_2)$ such that $\pi_1(Y)\xrightarrow{\cong}\pi_1(X)$ and $\dim\psi(Y)=2$. Then $\widehat Y=p^{-1}(Y)$ is connected, contains $x_1,x_2$, and its inclusion induces an isomorphism $\pi_1(\widehat Y)\xrightarrow{\cong}\pi_1(\widehat X)$. The primitive used on $\widehat Y$ is the restriction $u|_{\widehat Y}$.

The conclusion established for regular levels of $u|_{\widehat Y}$ above shows that $L_t\cap\widehat Y$ is path-connected. Thus $x_1$ and $x_2$ can be joined inside $L_t$, proving that $L_t$ is path-connected. Moreover, choosing a base point in $L_t\cap\widehat Y$, the composite $\pi_1(L_t\cap\widehat Y)\twoheadrightarrow \pi_1(\widehat Y)\xrightarrow{\cong}\pi_1(\widehat X)$ is surjective. Since this composite factors through $\pi_1(L_t)$, the inclusion induces a surjection $$\pi_1(L_t)\twoheadrightarrow\pi_1(\widehat X).$$ Therefore $(\widehat X,L_t)$ is $1$-connected for every $t\in\R$ regular.

\subsection{From regular level set to any level set}\label{subsec: From regular level set to any level set}
To pass from regular level sets to arbitrary level sets, the missing ingredient is an analogue of \cite[Lemma 9]{Simpson93Lef}, and we include its proof here for completeness. This is precisely the remaining case of \cref{prop: local sup level 1 connected}, namely $\eta_1=\eta_2=0$.

\begin{proposition}\label{prop: local 1-connectedness central real level I}
For a critical neighborhood $W_i$ of type I, set
$$
E_i^*=W_i^*\cap\{\Re f_i=0\}.
$$
Then the pair $(W_i^*,E_i^*)$ is $1$-connected.
\end{proposition}

\begin{proof}
Let
$I=\{is:-\delta<s<\delta\}\subset\Delta(0,\delta)$
be the imaginary axis, and write
$I=I^-\sqcup\{0\}\sqcup I^+$, where $I^\pm$ means the positive/negative imaginary segment.
Then $E_i^*=f_i^{-1}(I)\cap W_i^*$. By \cref{prop: local 1-connectedness for type I}, the fibers
$f_i^{-1}(z)\cap W_i$ are connected for every $z\ne0$. Since
$$
f_i:W_i^{**}\to\Delta^*
$$
is a locally trivial fibration, it follows that
$f_i^{-1}(I^+)\cap W_i^*$ and
$f_i^{-1}(I^-)\cap W_i^*$ are connected.

We first show that
$f_i^{-1}(0)\cap W_i^*$
is contained in the closure of
$f_i^{-1}(I^+)\cap W_i^*$. Let
$x\in f_i^{-1}(0)\cap W_i^*$. If $x$ lies in the interior of $W_i$, then the germ of $f_i$ at $x$ is nonconstant. Hence there is a holomorphic disc through $x$ on which $f_i$ is nonconstant, and the open mapping theorem gives points mapping to $I^+$ arbitrarily close to $x$. If $x\in T_i$, the same conclusion follows from the local triviality of
$f_i|_{T_i}:T_i\to\Delta(0,\delta)$.
Since the central fiber does not meet $R_i$, these are all possibilities. Thus
$$
f_i^{-1}(0)\cap W_i^*
\subset
\overline{f_i^{-1}(I^+)\cap W_i^*}.
$$
It follows that
$$
\bigl(f_i^{-1}(I^+)\cup f_i^{-1}(0)\bigr)\cap W_i^*
$$
is connected.

It remains to connect the negative part to the positive part. Choose
$y\in f_i^{-1}(I^-)\cap T_i$. Such a point exists because
$f_i|_{T_i}$ is a locally trivial fibration with nonempty fibers. Lifting a segment in $I$ joining $f_i(y)$ to a point of $I^+$ gives a path in
$T_i\cap E_i^*$
joining the negative part to the positive part. Hence $E_i^*$ is connected.

Finally, choose $z\in I^+$. Then $f_i^{-1}(z)\cap W_i\subset E_i^*$, and by \cref{prop: local 1-connectedness for type I} the composite
$\pi_1(f_i^{-1}(z)\cap W_i)
\to\pi_1(E_i^*)\to\pi_1(W_i^*)$
is surjective. Therefore
$\pi_1(E_i^*)\twoheadrightarrow\pi_1(W_i^*)$. Thus $(W_i^*,E_i^*)$ is $1$-connected.
\end{proof}

The same argument gives the corresponding statements for type II and mixed type I critical neighborhoods.

\begin{proposition}\label{prop: local 1-connectedness central real level II}
For a critical neighborhood $\fkW_l$ of type II, put
$$
\fkE_l^*=\fkW_l^*\cap\{\Re f_l=0\}.
$$
Then the pair $(\fkW_l^*,\fkE_l^*)$ is $1$-connected.
\end{proposition}

\begin{proposition}\label{prop: local 1-connectedness central real level mixed}
For a mixed type I critical neighborhood $W_i$, put
$$
E_i^\dagger=W_i^\dagger\cap\{\Re f_i=0\}.
$$
Then the pair $(W_i^\dagger,E_i^\dagger)$ is $1$-connected.
\end{proposition}

Now let $t$ be a singular value of $u$, and still put $L_t\coloneqq u^{-1}(t)\cap\widehat Y.$
We show that every path in $\widehat Y$ with endpoints in $L_t$ can be deformed, relative to its endpoints, into a path contained in $L_t$. 

In contrast to the regular case, a path may intersect $L_t$ infinitely many times, so we cannot simply perturb it to be transverse to $L_t$. Instead, we argue as follows: after a homotopy relative to the endpoints inside finitely many holomorphic coordinate balls in $\widehat Y$, we may assume that $\gamma$ is piecewise real-analytic. Since $u$ is locally the real part of a holomorphic function, on each such piece $u\circ\gamma-t$ either vanishes identically or has only finitely many zeros. Thus $\gamma^{-1}(L_t)$ is a finite union of points and closed intervals. Leaving the portions contained in $L_t$ fixed, we can therefore divide the remaining portions into finitely many subpaths whose endpoints lie in $L_t$ and whose interiors lie entirely either below or above $L_t$. We treat these subpaths separately. We can thus argue on each subpath as done in \cref{subsec: Construction of a flow on the infinite cyclic covering}, with \cref{prop: local sup level 1 connected,prop: local 1-connectedness central real level I,prop: local 1-connectedness central real level II,prop: local 1-connectedness central real level mixed} used whenever needed.

%\begin{remark}
%A similar situation is considered in \cite[\S 5.1.1]{Pybook}, where, after a subtle replacement, the path has only finitely many portions outside the singular level set. It is not clear to us whether that replacement path is homotopic to the original one.
%\end{remark}

%The first approach is to remove the singular obstruction before applying the argument as in the regular-level case. Since $u$ is proper on $\widehat{M}$, the singular level $L_t$ is relatively compact and meets only finitely many critical neighborhoods. In every Type I or Type II neighborhood $\widehat W_i$ satisfying $L_t\cap\widehat W_i=\{\Re f_i=0\}\cap\widehat W_i$, we first replace the portions of $\gamma$ passing through the central fiber $f_i^{-1}(0)$ by nearby paths contained in $\widehat W_i\setminus f_i^{-1}(0)$. The local replacement is possible thanks to \cref{prop: local 1-connectedness central real level I,prop: local 1-connectedness central real level II,prop: local 1-connectedness central real level mixed}. After performing these replacements in all relevant critical neighborhoods, we obtain a path homotopic to $\gamma$, relative to its endpoints, which avoids the central fibers. Along this modified path there is no longer any singular obstruction, and the same perturbation and flow argument as in the regular-level case applies. It therefore deforms the modified path into $L_t$, which is homotopic to the original path.

We can also argue using a varying flow time on the whole path, with a sign depending on which side of $L_t$ the point lies. This will also be used in \cref{subsec: The flow construction in general case} for the general case. For $x\in\widehat Y$, put
$$\tau_c(x)=\min\{c,|t-u(x)|\}$$
and define the signed flow time
$$\sigma_c(x)=\operatorname{sgn}(t-u(x))\,\tau_c(x),$$
where $\operatorname{sgn}(0)=0$. Thus $\sigma_c$ is continuous (near $L_t$ it is in fact smooth), it is positive on $\{u<t\}$ and negative on $\{u>t\}$, and it vanishes on $L_t$. We deform the whole path by
$$H(s,a)=\Theta_{s\sigma_c(\gamma(a))}(\gamma(a)),
\quad 0\le s,a\le1.$$
As in \cref{subsec: Construction of a flow on the infinite cyclic covering}, the choice of vector field ensures that the flow preserves $\widehat Y$; we use this without further comment below. Since $\sigma_c$ is continuous and the flow is smooth, this is a continuous homotopy fixing every point of $L_t$. On portions whose trajectories stay outside critical neighborhoods, \cref{prop: of flow}(iii) and its negative-flow analogue give
$$u(\Theta_s(x))=u(x)+s,$$
and therefore $|t-u(H(1,a))|=\max\{|t-u(\gamma(a))|-c,0\}$. 

Since
$$K_{\gamma}\coloneqq \left\{\Theta_{s\sigma_c(\gamma(a))}(\gamma(a)): 0\le s\le1,\ 0\le a\le1 \right\}$$
is compact, it intersects only finitely many lifted critical neighborhoods. Inside these neighborhoods the same flow is applied, but the preceding equality for $|t-u|$ need not hold. After the local replacements inside the critical neighborhoods thanks to \cref{prop: local sup level 1 connected,prop: local 1-connectedness central real level I,prop: local 1-connectedness central real level II,prop: local 1-connectedness central real level mixed}, we have
$$\min\{r_0+c,t\}\le u(\gamma_c(a))\le\max\{t,R_0-c\},\quad\text{for all }0\le a\le1,$$
where $r_0=\min_{a\in[0,1]}u(\gamma(a))\le t$ and $R_0=\max_{a\in[0,1]}u(\gamma(a))\ge t$. We denote the resulting path by $\gamma_c$.

%In conclusion, the argument in \cref{subsec: Construction of a flow on the infinite cyclic covering} works on both sides of $L_t$, with $|t-u|$ replacing the one-sided quantity $t-u$. At each stage, points outside the critical neighborhoods move toward $L_t$ by the prescribed amount, while the corresponding critical portions are replaced relative to their boundary points using local 1-connectedness, with prescribed behavior on $u$. Points which have already reached $L_t$ remain fixed. 

We repeat the procedure for the resulting path, using at each stage the corresponding compact set $K_\gamma, K_{\gamma_c}, K_{\gamma_{2c}},\ldots$. Although the lifted critical neighborhoods meeting these compact sets may vary from one step to the next, only finitely many occur at each step. As in the regular-level argument, after finitely many iterations the whole path lies in $L_t$. Since $\sigma_c=0$ on $L_t$, the endpoints remain fixed throughout. Hence every path in $\widehat Y$ with endpoints in $L_t$ can be deformed, relative to its endpoints, into $L_t$. 

The passage from $\widehat Y$ to $\widehat X$ follows the argument at the end of \cref{subsec: Construction of a flow on the infinite cyclic covering}. Although formulated there for regular level sets, this argument applies verbatim to arbitrary level sets using the surface statement established above. 
%In particular, no transversality of $Y$ to $L_t$ is required.

\subsection{The compact and quasi-projective cases: a comparison}\label{subsec: Theorem of Napier--Ramachandran}
This section discusses the differences between the compact and noncompact settings for \cref{thm: NR on Z covering} and \cref{thm: quasi-proj ker not fg implies fact}, with particular emphasis on the differences in their proofs.

For the theorem of Napier--Ramachandran (\cref{thm: NR on Z covering}), it was firstly stated and proved in \cite{NapierRamachandran01}, where the proof builds on the potential-theoretic methods developed in \cite{NapierRamachandran95}. Roughly speaking, they use Evans--Selberg potential theory to construct pluriharmonic functions on complete (noncompact) K\"ahler manifolds under suitable various conditions. By constructing pluriharmonic functions with prescribed asymptotic behavior along different ends (for example, \cite[Theorem 2.6]{NapierRamachandran95}) and studying the holomorphic foliations induced by their differentials, they use a cup-product lemma to obtain proper holomorphic mappings onto Riemann surfaces with connected fibers. This process is further developed in \cite{NapierRamachandran08} using the notion of ``filtered ends''. In particular, they proved:
\begin{theorem}[\cite{NapierRamachandran08}, Theorem 4.3]\label{thm: 3 filtered ends}
A connected complete K\"ahler manifold of bounded geometry with at least three filtered ends admits a proper holomorphic mapping onto a Riemann surface. 
\end{theorem}
For the infinite cyclic covering associated with an $\Z$-valued character whose kernel is not finitely generated, the covering has exactly two ends but infinitely many filtered ends, which is a restatement of \cref{prop: non fg implies level set disconnect}. Thus, three filtered-ends theorem with an argument of descending to the $\Z$-covering provides another perspective on \cref{thm: NR on Z covering}.

In the previous sections (\cref{subsec: On semi-Albanese map determined by a character,subsec: Choice of logarithmic holomorphic one-form,subsec: Reduction to surfaces,subsec: Construction of two types of critical neighborhoods,subsec: Construction of a flow on the infinite cyclic covering}), we provide a proof of this theorem in the quasi-projective setting. The reader may also use the argument in \cref{subsec: On semi-Albanese map determined by a character,subsec: Choice of logarithmic holomorphic one-form,subsec: Reduction to surfaces,subsec: Construction of two types of critical neighborhoods,subsec: Construction of a flow on the infinite cyclic covering} to recover a proof of \cref{thm: NR on Z covering} in this compact case.

\begin{remark}\label{rmk: other ways to prove NR Z cover thm}
The \cref{thm: NR on Z covering} can be proved using other methods, at least for the compact K\"ahler case.
\begin{enumerate}
\item One may use the min--max argument of Delzant--Py to prove \cref{thm: NR on Z covering}; see \cite[Chapter~5]{Pybook}. In fact, \cite[Theorem~5.1]{Pybook}, applied to the primitive function $u$ on $\widehat M$, gives the conclusion. As pointed out there, this argument applies only to the rational case.
\item One may also modify the argument of \cite[Lemma 3]{Zhou} to give another proof of \cref{thm: NR on Z covering}, working instead on the universal covering and using \cref{prop: non fg implies level set disconnect}. This approach applies to arbitrary real characters. Readers interested in this direction may similarly modify the argument of \cref{sec: Lefschetz alternative rational case} to work on the universal covering. We indicate the main idea. In this setting, it is unnecessary to keep track of the relative homotopy class of the path. If a level set in the universal covering is disconnected, then the Mayer--Vietoris sequence for reduced homology implies that either the corresponding sublevel set or the superlevel set is disconnected, or both. Suppose, for instance, that the superlevel set is disconnected. One then only needs to choose a path joining two of its components of the suplevel set and push the portion of the path lying in the sublevel set into the superlevel set, which produces a contradiction. See \cite[pp.~769--773]{Zhou}.
\end{enumerate}
\end{remark}

Here, we point a difference between the compact K\"ahler case and the (noncompact) smooth quasi-projective case: how to explain conditions \ref{it: condition A} and \ref{it: condition B} in \cref{subsec: Construction of two types of critical neighborhoods}. 

In the compact K\"ahler case, consider the singular holomorphic foliation $\calf$ defined by the closed holomorphic one-form $\xi$ whose real part represents $\chi:\pi_1(X)\twoheadrightarrow\Z$. Condition \ref{it: condition A} produces a compact regular leaf. The compact-leaf theorem (for example, see \cite{NapierRamachandran97,ABCKT,DelzantGromov05}) then implies that all singular leaves are compact, and the corresponding leaf space gives a holomorphic fibration onto a compact Riemann surface with connected fibers. Here singular leaves are understood as the equivalence classes generated by connected fibers of local holomorphic primitives
of $\xi$. This compact-leaf argument does not extend directly to arbitrary smooth quasi-projective varieties, as the following example shows.

%For the compact K\"ahler case, if we consider the holomorphic foliation $\mathcal{F}$ defined by the (closed) holomorphic one-form $\xi$ whose real part corresponds to $\chi:\pi_1(X)\twoheadrightarrow \Z$, then condition \ref{it: condition A} implies that we have a compact leaf of $\calf$, due to a well-known result in K\"ahler geometry , the existence of \textit{one} compact leaf tells us all leaves are compact. Using equivalent relation defined by these leaves, we get the modulo space which is a (compact) Riemann surface; moreover, Stein factorization gives as a holomorphic fibration with connected fibers. Indeed, if we consider the set consisting of compact leaves, it is open due to topological reasons; while for being closed, we need the K\"ahler condition and also the compactness. This is precisely the point that does not carry over directly to the noncompact quasi-projective setting, where compact leaves may escape toward infinity.
\begin{example}
Let $X=\bigl(\mathbb P^1\times\mathbb C^*\bigr)\setminus\{(\infty,1)\}$
and let $\xi=\frac{1}{2\pi i}\frac{dz_2}{z_2}$.
Then $\operatorname{Re}\xi$ represents an epimorphism $\pi_1(X)\to\Z $, and the leaves of the foliation defined by $\xi$ are the fibers of the projection $X\to\mathbb C^*$. For $z\ne1$, the leaf $\mathbb P^1\times\{z\}$ is compact, whereas the leaf over $z=1$ is $\bigl(\mathbb P^1\setminus\{\infty\}\bigr)\times\{1\}\cong\mathbb C$. Thus the union of the compact leaves is not closed in $X$. %Geometrically, the compact leaves degenerate to a noncompact one because a point of the fiber escapes through the boundary of the quasi-projective manifold.
\end{example}

With the function $\psi$ defined in \cref{subsec: On semi-Albanese map determined by a character}, for compact K\"ahler case, we can use the above foliation argument to get the following implication, when the set of (type I) critical neighborhoods are not empty:
\begin{itemize}
    \item Condition \ref{it: condition A} $\implies$ $\dim \psi(X)=1$ $\Longleftrightarrow$ condition (i) \cref{thm: Lefschetz alternative rational case}.
    \item Condition \ref{it: condition B} $\impliedby$ $\dim \psi(X)\ge 2$ $\implies$ condition (ii) \cref{thm: Lefschetz alternative rational case}.
\end{itemize}
While for noncompact quasi-projective manifolds, the same implications are from using \cref{lem: good choice of holomorphic one-form}(b) and \cref{lem: good choice of holomorphic one-form real case}(b) instead.

In one sentence, for compact K\"ahler manifolds the alternative may be detected either through the existence of a compact leaf, as in \cite{NapierRamachandran95,NapierRamachandran97,DelzantGromov05,Zhou,Liu23}, or through whether $\dim\psi(X)=1$ or $\dim\psi(X)\ge2$, as in \cite{Simpson93Lef}; for noncompact quasi-projective manifolds, however, the latter dimension criterion is the only one available to us.

\begin{remark}
Note that in \cite{Simpson93Lef}, different components of $\Sigma$ with the same $\psi$-values are put in a common critical neighborhood; while the critical neighborhood (of type I) in \cref{subsec: Construction of two types of critical neighborhoods} is chosen around the individual connected components $\Sigma_i$.

In the compact K\"ahler case, assume $\Sigma\ne\varnothing$ and interpret conditions \ref{it: condition A} and \ref{it: condition B} using Simpson's critical neighborhoods, then:
\begin{itemize}
    \item Condition \ref{it: condition A}
    $\Longleftrightarrow$ $\dim\psi(X)=1$
    $\Longleftrightarrow$ \textup{(i)} of
    \cref{thm: Lefschetz alternative rational case}.
    \item Condition \ref{it: condition B}
    $\Longleftrightarrow$ $\dim\psi(X)\ge2$
    $\implies$ \textup{(ii)} of
    \cref{thm: Lefschetz alternative rational case}.
\end{itemize}
\end{remark}

The alternatives \textup{(i)} and \textup{(ii)} may overlap, see \cref{rmk: overlap for Lefschetz alternative}. In particular, \textup{(ii)} does not imply $\dim\psi(X)\ge2$ in general. If we further require in (i), the orbicurve $C$ is hyperbolic, then (i) and (ii) are exclusive, and 
\begin{itemize}
    \item Condition \textup{(i)} in \cref{thm: Simpson's Lefschetz} with $C$ hyperbolic. $\Longleftrightarrow$ $\ker(\chi)$ is not finitely generated.
    \item Condition \textup{(ii)} in \cref{thm: Simpson's Lefschetz} $\Longleftrightarrow$ $\ker(\chi)$ is finitely generated.
\end{itemize}
For rational characters in the smooth quasi-projective setting, the same equivalences hold with \textup{(ii)} understood for the prescribed logarithmic representative $\alpha$; see \cref{subsec: Lefschetz alternative revisited via the BNS invariant}. The rationality assumption and the restriction on the representative are essential here.

\section{From the rational case to the general case}\label{sec: Lefschetz alternative general case}

In this section, we modify the argument in \cref{sec: Lefschetz alternative rational case} to prove the Lefschetz alternative for an arbitrary nonzero real character in the smooth quasi-projective setting. We will explicitly identify the points at which essential modifications are required and explain how the relevant properties differ from those in the rational case.

Let $X$ be a connected smooth complex quasi-projective variety, let $G=\pi_1(X)$, and let $\chi:G\to\R$ be a nonzero character. Denote by $r$ the $\Z$-rank of $\chi(G)\subset\R$. When $r=1$, this is the rational case discussed in \cref{sec: Lefschetz alternative rational case}; when $r\ge2$, we are in the irrational case, which is the main subject of present section.

The main difference between the two settings is that the primitive function $u:\widehat M\to\R$ is not proper in the irrational case. Thus the finite-generation argument used in the rational case no longer applies: the level sets need not have the homotopy type of a finite CW complex. We explain this further in \cref{subsec: Lefschetz alternative revisited via the BNS invariant}.

Choose an identification $\chi(G)\cong\Z^r$. By abuse of notation, we also denote by $\chi$ the induced epimorphism $G\to\Z^r$ through which the original character factors:
\begin{equation}\label{diag: chi through free abelian}
\begin{tikzcd}
G\ar[r,"\chi",two heads]\ar[dr,"\chi"'] & \Z^r\ar[d,hook]\\
& \R.
\end{tikzcd}
\end{equation}

We now prove the following.

\begin{theorem}\label{thm: Lefschetz alternative general case}
Let $X$ be a smooth connected complex quasi-projective variety and let $\chi:\pi_1(X)\to\R$ be a nonzero group homomorphism. Let $p:\widehat X\to X$ be the regular covering corresponding to $\ker(\chi)$, so that $\pi_1(\widehat X)=\ker(\chi)$. Then one of the following holds:
\begin{enumerate}[label=(\roman*)]
    \item There exists a surjective algebraic morphism with connected general fiber $f:X\to C$
    onto a smooth quasi-projective orbicurve such that $\chi$ factors through $\pi_1^\orb(C)$.
    \item There exists a smooth closed $1$-form $\alpha$ on $X$ representing $\chi$ such that, for a primitive function $
    u:\widehat X\to\R$ satisfying 
    $du=p^*\alpha$, the pair $(\widehat X,u^{-1}(t))$ is $1$-connected for every $t\in\R$.
\end{enumerate}
\end{theorem}

\subsection{Modifications before lifting}\label{subsec: Modifications before lifting}

For \cref{subsec: On semi-Albanese map determined by a character}, we still take the semi-abelian variety $A_0$, the quasi-Albanese morphism $a:X\to A_0$ associated to $X$, and write $\lambda_0\in H^1(A_0;\R)$ such that $a^*\lambda_0=\chi$. Let $B$ be the sum of connected positive dimensional algebraic subgroups where $\lambda_0$ vanishes, and $\psi:X\to A\coloneqq A_0/B$ be the quasi-Albanese map associated to $\chi$. The \cref{lem: nonvanishing on postive subgp after quotient} still holds by our construction, while \cref{lem: rank of torus le 1} is certainly not true. In fact, due to the same argument, we now have $\ell\le r$. To avoid any possible confusion, we restate the lemmas from \cref{sec: Lefschetz alternative rational case} that remain valid in the present setting, at the risk of some repetition.

\begin{lemma}\label{lem: nonvanishing on postive subgp after quotient real case}
For every positive dimensional connected algebraic subgroup $H\subset A$, the restriction of $\lambda$ to $H$ is not identically zero.
\end{lemma}

Now we make the modification of \cref{subsec: Choice of logarithmic holomorphic one-form}. There is still a natural compactification of $A$ as follows. Recall that we have a short exact sequence of algebraic groups
$$1\to T=(\C^*)^\ell\to A\xrightarrow{\pi}A_{\ab}\to1.$$
Consider the action of $T$ on $\BBP^\ell$ defined by
$$(t_1,\cdots,t_\ell)\cdot[z_0:z_1:\cdots:z_\ell]=[z_0:t_1z_1:\cdots:t_\ell z_\ell].$$
The open orbit satisfying $z_0z_1\cdots z_\ell\ne0$ is isomorphic to $T\subset\BBP^\ell$. Treating $A$ as a principal $T$-bundle over $A_{\ab}$, we take the corresponding $\BBP^\ell$-bundle over $A_{\ab}$, which is smooth projective, denoted as $\bar A$. Then $\bar A\setminus A$ consists of finitely many divisors:
$$D_A=\bigcup_{0\le i\le\ell}D_i\coloneqq\bigcup_{0\le i\le\ell}\{z_i=0\}.$$
Unlike in the rank-one case, these boundary divisors may intersect; nevertheless, their union is a simple normal crossings divisor. Under the same notation as in \cref{subsec: Choice of logarithmic holomorphic one-form}, we have
$$\bar A=\BBP_{A_{\ab}}\bigl(\mathcal O_{A_{\ab}}\oplus L_1\oplus\cdots\oplus L_\ell\bigr),$$
where $L_i\in\Pic^0(A_{\ab})$ is the line bundle corresponding to the $i$-th coordinate character of $T=(\C^*)^\ell$.

Take $\gamma_k$ be the $k$-th generator of $\pi_1(T)\cong \Z^\ell$ and take $a_k=\lambda(\gamma_k)\in \R$. For each $1\le k\le\ell$, it is standard that there exists a global logarithmic $1$-form $\eta_k$ on $A$, extending logarithmically to $\bar A$, such that
$$\Res_{D_j}\eta_k=
\begin{cases}
1,&j=k,\\
-1,&j=0,\\
0,&j\ne0,k.
\end{cases}$$
%For example, this follows from the existence of holomorphic flat connections on the line bundles $L_k\in\Pic^0(A_{\ab})$.
Locally, after trivializing $L_k$, the form $\eta_k$ is given by $d\log(w_k/w_0)$ up to a holomorphic $1$-form pulled back from the base $A_\ab$.

Then
$$\lambda-\left[\Re \sum_{k=1}^\ell\frac{a_k}{2\pi\sqrt{-1}}\eta_k\right]$$
has trivial image in $\Gr_2^W H^1(A;\R)$, and we can find $\lambda_\ab\in H^1(A_\ab;\R)$ such that 
$$\lambda=\sum_{k=1}^\ell\left[\Re \frac{a_k}{2\pi\sqrt{-1}}\eta_k\right]+\pi^*\lambda_\ab.$$
Let $\omega_\ab$ be a holomorphic one-form on $A_\ab$ such that $[\Re \omega_\ab]=\lambda_\ab$ and 
\begin{equation}\label{eq: log one form irrational}
\omega=\sum_{i=1}^\ell \frac{a_i}{2\pi \sqrt{-1}}\eta_i +\pi^*\omega_\ab
\end{equation}
is the logarithmic one-form on $A$ with $[\Re \omega]=\lambda$, and $\omega$ has purely imaginary residues at each divisor in $D_A$. Due to \cref{lem: nonvanishing on postive subgp after quotient real case}, each $a_k\ne 0$ and $a_1+\cdots+a_\ell\ne 0$. In conclusion, the residue of $\omega$ along each irreducible component of $D_A$ is purely imaginary and nonzero. Also, due to the same reason as in \cref{lem: good choice of holomorphic one-form}(b), we have
\begin{lemma}\label{lem: good choice of holomorphic one-form real case}
\begin{enumerate}
\item There is a closed algebraic one-form $\omega$ on $A$ such that
$$[\Re\omega]=\lambda,$$
and $\omega$ extends logarithmically to the standard projective compactification of $A$ with nonzero purely imaginary residues at each boundary divisor. 
\item For every nonconstant algebraic map $u:C\to A$ from a smooth algebraic curve,
$$u^*\omega\not\equiv0.$$
\end{enumerate}
\end{lemma}

For \cref{subsec: Reduction to surfaces}, there is no essential change, except that $\Phi_A$ introduced there is no longer well-defined. Also, from now on we assume that $\dim\psi(X)\ge2$. With the notation introduced above, let $\bar\psi:\bar Y\to\bar A$
be the extension to smooth compactifications with SNC divisors and put $P=(\bar\psi^{-1}(D_A))_{\red}$. Write $P=\bigcup_\nu P_\nu$ for its irreducible components and set $\xi=\bar\psi^*\omega$ on $Y^\circ=\bar Y\setminus P$.

Near a point of $P_\nu$, choose local coordinates such that
$$P=\{z_1\cdots z_k=0\},\quad k=1,2.$$
After choosing a local trivialization of the $\BBP^\ell$-bundle $\bar A\to A_{\ab}$, let $t_i=w_i/w_0$, $1\le i\le\ell$, denote the corresponding torus characters. Then locally
$$\bar\psi^*t_i=v_i z_1^{n_{i1}}\cdots z_k^{n_{ik}},$$
where $v_i$ is a nonvanishing holomorphic function and $n_{ij}\in\Z$. Hence
$$\bar\psi^*\frac{dt_i}{t_i}=\frac{dv_i}{v_i}+\sum_{j=1}^k n_{ij}\frac{dz_j}{z_j}.$$
Absorbing all holomorphic terms into a locally defined holomorphic function $g$, we obtain
\begin{equation}\label{eq: xi real case}
\xi=\bar\psi^*\omega=\sum_{j=1}^k\frac{c_j}{2\pi\sqrt{-1}}\frac{dz_j}{z_j}+dg,
\end{equation}
where $c_j=\sum_{i=1}^\ell a_i n_{ij}\in\R$. Writing $z_j=r_je^{\sqrt{-1}\theta_j}$, we get
\begin{equation}\label{eq: alpha real case}
\alpha\coloneqq\Re\xi=\sum_{j=1}^k\frac{c_j}{2\pi}d\theta_j+\Re dg.
\end{equation}

We claim that $c_j\ne0$ for every $1\le j\le k$. For a fixed $j$, the vector $(n_{1j},\ldots,n_{\ell j})\in\Z^\ell$ is nonzero. It therefore determines a nonconstant one-parameter subgroup
$$\rho_j:\C^*\to T=(\C^*)^\ell,\quad s\mapsto(s^{n_{1j}},\ldots,s^{n_{\ell j}}).$$
If $\gamma_i$ denotes the $i$-th standard generator of $\pi_1(T)$, then $(\rho_j)_*(1)=\sum_{i=1}^\ell n_{ij}\gamma_i$, and hence $\lambda((\rho_j)_*(1))=\sum_{i=1}^\ell a_i n_{ij}=c_j$. If $c_j=0$, then $\rho_j^*\lambda=0$, and therefore $\lambda$ vanishes identically on the positive-dimensional connected algebraic subgroup $\im\rho_j\subset A$. This contradicts \cref{lem: nonvanishing on postive subgp after quotient real case}. Thus $c_j\ne0$.

\begin{lemma}
In the local form \cref{eq: xi real case} of $\xi$, we have $c_j\ne0$ for each $1\le j\le k$, where $k=1$ or $2$.
\end{lemma}

Now take the same real-oriented blowup, so that the angular coordinates remain well-defined when the radial coordinates vanish. Near a point where
$$P=\{z_1\cdots z_k=0\},\quad z_j=r_je^{\sqrt{-1}\theta_j},
\ k=1,2,$$
we have \cref{eq: alpha real case} for $\alpha$. Consider the boundary corresponding to $P_j$, given by $r_j=0$. Since $g$ is holomorphic across $P_j$, its restriction to $r_j=0$ is independent of the angular variable $\theta_j$. Hence $\frac{\partial}{\partial\theta_j}\Re g=0$ along $r_j=0$. Therefore
$$
\alpha\left(\frac{\partial}{\partial\theta_j}\right)=\frac{c_j}{2\pi}\ne0.
$$
It follows that $\alpha$ does not vanish on the boundary of the real-oriented blowup. Since $\partial M$ is compact, for any fixed smooth Riemannian metric on $M$, there exists a constant $\varepsilon_0>0$ such that $\norm{\alpha_x}\ge\varepsilon_0$ for every $x\in\partial M$. By continuity, after shrinking to a collar neighborhood $N$ of $\partial M$, we may assume that
\begin{equation}\label{eq: uniform lower bound}
\norm{\alpha_x}\ge\frac{\varepsilon_0}{2} \quad\text{for every $x\in N$.}
\end{equation}
This is the uniform bound we will use to construct a flow later.

The argument in \cref{subsec: Construction of two types of critical neighborhoods} goes through verbatim in the general setting, since we are working downstairs and nothing related to the rationality of $\chi$ is used there. In conclusion, we still have finitely many Type I and Type II critical neighborhoods, with the mixed case classified as Type I.

\subsection{The flow construction in the general case}\label{subsec: The flow construction in general case}

This is the part where we need a major revision of \cref{subsec: Construction of a flow on the infinite cyclic covering}.

Since $\alpha$ is closed on $M$, it defines
$$\chi_M:\pi_1(M)\to\R,\quad \chi_M([\gamma])=\int_\gamma\alpha.$$
Under the identification $\pi_1(M)\cong\pi_1(Y^\circ)$, the inclusion $j:Y\hookrightarrow Y^\circ$ induces a surjection on $\pi_1$, whose kernel is normally generated by small positively-oriented loops around $Q$. Since $\alpha$ extends across $Q$, these loops have zero $\alpha$-period, and hence $\chi=\chi_M\circ j_*$ and $\im\chi_M=\im\chi\cong\Z^r$.
Let $p:\widehat M\to M$ be the covering corresponding to $\ker\chi_M$. Then $p^*\alpha=du$ for some $u:\widehat M\to\R$, and $\widehat Y:=p^{-1}(Y)$ is precisely the covering of $Y$ corresponding to $\ker\chi$.

We first consider the lifts of the two types of critical neighborhoods. For Type I, it can be argued as before that
$$p^{-1}(W_i)=\bigsqcup_{\gamma\in\Z^r}\widehat W_{i,\gamma}=\bigsqcup_{\gamma\in\Z^r}\gamma\widehat W_{i,0}.$$
Indeed, on $W_i$ we have a local primitive $f_i$ such that $\xi=df_i$, and hence $\alpha=d\Re f_i$. Thus, for any loop $\gamma$ in $W_i$, we have
$$\chi([\gamma])=\int_\gamma\alpha=\int_\gamma d\Re f_i=0,$$
i.e., $\pi_1(W_i)\subset\ker(\chi)$, which gives the conclusion above. For a Type II domain, the situation is the same, and we have
$$p^{-1}(\fkW_l)=\bigsqcup_{\gamma\in\Z^r}\widehat{\fkW}_{l,\gamma}=\bigsqcup_{\gamma\in\Z^r}\gamma\widehat{\fkW}_{l,0}.$$
The mixed case is still classified as Type I and has the same lifting property.

Now consider the vector field $\bfv$ with properties as in \cref{subsec: Construction of a flow on the infinite cyclic covering}. Precisely, take
$$W_i'''\Subset W_i''\Subset W_i\cap f_i^{-1}(\overline{\Delta(0;c)})\Subset W_i'\Subset W_i$$
as shown there, and similarly
$$\fkW_l'''\Subset\fkW_l''\Subset \fkW_l\cap f_l^{-1}(\overline{\Delta(0;c)})\Subset \fkW_l'\Subset\fkW_l.$$
Choose a smooth vector field $\bfv$ on $M$ such that
\begin{enumerate}[label=(\alph*)$^*$]
    \item $\bfv$ vanishes in the neighborhood $W_i'''$ for each $i$ and $\fkW_l'''$ for each $l$.
    \item Outside Type I critical neighborhoods $\bigsqcup_{i=1}^M W_i''$ and Type II critical neighborhoods $\bigsqcup_{l=1}^N \fkW_l''$, we have
    $$\alpha(\bfv)=1.$$
    \item $\bfv$ is tangent to every boundary face $P$ of the real-oriented blowup $M$ and satisfies $\alpha(\bfv)=1$. Such a choice is possible by the uniform lower bound in \cref{eq: uniform lower bound}; in fact, locally near a smooth point of $P_j=\{r_j=0\}$, we can simply take $\bfv=\frac{2\pi}{c_j}\frac{\partial}{\partial\theta_j}$. Near a corner, the corresponding choices can be interpolated within the tangent bundle of the boundary while preserving $\alpha(\bfv)=1$.
    \item On each boundary $T_i$ of a Type I critical neighborhood $W_i$, $\bfv$ is tangent to $T_i$ and satisfies $\alpha(\bfv)=1$.
    \item On each boundary $\fkT_l$ of a Type II critical neighborhood $\fkW_l$, $\bfv$ is tangent to $\fkT_l$ and satisfies $\alpha(\bfv)=1$.
    \item On $Q''$ away from Type I and Type II critical neighborhoods, $\bfv$ is tangent to $Q''$ with $\alpha(\bfv)=1$.
    \item On the annular regions
    $W_i\cap f_i^{-1}(\overline{\Delta(0;5c)})\setminus W_i''$ and $\fkW_i\cap f_i^{-1}(\overline{\Delta(0;5c)})\setminus\fkW_i''$, the vector field may be chosen so that
    $$d\Re f_i(\bfv)=1,d\Im f_i(\bfv)=0.$$
\end{enumerate}
Let $\widehat{\bfv}$ be the lift of $\bfv$ to $\widehat M$. Since $M$ is compact and $\bfv$ is tangent to $\partial M$, the vector field $\bfv$ is complete, and hence so is $\widehat{\bfv}$. Denote by $\Theta_s$ the flow generated by $\widehat{\bfv}$. For $t\in\R$, set
$$L_t:=\{x\in\widehat Y:u(x)=t\}.$$

We now use the second method sketched in \cref{subsec: From regular level set to any level set}. For $c>0$, define the continuous function on $\widehat M$:
$$\sigma_c(x)=\begin{cases}
\min\{c,t-u(x)\}, & \quad \text{if }u(x)<t,\\
0, & \quad \text{if }u(x)=t,\\
-\min\{c,u(x)-t\}, & \quad \text{if }u(x)>t.
\end{cases}$$
For a path $\gamma:[0,1]\to\widehat Y$, define
$$H(s,a):=\Theta_{s\sigma_c(\gamma(a))}(\gamma(a)),\quad 0\le s,a\le1.$$
Since $\sigma_c$ is continuous and $\widehat{\bfv}$ is smooth, $H$ is a continuous homotopy. Moreover, if $\gamma(a)\in L_t$, then $\sigma_c(\gamma(a))=0$, and hence $H(s,a)=\gamma(a)$ for all $s\in[0,1]$. Thus the points of $\gamma$ already lying in $L_t$ remain fixed throughout the deformation.

Denote the compact subset
$$K_\gamma=\{\Theta_{s\sigma_c(\gamma(a))}(\gamma(a)):0\le s\le1,\ 0\le a\le1\}.$$
Although the inverse images of the Type I and Type II critical neighborhoods may have infinitely many connected components in $\widehat M$, these components form a locally finite collection, since the deck group $\Z^r$ acts properly discontinuously on $\widehat M$. Consequently, the compact set $K_\gamma$ meets only finitely many of them, relabeled as
$$\widehat W_1,\ldots,\widehat W_{n_I};\widehat{\fkW}_1,\ldots,\widehat{\fkW}_{n_{\mathrm{II}}}.$$
Similarly, the inverse images of the polar boundary $P$ and the nonpolar boundary $Q$ may have infinitely many lifts. In the irrational case, these faces may moreover occur at arbitrarily close $u$-levels, since the subgroup $\chi(G)\subset\R$ need not be discrete. Nevertheless, the lifted boundaries form a locally finite family by the proper discontinuity of the deck action. Consequently, the compact subset $K_\gamma$ meets only finitely many of them.

Away from these critical neighborhoods, we have from \cref{prop: of flow} that
$$du(\widehat{\bfv})=p^*\alpha(\widehat{\bfv})=\alpha(\bfv)=1.$$
Hence, as long as the corresponding flow segment remains outside the critical neighborhoods $W'\cup\fkW'$,
$$u(\Theta_s(x))=u(x)+s.$$
The flow moves points toward the level $L_t$; more precisely, outside of the critical neighborhoods, we have
\begin{equation}\label{eq: upper and lower bounds}
\min\{r_0+c,t\}\le u(\gamma_c(a))\le\max\{t,R_0-c\},\quad\text{for all }0\le a\le1.
\end{equation}
where
$$r_0=\min_{a\in[0,1]}u(\gamma(a))\le t,\quad R_0=\max_{a\in[0,1]}u(\gamma(a))\ge t.$$

Whenever the deformation encounters one of the finitely many lifted Type I or Type II critical neighborhoods, we use \cref{prop: local sup level 1 connected} to replace the corresponding portion of the path by one satisfying simultaneously \cref{eq: upper and lower bounds}. If both endpoints lie in $L_t$, corresponding to the real level $\{\Re f_i=\eta\}$, we use \cref{prop: local 1-connectedness for type I,prop: local 1-connectedness for type II,prop: local 1-connectedness for mixed} when $\eta\ne0$, and
\cref{prop: local 1-connectedness central real level I,prop: local 1-connectedness central real level II,prop: local 1-connectedness central real level mixed} when $\eta=0$. These replacements are performed in the same homotopy class relative to the endpoints. Moreover, the flow
preserves $\widehat Y$ by $(c)^*$ and $(f)^*$, and the local replacements take place in the corresponding neighborhoods with $Q$ deleted, say $W_i^*,\fkW_l^*$ and $W_i^\dagger$.

After applying the homotopy $H$ and making the local replacements described above, denote the resulting path by $\gamma_c$. At each iteration, the corresponding compact subsets $K_\gamma,K_{\gamma_c},K_{\gamma_{2c}},\cdots$ meets only finitely many lifted Type I and Type II critical neighborhoods. Whenever such a neighborhood is encountered, we make the same local replacement as above. The collection of critical neighborhoods involved may vary from one iteration to the next, but it is finite at every step. Moreover, by $(c)^*$ and $(f)^*$, together with the construction of the local replacements, the resulting path remains entirely in $\widehat Y$ throughout the deformation.

Thus, unless a point has already reached $L_t$, its distance in the $u$-coordinate from $t$ decreases by at least $c$ at each step. Since $c>0$ is fixed, after finitely many iterations we obtain a path homotopic to $\gamma$, relative to its endpoints, which lies entirely in $L_t$. This proves the pair $(\widehat Y, L_t)$ is 1-connected.

By the same argument as at the end of \cref{subsec: Construction of a flow on the infinite cyclic covering}, the $1$-connectedness of $(\widehat Y,L_t\cap\widehat Y)$ transfers to $(\widehat X,L_t)$. Consequently, if $\dim\psi(X)\ge2$, then alternative~\textup{(ii)} of \cref{thm: Lefschetz alternative general case} holds. If $\dim\psi(X)=1$, we argue as in the proof of \cref{thm: Lefschetz alternative rational case} to obtain alternative~\textup{(i)}.

\begin{remark}
Notice that, when $r\ge2$, the primitive function $u:\widehat M\to\R$ is not proper.
%Indeed, the subgroup $\chi(G)\subset\R$ has $\Z$-rank at least two and hence is dense in $\R$. For a fixed $x\in\widehat M$, we can choose distinct deck transformations $\gamma_n\in\Z^r$ such that $\chi(\gamma_n)\to0$. Since $$u(\gamma_nx)=u(x)+\chi(\gamma_n),$$ the points $\gamma_nx$ eventually lie in the inverse image of a fixed compact interval in $\R$, whereas they have no convergent subsequence by proper discontinuity of the deck action. Hence $u$ is not proper.
In particular, unlike in the rational case, the level sets of $u$, even regular ones, need not have finite CW type or finitely generated fundamental group. We postpone the interpretation of this phenomenon in terms of the BNS invariant to \cref{subsec: Lefschetz alternative revisited via the BNS invariant}.
\end{remark}

\section{Applications on quasi-projective groups}\label{sec: applications}

\subsection{Results on BNS invariants}\label{subsec: Applications on BNS invariants}
The application of Simpson's Lefschetz theorem to the BNS invariant of K\"ahler groups was given by Delzant in \cite{Delzant10MA}. Using our generalization of the Lefschetz theorem to smooth quasi-projective varieties, we obtain an analogous description of the BNS invariant for quasi-projective groups. We refer the reader to \cite{BNS87,BSunpublished,Strebel12LN} for the general results on BNS invariants used in this section. Throughout this section, we assume that $H^1(G;\R)\ne0$; otherwise, the character sphere $\mathrm{S}(G)$ is empty and there is nothing to consider.

Let $G$ be a finitely generated group and $\chi\colon G\to\R$ be a nontrivial group homomorphism. Denote $G_{\chi}=\{g\in G\colon \chi(g)\geqslant0\}$. For any finite generating set $\mathcal G$, denote the associated Cayley graph of $G$ by $C\coloneqq\mathrm{Cayley}(G;\mathcal G)$. We define $C_{\chi}$ to be the full subgraph of $C$ spanned by $G_\chi$. We consider whether $C_{\chi}$ is connected and record those $\chi$ for which this happens. Notice that this property does not depend on the generating set $\mathcal G$ (see \cite[Proposition 11.1]{Pybook}) and is invariant under multiplication of $\chi$ by a positive real constant. Denote by $[\chi]$ the equivalence class of $\chi$ under multiplication by positive real constants, and write $\mathrm{S}(G)=(\mathrm{Hom}(G;\R)\setminus\{0\})/\R^{+}$ for the collection of all such equivalence classes. For a subgroup $H\le G$, define
$$\rms(G,H)\coloneqq\{[\chi]\in\rms(G):\chi|_H=0\}.$$

\begin{definition}
With the notation above, for a finitely generated group $G$, we define its BNS invariant by
$$\Sigma^1(G)=\{[\chi]\in\mathrm{S}(G)\colon C_{\chi}\text{ is connected}\}.$$
\end{definition}

Note that $\Sigma^1(G)$ is an open subset of $\rms(G)$; see \cite[Theorem~A]{BNS87}. Moreover, we give a list of results we need on the BNS invariant:
\begin{proposition}\label{prop: of bns}
Let $G$ be a finitely generated group.
\begin{enumerate}
    \item \textbf{(\cite[Theorem H]{BNS87})} Let $\pi:G\to Q$ be a group epimorphism. Then
    $$\pi^*\left(\Sigma^1(Q)^\complement\right)\subset\Sigma^1(G)^\complement.$$
    \item \textbf{(\cite[Theorem B1]{BNS87})} Let $\pi:G\to A$ be a group epimorphism with $A$ abelian. Then $\ker(\pi)$ is finitely generated if and only if
    $$\rms(G,\ker(\pi))\subset\Sigma^1(G).$$
    \item If $G$ is abelian, then
    $$\Sigma^1(G)=\rms(G).$$
\end{enumerate}
\end{proposition}

In \cite{LiuLiu25}, building on a result of Suciu \cite{Suciu21}, we described the BNS invariant of a K\"ahler group in terms of the \textit{integral} tropicalization of its rank-one jump ideal. More precisely, if $X$ is a compact K\"ahler manifold with $G=\pi_1(X)$, then
$$\Sigma^1(G)^\complement=\rms\bigl(\Trop_{\Z}(\mathcal{J}^1(G;\Z))\bigr).$$
Here $\rms$ denotes the projection from $H^1(G;\R)\setminus\{0\}$ onto the character sphere $\rms(G)$, and the tropicalization over $\Z$ is defined in \cite[Proposition 3.9]{LiuLiu25} using the trivial, mod-$p$, and $p$-adic valuations. The analogous description in terms of orbifold fibrations also holds for quasi-projective groups, by the following result.

\begin{theorem}\label{thm: BNS quasiprojective}
Let $X$ be a connected smooth complex quasi-projective variety with $G=\pi_1(X)$. Then
$$
\Sigma^1(G)^\complement=
\rms\left(\bigcup_f\im\left(f^*:H^1(C,\mathbf m;\R)\to H^1(X;\R)\right)\right),
$$
where the union runs over all hyperbolic orbifold fibrations $f:X\to(C,\mathbf m)$. Moreover, this is a finite union.
\end{theorem}

\begin{proof}
Set $T$ to be the right hand side of above equation, and we prove the two inclusions.

First, let $f:X\to(C,\mathbf m)$ be a hyperbolic orbifold fibration and put $Q=\pi_1^{\orb}(C,\mathbf m)$. The induced homomorphism $f_*:G\twoheadrightarrow Q$ is surjective. For a hyperbolic Riemann orbifold, we have $\Sigma^1(Q)^\complement=\rms(Q)$ by \cref{prop: bns empty for hyperbolic RS}. By \cref{prop: of bns}(a), we obtain $$f^*\rms(Q)=f^*\Sigma^1(Q)^\complement\subseteq\Sigma^1(G)^\complement.$$
Taking the union over all hyperbolic orbifold fibrations gives $T\subseteq\Sigma^1(G)^\complement$; this union is finite due to \cref{thm: finiteness of hyperbolic fibration when noncompact}.

Conversely, let $[\chi]\in\Sigma^1(G)^\complement$.
Suppose that $\dim\psi(X)\ge2$, and take the logarithmic representative $\alpha=\Re\omega$ constructed in \cref{sec: Lefschetz alternative general case}. Let $K=\mathrm{Bl}_{P\cup Q}^\mathrm{ro}(\bar X)$ be the full real-oriented blowup of an SNC compactification of $X$. Since the logarithmic form $\omega$ extends holomorphically across $Q$ and has purely imaginary residues along $P$, its real part $\alpha$ extends smoothly to $K$. Moreover, $K$ is a compact manifold with corners whose interior is $X$. Since the inclusion $X\hookrightarrow K$ is a homotopy equivalence, we have $\pi_1(K)\cong G$.

By our construction, every level set of primitive function $u$ in the $\chi$-cover of $X$ is connected and has $\pi_1$ surjects onto $\ker\chi$. Its inverse image in the maximal abelian cover of $K$ is therefore connected and has $\pi_1$ surjects onto the fundamental group of that cover. Consequently, the inclusion of the (unique) unbounded superlevel component of $u$ on the maximal abelian cover of $K$ also induces a surjection on fundamental groups. Applying \cite[Theorem G and \S5]{BNS87} on the compact model $K$ gives $[\chi]\in\Sigma^1(G)$, a contradiction. Hence $\dim\psi(X)=1$, and the Stein factorization  yields an algebraic orbifold fibration in \textup{(i)}.

We claim that $(C,\mathbf m)$ is hyperbolic. Indeed, the same argument as in the proof of \cref{thm: quasi-proj ker not fg implies fact} applies. If $(C,\mathbf m)$ were not hyperbolic, then, since $\chi_C\ne0$, the only possibilities are $2g+r=2$. Moreover, by (\ref{eq: char of orbifold}), there are no multiple fibers and $Q$ is free abelian. By \cref{prop: of bns}(b)-(c), $\ker(\chi_C)$ is finitely generated. Since the kernel of the orbifold fibration is also finitely generated, the exact sequence
$$1\to N\to\ker(\chi)\to\ker(\chi_C)\to1$$
would imply that $\ker(\chi)$ is finitely generated. By \cref{prop: of bns}(b), this gives $\rms(G,\ker\chi)\subseteq\Sigma^1(G)$, contradicting $[\chi]\in\Sigma^1(G)^\complement$. Therefore $(C,\mathbf m)$ is hyperbolic. Hence $\Sigma^1(G)^\complement\subseteq T$, which proves the desired equality.
\end{proof}

\begin{corollary}\label{Cor: antipodal symmetric}
The BNS invariant and its complement for the fundamental group of a connected smooth complex quasi-projective variety are both antipodally symmetric.
\end{corollary}

\subsection{Application on quasi-projective groups}\label{subsec: Application on quasi-projective groups}

Recall that a quasi-projective group is a group that can be realized as the fundamental group of a smooth complex quasi-projective variety. Not every finitely presented group is quasi-projective; see, for example, \cite{Morgan,KapovichMillson} for known obstructions.

Combining the antipodal symmetry in \cref{Cor: antipodal symmetric} with the argument of \cite[Theorem~11]{Liu26CR}, we obtain the following obstruction to a group being quasi-projective.

\begin{theorem}\label{thm: HNN}
Any group admitting a properly ascending or descending HNN decomposition with finitely generated base group is not a quasi-projective group.
\end{theorem}

Indeed, the canonical character of such an HNN extension is integral, and \cref{thm: Lefschetz alternative rational case} alone is sufficient to prove this. Hence, the proof of \cite[Theorem 11]{Liu26CR} applies verbatim, and we omit the details. The corresponding result for K\"ahler groups is due to Napier--Ramachandran \cite{NapierRamachandran08}. Their proof, however, does not directly extend to the quasi-projective setting.

Next, we consider Thompson--Stein groups. Let $n_1,\ldots,n_k\ge2$ and
$$\ell\in \Z\left[\frac{1}{n_1\cdots n_k}\right],\quad\ell>0.$$
We denote by $\langle n_1,\ldots,n_k\rangle=\left\{
n_1^{a_1}\cdots n_k^{a_k}: a_1,\ldots,a_k\in\Z\right\}$
the multiplicative subgroup of $\R_{>0}^*$ generated by $n_1,\ldots,n_k$. The Thompson--Stein group
$$F\coloneqq F\left(
\ell,
\Z\left[\frac{1}{n_1\cdots n_k}\right],
\langle n_1,\ldots,n_k\rangle
\right)$$
is the group of all orientation-preserving piecewise-linear homeomorphisms $f:[0,\ell]\to[0,\ell]$ having finitely many breakpoints, such that all breakpoints and their images belong to
$$
\Z\left[\frac{1}{n_1\cdots n_k}\right],
$$
and every slope of $f$ belongs to $\langle n_1,\ldots,n_k\rangle$. This group is known to be finitely presented.%, with an explicit presentation given by Wladis in \cite{Wladis}.

By \cite[Theorem 8.1]{BNS87} (which is true for any $\ell$ as noticed by Fried--Vidussi in \cite[Appendix]{FriedlVidussi16}), we have
\begin{equation}\label{eq: BNS complement of Thompson-Stein gp}
\Sigma^1(F)^\complement=\{[\log \lambda],[\log \rho]\}
\end{equation}
where $\lambda$ and $\rho$ are given by derivatives at the left and right end-points of $[0, \ell]$ respectively:
$$\lambda: F\to \R^+,\quad \lambda(f)=\frac{df}{dt}(0^+)$$
$$\rho: F\to \R^+,\quad \rho(f)=\frac{df}{dt}(\ell^-)$$

Denote $r=\mathrm{rank}_\Z\langle n_1,\cdots,n_k\rangle\le k$. When $r\ge 2$, we can see that $\log(\lambda)$ and $\log(\rho)$ are both irrational. Indeed, for $f\in F$,
write its initial uniquely as
$$f'(0^+)=n_1^{a_1}\cdots n_k^{a_k}$$
then define characters $\chi^0_i:F\to \Z$ defined as $f\mapsto a_i$, $1\le i\le k$. Then 
$$\log(\lambda)=\sum_{i=1}^k \chi^0_i\cdot \log(n_i)$$
If $\log\lambda$ is rational, choose element $f_i\in F$ such that the initial slope is $n_i$, then we need $\log(n_i)/\log(n_j)\in \Q$ for each $1\le i\ne j\le k$
, or equivalently, $n_i^a=n_j^b$ for some integers $a$ and $b$. If $r\ge 2$, then we know $\log\lambda$ is irrational; similarly, $\log \rho$ is also irrational. Moreover, these two points cannot be antipodal. Indeed, one can explicitly construct an element $f\in F$ whose initial slope is $1$ while its terminal slope is nontrivial (see \cite[\S 8]{BNS87} or \cite{BieriStrebel16book}). Hence $\log\lambda$ and $\log\rho$ cannot be proportional; in particular, these two cannot be antipodal. 

Now, \cref{Cor: antipodal symmetric} yields the following:
\begin{proposition}\label{prop: Thompson not quasi-proj}
Any Thompson--Stein group is not quasi-projective.
\end{proposition}

For the case $k=1$, write
$$F_m\coloneqq F\left(\ell,\Z \left[\frac{1}{m}\right],\langle m\rangle\right).$$
When $m=2$, this recovers the ordinary Thompson group $F$, which was shown not to be K\"ahler by Napier--Ramachandran \cite{NapierRamachandran06}. More generally, Friedl--Vidussi proved in \cite[Theorem~3.3]{FriedlVidussi16} proves that all Thompson--Stein groups are not K\"ahler, using rank-gradient arguments together with the BNS invariants. The same antipodal-symmetry argument also gives another proof of \cite[Theorem~3.3]{FriedlVidussi16}. We also remark that for the case $k=1$, the proof of \cref{prop: Thompson not quasi-proj} only requires the rational Lefschetz alternative, say \cref{thm: Lefschetz alternative rational case}.

\subsection{Lefschetz alternative revisited via the BNS invariant}\label{subsec: Lefschetz alternative revisited via the BNS invariant}

With the above explanation via the BNS invariant, especially \cref{Cor: antipodal symmetric}, we are able to understand the relationship between \cref{thm: Lefschetz alternative general case} and \cref{thm: quasi-proj ker not fg implies fact}.

When $\chi:G\to\R$ has image of rank $1$, up to rescaling, we can assume that $\chi:G\twoheadrightarrow\Z$. To determine whether $\ker\chi$ is finitely generated, by \cref{prop: of bns}(b), it is the same to determine whether
$$\rms(G,\ker\chi)\subset\Sigma^1(G).$$
But now $\rms(G,\ker\chi)=\{[\chi],[-\chi]\}$. Also, by \cref{Cor: antipodal symmetric}, these two points are either both contained in $\Sigma^1(G)$ or both contained in $\Sigma^1(G)^\complement$. Hence, to determine whether $\ker\chi$ is finitely generated, one only needs to check whether $[\chi]\in\Sigma^1(G)$. Now, with the complement as explained in \cref{thm: BNS quasiprojective}, we have the following dichotomy:
\begin{itemize}
    \item $\ker\chi$ is not finitely generated $\Longleftrightarrow$ $[\chi]\in\Sigma^1(G)^\complement$ $\Longleftrightarrow$ condition \textup{(i)} of \cref{thm: Lefschetz alternative rational case} with $C^\orb$ hyperbolic.
    \item $\ker\chi$ is finitely generated $\Longleftrightarrow$ $[\chi]\in\Sigma^1(G)$ $\Longleftrightarrow$ condition \textup{(ii)} of \cref{thm: Lefschetz alternative rational case} for $\alpha$ the real part of a suitable chosen logarithmic one-form as in \cref{subsec: Choice of logarithmic holomorphic one-form}.
\end{itemize}

The above consideration also explains why we cannot expect analogues of \cref{thm: quasi-proj ker not fg implies fact} or \cref{thm: NR on Z covering} to hold for irrational characters. More precisely, we can find an irrational character $\chi:G\to\R$ such that $\ker\chi$ is not finitely generated, while $\chi$ does not factor through any holomorphic fibration, even in the compact K\"ahler case.

The process of finding such an example is straightforward. We first choose a K\"ahler group $G$ for which both $\Sigma^1(G)$ and $\Sigma^1(G)^\complement$ are nonempty. Since $\Sigma^1(G)$ is open, we can choose a rational character $\varphi$ such that $[\varphi]\in\Sigma^1(G)$. On the other hand, $\Sigma^1(G)^\complement$ is a union of positive-dimensional rational subspaces arising from holomorphic fibrations; see \cite{Delzant10MA} for the compact K\"ahler case, while \cref{thm: BNS quasiprojective} gives the corresponding statement in the more general smooth quasi-projective setting. Hence $\Sigma^1(G)^\complement$ also contains a nonzero rational character, and we choose a rational character $\psi$ such that $[\psi]\in\Sigma^1(G)^\complement$. Now consider
$$\chi=\varphi+a\psi$$
with $a$ an irrational number sufficiently small so that $[\chi]\in\Sigma^1(G)$, due to openness of $\Sigma^1(G)$. By our construction, the equality $\chi(g)=0$ implies $\varphi(g)=\psi(g)=0$. Hence
$$\ker(\chi)=\ker(\varphi)\cap\ker(\psi)\subset\ker(\psi).$$
In particular,
$$[\psi]\in\rms(G,\ker\chi)\cap\Sigma^1(G)^\complement,$$
and therefore $\ker\chi$ is not finitely generated by \cref{prop: of bns}(b). On the other hand, $[\chi]\in\Sigma^1(G)$, so by \cref{thm: BNS quasiprojective}, $\chi$ cannot factor through any hyperbolic holomorphic fibration. It cannot factor through a nonhyperbolic one either, with the same argument as in the proof of \cref{thm: BNS quasiprojective}

This example also reminds us that in the proof of \cref{thm: Lefschetz alternative general case}, the level set constructed there, say $L_t\subset\widehat Y$, need not have the homotopy type of a finite CW complex. Indeed, the example above gives an irrational character $\chi$ satisfying condition \textup{(ii)} of \cref{thm: Lefschetz alternative general case} for which $\ker(\chi)$ is not finitely generated. Since
$$\pi_1(L_t)\twoheadrightarrow\pi_1(\widehat Y)=\ker(\chi),$$
it follows that $\pi_1(L_t)$ is not finitely generated. In particular, $L_t$ cannot have the homotopy type of a finite CW complex.

\section*{Acknowledgments}
The author would like to thank Chen Zhou for helpful discussions concerning his paper and for pointing out the example appearing in \cref{rmk: accumulating zero sets}. The author would also like to thank Yongqiang Liu for many useful discussions; the collaboration with him provided an important motivation for writing this paper. The author is partially supported by the National Natural Science Foundation of China (Grant No.~12501105).

\bibliographystyle{ssmfalpha}
\bibliography{reference}

\end{document}